\documentclass[a4paper,12pt]{article}
\usepackage{amsmath}
\usepackage[vlined]{algorithm2e}
\usepackage{diagbox}
\usepackage{enumerate}
\usepackage{appendix}
\usepackage{cases}
\usepackage{longtable}
\usepackage{caption}
\usepackage{anyfontsize}
\usepackage[hidelinks]{hyperref}

\title{{\fontsize{20}{60}\selectfont Unifying  a posteriori error analysis of five piecewise quadratic discretisations for the
biharmonic equation}}
\author{Carsten Carstensen\footnote{Department of Mathematics, 
Humboldt-Universit\"{a}t zu Berlin, 10099 Berlin, Germany.
Distinguished Visiting Professor, Department of Mathematics, Indian Institute of 
Technology Bombay, Powai, Mumbai-400076, India.  cc@math.hu-berlin.de}        
\quad\text{and}\quad Benedikt Gr\"a{\ss}le\footnote{Department of Mathematics, 
Humboldt-Universit\"{a}t zu Berlin, 10099 Berlin, Germany.
graesslb@math.hu-berlin.de}\quad\text{and}\quad
Neela Nataraj\footnote{Department of Mathematics, Indian Institute of Technology Bombay, Powai, Mumbai 400076, India. neela@math.iitb.ac.in}
}
\date{}
\SetFuncSty{textsc}
        \SetKwFunction{frun}{Run}
        \SetKwHangingKw{arun}{\frun}
        \SetKwFunction{fset}{Set}
        \SetKwHangingKw{aset}{\fset}
        \SetKwFunction{fselect}{Select}
        \SetKwHangingKw{aselect}{\fselect}
        \SetKwFunction{fcompute}{Compute}
        \SetKwHangingKw{acompute}{\fcompute}
        \SetKwFunction{fsolve}{Solve}
        \SetKwHangingKw{asolve}{\fsolve}
        \SetKwFunction{festimate}{Estimate}
        \SetKwHangingKw{aestimate}{\festimate}
        \SetKwFunction{fmark}{Mark}
        \SetKwHangingKw{amark}{\fmark}
        \SetKwFunction{frefine}{Refine}
        \SetKwHangingKw{arefine}{\frefine}
        \SetKwFunction{compute}{Compute}
        \SetKwFunction{set}{Set}

{\algorithm}%
   {\endalgorithm}%

\makeatletter
\newcommand\Label[1]{&\refstepcounter{equation}(\theequation)\ltx@label{#1}&}
\makeatother
\newcommand{\W}[1]{{\color{black}#1}}
\newcommand{\WW}[1]{{\color{black}#1}}
\newcommand{\X}[1]{{\color{black}}}

\usepackage{mathtools}
\DeclareFontEncoding{LS1}{}{}
\DeclareFontSubstitution{LS1}{stix}{m}{n}
\DeclareSymbolFont{symbols2}{LS1}{stixfrak} {m} {n}
\DeclareMathSymbol{\operp}{\mathbin}{symbols2}{"A8}

\usepackage{amsmath,amsthm,amssymb,enumerate}
\usepackage[T1]{fontenc}
\usepackage{gensymb}
\usepackage[square,sort&compress,comma,numbers]{natbib}
\usepackage[sc]{mathpazo}
\usepackage{multirow} 
\usepackage{newtxtext,newtxmath}
\usepackage{subfig}
\usepackage{graphicx}
\usepackage{epstopdf}
\usepackage{cancel}
\usepackage[margin=3cm]{geometry}
\usepackage{tikz}
\usepackage{caption}
\usetikzlibrary{shapes,calc}
\usepackage{verbatim}
\usepackage{mathrsfs}
\usepackage{accents}
\usepackage[utf8]{inputenc}

\usetikzlibrary{decorations.pathmorphing}
\usetikzlibrary{decorations.pathreplacing}
\usetikzlibrary{positioning}
\usetikzlibrary{shapes}
\usetikzlibrary{arrows}
\usetikzlibrary{patterns}
\usetikzlibrary{fadings}
\usetikzlibrary{plotmarks}
\usetikzlibrary{calc}
\usetikzlibrary{intersections}
\tikzstyle{every picture}+=[font=\footnotesize]
\usepackage{paralist}
\usepackage{bbm}
\usepackage{latexsym}           %
\usepackage{enumerate}
\usepackage{enumitem}

\setlist{noitemsep, topsep=0.8ex, partopsep=0pt%
	, leftmargin=3em}
\setlist[1]{labelindent=\parindent}

\newlist{axioms}{enumerate}{1}
\setlist[axioms]{font=\bfseries}

\newlist{alphenum}{enumerate}{1}
\setlist[alphenum]{label=\textbf{(\alph*)}, leftmargin=4em}

\newlist{alphienum}{enumerate}{1}
\setlist[alphienum]{label=\textit{(\alph*)}}

\newlist{romanenum}{enumerate}{1}
\setlist[romanenum]{label=\textit{(\roman*)}}

\newlist{romaninenum}{enumerate*}{1}
\setlist[romaninenum]{label=\textit{(\roman*)}}
\usepackage[vlined]{algorithm2e}
\SetKwIF{If}{ElseIf}{Else}{if}{}{else if}{else}{endif}
\SetKwFor{For}{for}{}{endfor}
\usepackage[noabbrev, capitalise]{cleveref}
\usepackage{tabu}
\crefname{equation}{\unskip}{\unskip}
\creflabelformat{equation}{#2(#1)#3}
        \SetFuncSty{textsc}
        \SetKwFunction{frun}{Run}
        \SetKwHangingKw{arun}{\frun}
        \SetKwFunction{fset}{Set}
        \SetKwHangingKw{aset}{\fset}
        \SetKwFunction{fselect}{Select}
        \SetKwHangingKw{aselect}{\fselect}
        \SetKwFunction{fcompute}{Compute}
        \SetKwHangingKw{acompute}{\fcompute}
        \SetKwFunction{fsolve}{Solve}
        \SetKwHangingKw{asolve}{\fsolve}
        \SetKwFunction{festimate}{Estimate}
        \SetKwHangingKw{aestimate}{\festimate}
        \SetKwFunction{fmark}{Mark}
        \SetKwHangingKw{amark}{\fmark}
        \SetKwFunction{frefine}{Refine}
        \SetKwHangingKw{arefine}{\frefine}
        \SetKwFunction{compute}{Compute}
        \SetKwFunction{set}{Set}

\usepackage{tikz}
\usetikzlibrary{matrix}
\usepackage{tikz-cd}
\theoremstyle{definition}
\newtheorem{defn}{Definition}[section]

\newtheorem{rem}{Remark}[section]
\newtheorem{example}{\bf Example}[section]

\numberwithin{equation}{section}

\newcommand{\cT}{\mathcal T}

\newcommand{\pw}{\mathrm{pw}}
\newcommand{\nc}{\mathrm{M}}
\newcommand{\jc}{\mathrm{J}}

\newcommand{\qo}{\mathrm{qo}}

\newcommand{\NC}{\text{pw}}
\newcommand{\C}{\mathrm{C}}
\newcommand{\w}{\mathrm{P}}

\newcommand{\ip}{{\mathrm{IP}}}

\newcommand{\half}{\frac{1}{2}}
\newcommand{\trinl}{\ensuremath{|\!|\!|}}
\newcommand{\trinr}{\ensuremath{|\!|\!|}}

\newcommand{\dx}{{\rm\,dx}}

\newcommand{\ds}{{\rm\,ds}}

\newcommand{\dg}{{\rm dG}}

\DeclareMathOperator{\E}{\mathcal{E}}
\newcommand{\M}{\mathrm{M}}
\newcommand{\CR}{\mathrm{CR}}
\newcommand{\N}{\mathbb{N}}

\newcommand{\T}{\mathcal{T}}

\def\mean#1{\left<#1\right>_E}
\def\jump#1{\left[ #1\right]_E}

\newcommand{\V}{\mathcal{V}}

   \usepackage{xparse}
   
   \newcounter{const}
\NewDocumentCommand{\constant}{o}
 {
  \IfValueTF{#1}%
  {C_{#1}}%
  {\refstepcounter{const}%
  C_{\theconst}}%
 }

\usepackage{tikz,pgfplots}
\usetikzlibrary{calc,angles,positioning,intersections,quotes,decorations.markings}
\usepackage{tkz-euclide}
\pgfplotsset{compat=1.5}

\newcommand{\trb}[1]{|\!|\!|#1|\!|\!|}
\renewcommand{\nc}{{\textup{nc}}}
\newcommand{\Vnc}{V_{\nc}}
\newcommand{\hatF}{\widehat F}
\newcommand{\hatFapx}{\widehat F_{\textup{apx}}}
\newcommand{\hatFh}{\hatFapx}
\newcommand{\hatv}{\widehat v}

\newcommand{\id}{\mathrm{id}}
\newcommand{\JIM}{J_h}

\renewcommand{\div}{\mathrm{div}}
\newcommand{\APX}{\mathrm{apx}}
\newcommand{\OSC}{\mathrm{osc}}
\newcommand{\tg}{G}
\newcommand{\tf}{F}
\newcommand{\hatg}{\widehat \vartheta}
\newcommand{\hatx}{\widehat \xi}
\newcommand{\muold}{\|u_h-J_hu_h\|_h}

\theoremstyle{remark}
\theoremstyle{plain}
\newtheorem{theorem}{Theorem}[section]
\newtheorem{lemma}[theorem]{Lemma}
\newtheorem{proposition}[theorem]{Proposition}
\newtheorem{cor}[theorem]{Corollary}

\usepackage{totcount}
\newtotcounter{nconst}
\makeatletter
\newcount\rc@count
\let\rc@clearconstantlist\empty

\newcommand\rc@clearconstant[1]{\global\expandafter\let\csname rc@const@#1\endcsname\undefined}

\newcommand\resetconstants[1]{%
    \def\rc@constname{#1}
    \global\rc@count=1\relax 
    \bgroup 
        \let\\\rc@clearconstant 
        \rc@clearconstantlist
        \global\let\rc@clearconstantlist\empty 
    \egroup
}
\newcommand\const[1]{%
    \@ifundefined{rc@const@#1}{%
        \expandafter\xdef\csname rc@const@#1\endcsname{%
           \noexpand\rc@useconst{\rc@constname}{\the\rc@count}%
        }%
        \g@addto@macro\rc@clearconstantlist{\\{#1}}%
        \global\advance\rc@count1\relax
    }{}%
	\setcounter{nconst}{\the\rc@count-1}
    \csname rc@const@#1\endcsname
}

\newcommand\rc@useconst[2]{\ensuremath{#1_{#2}}}
\makeatother

	\resetconstants{C}

\begin{document}	
\maketitle
\begin{abstract}
An abstract  property \eqref{eqn:H} is the key to a complete a~priori error 
analysis in the (discrete) energy norm
for several nonstandard finite element methods in the recent work [Lowest-order 
equivalent nonstandard finite
element methods for biharmonic plates,  Carstensen and Nataraj, M2AN, 2022].
This paper investigates the impact of \eqref{eqn:H} to the a posteriori error analysis 
and establishes known and novel explicit residual-based a posteriori error estimates.
The abstract framework applies to Morley, two versions of discontinuous Galerkin, $C^0$ interior penalty, as well as weakly over-penalized symmetric interior penalty
schemes for the biharmonic equation with a general source term in $H^{-2}(\Omega)$.
\end{abstract}


\section{Introduction}\label{sectionintroduction}
The concept of a quasi-optimal
smoother and the key assumption \eqref{eqn:H}
from \cite{carstensen_lowest-order_2022} allow for an abstract a posteriori error analysis 
for five lowest-order schemes for the biharmonic problem.
This paper unifies
and completes \cite{DVNJSR2007, HuShi09, BGS10, Georgoulis2011,
BrenGudiSung10, MozoSuli07} and provides
novel reliable and efficient a posteriori error estimators for a right-hand side $F\in H^{-2}(\Omega)$.

\subsection{Overview}%
\label{sub:Overview}
The traditional view on  a posteriori error control is that the well-posedness of the linear problem on the continuous
level directly leads from the error
to residuals and their dual norms.
In the simplest setting of a Hilbert space $(V,a)$ with induced norm $\trb{\bullet}\coloneqq a(\bullet, \bullet)^{1/2}$, the weak solution $u\in V$ is the Riesz representation of a given
source $F\in V^*$:  $u\in V$ solves
\begin{align}\label{eqn:WP}
	a(u,v)&=F(v) \quad \text{ for all }v\in V.
\end{align}
Given any conforming companion $J_hu_h\in V$ to some discrete approximation $u_h\in V_h$, where $V_h\not\subseteq V$  is
\W{typically not a subset of} $V$ and $J_hu_h\in V$ is a postprocessing of $u_h$, the norm of the error $e\coloneqq u-J_hu_h\in V$ is the norm of the residual  $F-
a(J_hu_h,\bullet)\in V^*$: The Riesz isometry between the residual and its Riesz representation $e\in V$ reads
\begin{align} \label{eqn:error}
	\trb{e} = \trb{F - a(J_hu_h, \bullet)}_{*}\coloneqq \sup_{v\in V\setminus\{0\}}\frac{|F(v) - a(J_hu_h,
	v)|}{\trb{v}}.
\end{align}
The a posteriori error control is left with the task of deriving computable upper and lower bounds of the
dual norm $\trb{F-a(J_hu_h, \bullet)}_*$.
The known data are $F\in V^*$ and $J_hu_h\in V$ and the techniques to derive bounds are
very different from those of an a priori error analysis. 

The paradigm change in this paper employs a recent tool %
\eqref{eqn:H} (stated in Section~\ref{sub:Quasi-optimal approximations} below) from the \emph{a~priori error analysis}
\cite{carstensen_lowest-order_2022} to arrive at an \emph{a posteriori error bound}
\begin{align}\label{eqn:intro_error_2}
	\trb{e}^2 \W{\leq C\big(}\|u_h-J_hu_h\|_h^2+Res((1-J_h\W{I})e)+ \text{data approximation error}\W{\big)}
\end{align}
with some operator $\W{I}\in L(V; V_h)$ and a norm
$\|\bullet\|_h$ on $V+V_h$.
The main advantage of the master estimate \eqref{eqn:intro_error_2} over the error identity \eqref{eqn:error} is the known
structure $(1-J_h\W{I})e\in V$ of the
test function.
The a posteriori error analysis based on \eqref{eqn:intro_error_2} then only requires to study the properties of the
operators  $(1-J_h\W{I})\in L(V; V)$.
This allows explicit estimates 
of the error term $Res((1-J_h\W{I})e)$ with universal arguments for generic $u_h\in V_h$ and, most importantly, independent of the discrete
system that defines $u_h\in V_h$.

The application to the biharmonic equation \eqref{eqn:WP} provides novel simultaneous insight in the residuals and
estimators for the piecewise quadratic discrete solution $u_h\in P_{\hspace{-.13em}2}(\T)$ to the  Morley,
two variants of  discontinuous Galerkin (dG),
the $C^0$ interior penalty ($C^0$IP),  and the weakly over-penalized
symmetric interior penalty (WOPSIP) method.
\W{Table~\ref{tab:spaces} below displays the discrete spaces $V_h$ and operators $\W{I},J_h$ introduced 
in Section \ref{sec:Examples of second-order finite element schemes}.
The multiplicative constant $C$ in \eqref{eqn:intro_error_2} exclusively
depends on the
shape regularity of the underlying triangulation.
\begin{table}[bp]
\centering {\footnotesize
	\W{
\begin{tabular}{|c|c|c|cc|c|c|}
\hline
{}{}& Morley & $C^0$IP & \multicolumn{1}{c|} {dG I} & dG II & WOPSIP & Reference\\ \hline
\begin{minipage}{2cm}
\vspace{3pt}
\centering
$V_h$
\end{minipage}
&  ${\rm M}(\cT) $  &  $S^2_0(\cT)$     & \multicolumn{1}{c|}{ $P_2(\cT) $ }         & $P_2(\cT) $ & $P_2(\cT) $
& \eqref{eqn:V_h_def} \\ \hline

$I=I_h I_\M:V\to V_h$  & ${\rm id}$    & $I_{\rm C}I_{\rm M}$     & \multicolumn{1}{c|}{$I_{\rm M}$}& $I_{\rm M}$  & $I_{\rm M}$
			  & $I_\M$ in Def.~\ref{defccMorleyinterpolation}, $I_\C$ in \eqref{eq:ic}  \\ \hline
$J_h=J I_\M:V_h\to V$  & $J$    & $JI_{\rm M}$     & \multicolumn{1}{c|}{$JI_{\rm M}$}& $JI_{\rm M}$  & $JI_{\rm M}$      &
$J$ in Lemma \ref{lem:MorleyCompanion}\\ \hline
\end{tabular}
\caption{\W{Five discretizations}}
\label{tab:spaces}}
}
\end{table}
}
The discussion includes the standard and modified schemes that come with and without a smoother $J_h$ on the right-hand
side.
This paper completes the a posteriori error analysis for these lowest-order discretisations and provides novel reliable and efficient a posteriori error estimators for a rather general class of general sources $F\in
V^*$. %

\subsection{Outline}%
Section \ref{sec:a priori} introduces the abstract discretisation scheme with the key assumption \eqref{eqn:H} for the a
priori analysis in \cite{carstensen_lowest-order_2022}.
Section \ref{sec:Abstract a posteriori error analysis} discusses a known abstract error identity and its
application in the a posteriori error analysis. This is followed by the  concept of a quasi-optimal smoother and the a
priori key property \eqref{eqn:H} that lead to an explicit a posteriori error bound with a particular structure of the
test function as in \eqref{eqn:intro_error_2}.
Section \ref{sec:Examples of second-order finite element schemes} provides examples for the abstract setting in terms of
five lowest-order schemes for the biharmonic equation.  
Section \ref{sec:Building block for explicit residual-based a posteriori error estimators} establishes explicit
estimates
for the error contributions of the a posteriori error bound from Section \ref{sec:Abstract a posteriori error analysis}.
Section \ref{sec:Unified a posteriori error control} presents a unified a~posteriori error control for five
lowest-order schemes for the biharmonic equation in a simplified setting with a right-hand side $F\in L^2(\Omega)$ and recovers \cite{DVNJSR2007, HuShi09, BGS10, Georgoulis2011,
BrenGudiSung10, MozoSuli07}.
The restriction to sources in $L^2$ underlines the state of the art before this paper and thereby highlights the new paradigm through
comparison with known results.
The emphasis in 
Section \ref{sec:General sources} is on a class of general sources $F \in V^*$ with a novel a posteriori error estimator of
the residual that is reliable and efficient up to data-oscillations.
\W{Appendix \ref{apx:A posteriori error control of a piecewise polynomial source} shades a different light on the
	discussion in Section \ref{sec:General sources} and provides lower and
upper bounds for the dual norm of a functionals $F\in V^*$.}

The presentation is laid out in two dimensions with shape-regular triangulations into triangles and second-order
discretizations for simplicity; but the arguments apply to 3D as well, cf.~\cite{CCP_new} for a companion operator $J_h$ in 3D.
The abstract results of this paper will be applied to an a posteriori error analysis of semilinear problems
\cite{carstensen_semilinear202+},
where a linearisation enforces (piecewise polynomial) 
$F\in H^{-2}(\Omega)\setminus L^2(\Omega)$ in future research.

\subsection{General notation}%
\label{sub:General notation}
Standard notation on Lebesgue and Sobolev spaces, 
their norms, and $L^2$ scalar products  applies throughout the paper
such as the abbreviation $\|\bullet\|$ for $\|\bullet\|_{L^2(\Omega)}$. 
Recall that the energy norm  
	$\trinl \bullet \trinr:=\|D^2\bullet\|$ is a norm on $H^2_0(\Omega)$.
	Throughout this paper, $\T$ denotes a shape-regular triangulation of a polygonal and bounded (possibly
	multiply-connected) Lipschitz domain $\Omega\subset \mathbb
R^2$ into
triangles.
Let $\mathcal{V}(\Omega)$ and $\mathcal{E}(\Omega)$ denote the set of interior vertices and edges in the triangulation
$\T$ and let
$\mathcal{V}(\partial\Omega)$ and $\mathcal{E}(\partial\Omega)$ denote the boundary vertices and edges.
The gradient and Hessian operators $\nabla_\pw\coloneqq D_\pw$ and $D^2_\pw$ act piecewise on the space
$H^m(\T)\coloneqq \prod_{T\in\T}H^m(T)$ of piecewise Sobolev functions for $m=1,2$ with the abbreviation $H^m(K)\coloneqq
H^m(\mathrm{int}\;K)$ for a triangle or edge $K\in \T\cup \E$ with relative interior $\mathrm{int}(K)$.
The space $P_{\hspace{-.13em}k}(K)$ of polynomials of total degree at most $k\in\mathbb N_0$ on $K\in\T\cup\E$ with diameter $h_K$ defines
the space of piecewise polynomials
\begin{align*}
	P_{\hspace{-.13em}k}(\T)&\coloneqq\{p\in L^\infty(\Omega) : p|_{T}\in P_{\hspace{-.13em}k}(T)\text{ for all } T\in \T\}.
\end{align*}
The mesh-size $h_\T\in P_{\hspace{-.13em}0}(\T)$ is the piecewise constant function with $h_\T|_T\equiv h_T\coloneqq\mathrm{diam}(T)$ for all $T\in\T$.
Throughout this paper, let $H^k(\Omega; X),
H^k(\T; X)$, resp.~$P_{\hspace{-.13em}k}(\T; X)$ denote the space of (piecewise) Sobolev functions resp.~polynomials with values
in $X=\mathbb R^2, \mathbb R^{2\times2}, \mathbb S$ for $k\in\mathbb N_0$; $\mathbb S\subset\mathbb R^{2\times2}$ is the set of symmetric $2 \times 2$ matrices.
The spaces $H^{-k}(\Omega)\coloneqq(H^k_0(\Omega))^*$ are the dual spaces of $H^k_0(\Omega)$ for $k\in \mathbb N$.
Given any function $v \in L^2(E)$ on an edge $E\in \E$, 
define the integral mean $ \fint_E v \dx:= h_E^{-1} \int_E v \dx$. 
The notation $A \lesssim B$ abbreviates $A \leq CB$ for some positive generic constant $C$, 
which exclusively depends on the shape-regularity of the underlying triangulation $\T$;
 $A\approx B$ abbreviates $A\lesssim B \lesssim A$. 

\section{Unified a priori error analysis}\label{sec:a priori}
Nonstandard schemes compute discontinuous approximations in general and require a
smoother to map the discrete functions  into the continuous space $V$.

\subsection{Discretisation}\label{sec:abstractSeconddiscretisation}
Given the Hilbert space ($V, a$) from the continuous problem \eqref{eqn:WP}, consider some bigger Hilbert space $(\widehat V, \widehat
a)$ that contains $V\subset\widehat V$ as well as the discrete spaces $V_h, V_\nc\subset \widehat V$.
Let $\widehat a\coloneqq a_\pw + j_h$ be the sum of 
the
semi-scalar products
$a_\pw,j_h:\widehat V\times\widehat V\to\mathbb R$ where,  $a_\pw$ extends $a=a_{\pw}|_{V\times V}$ and is a scalar product with induced norm $\trb{\bullet}_\pw\coloneqq a_\pw(\bullet, \bullet)^{1/2}$
in $V+V_\nc$. 
The semi-scalar product $j_h:\widehat
V\times\widehat V\to\mathbb R$ represents jumps that vanish in $V + V_\nc$, i.e., $j_h(v, \bullet) = 0 $ for any $v\in
V+V_\nc$.
The induced norm on $\widehat V$ reads \begin{align}\label{eqn:pw_h} 
	\|\bullet\|_h\coloneqq \big(\trb{\bullet}^2_\pw + j_h(\bullet,
	\bullet)\big)^{1/2}\quad \text{and satisfies}\quad \trinl \bullet \trinr_\pw=\|\bullet \|_h \text{ in } V+\Vnc.
\end{align} 
The discretisation consists of a finite-dimensional trial and test space $V_h$
with respect to a shape-regular triangulation $\T$ of $\Omega$ and the  (possibly unsymmetric) bilinear form
\[
a_h:({ V}+V_h+\Vnc)\times ({ V}+V_h+\Vnc)\rightarrow {\mathbb R}.
\]
We assume
 that $a_h$ is $V_h$-elliptic and bounded on $V_h$ with respect to 
$\|\bullet \|_h$ in the sense that some universal constants $ 0< \alpha\leq M<\infty $ satisfy, for all  $v_h,w_h \in
V_h$, that
\begin{align}\label{eq:a_ellipticity} 
	\alpha \| v_h \|_h^2 \le a_h(v_h,v_h)\quad\text{ and }\quad a_h(v_h, w_h) \le M \|v_h \|_h \|w_h\|_h  .
\end{align}
Since $V_h\subsetneq V$ is \emph{not} a subset of $V$, the evaluation $F(v_h)$ at $v_h\in V_h$ is \emph{not} well-defined
	for general $F\in V^*$.
	Therefore many of the earlier contributions, in particular to the a posteriori error control, merely consider $F\in
	L^2(\Omega)$ whenever $\widehat V\subset L^2(\Omega)$.
	The series of papers \cite{veeser_zanotti1, veeser_zanotti3,
veeser_zanotti2} advertise a smoother $Q\in L(V_h;V)$ to evaluate the modified source
	$F(Qv_h)$ on the discrete level.
This paper complements those contributions on the a priori error analysis by reliable and efficient a posteriori error
estimates. 
This is itself highly relevant in scientific computing and a first step towards adaptive mesh-refining\X{ \cite{carstensen_optimal_2022}}.

	To be more general, this paper considers a rather general class of sources that allow an extension $\hatF\in \widehat V^*$ of $F=\hatF|_V$.
The Lax-Milgram lemma ensures the existence of a unique discrete solution $u_h\in V_h$ to 
\begin{equation}\label{eq:discrete2} 
	a_h(u_h, v_h) = \hatF(Qv_h) \quad\text{for all  } v_h \in V_h
\end{equation}
for the two cases $Q=\id$ (no smoother, but depending on $\hatF$) and $Q=J_h$ for a smoother $J_h\in L(V_h;V)$.
The history of $J_h$ is related to averaging techniques and dates back to the analysis of the Crouzeix-Raviart method 
\cite{ccdgmsMathComp, aCCP, veeser_zanotti2} for the
reliable error
control \cite{carstensen_computational_2014}.
An earlier motivation was the construction of intergrid transfer operators in the convergence analysis of multigrid methods for
nonconforming schemes \cite{BrennerMR1620215}.

The first results will be derived for $\hatF\equiv f\in L^2(\Omega)$ to recover known results in a unified framework, while
Section \ref{sec:General sources} specifies a large class of extended sources $\hatF$ and provides novel a posteriori
error estimates with and without smoother.

\subsection{Quasi-best approximation}%
\label{sub:Quasi-optimal approximations}
The abstract framework from \cite{carstensen_lowest-order_2022} provides a tool for the a priori analysis therein. %
\begin{defn}[quasi-optimal smoother]
	An operator $J_h\in L(V_h;V)$ is called a quasi-optimal smoother if there exists a constant $C_{\rm J} \ge
	0$ such that 
\begin{align}\label{quasioptimalsmoother}
	\|v_h -J_hv_h\|_h \le C_{\rm J} \min_{v\in V}\|v-v_h\|_{h} \quad\text{ for all } v_h \in V_h.
\end{align}
\end{defn}
All the examples in  \cite{veeser_zanotti1, veeser_zanotti3, veeser_zanotti2} discuss $ J_h \in L(V_h; V)$ with $J_h=\id$ in $V_h\cap V$.
The framework in \cite{carstensen_lowest-order_2022}  introduces a smoother  that satisfies \eqref{quasioptimalsmoother}
and is quasi-optimal with a constant $ C_{\rm J}\approx 1$.
The interpretation is that $J_hv_h\in V$ is a good approximation of $v_h\in V_h$ and provides a bridge between the
discrete objects in $V_h$ and $V$.

The key assumption \eqref{eqn:H} connects the bilinear forms $a$ from \eqref{eqn:WP} and $a_h$ from \eqref{eq:discrete2} and requires
the existence of $\Lambda_{\rm H}\geq 0$ with 
\begin{align*}\tag{\bf H}\label{eqn:H}
  {a}_h(w_h, v_h)  -  {a}(J_h w_h,J_hv_h)   \le \Lambda_{\rm H} \|w_h-J_hw_h\|_{h} \|v_h\|_{h}
					  &\text{ for all }w_h,v_h \in V_h.
\end{align*}
This assumption leads to quasi-optimality of $u_h$ in the discrete norm $\|\bullet\|_h$ 
and holds for a class of problems including the
examples in \cite{carstensen_lowest-order_2022} except WOPSIP.  A key step is therefore the design of a quasi-optimal smoother, e.g., $J_h=J\circ I_\nc$ with the conforming
companion $J$ and a generalised interpolation operator $I_\nc$. 
\begin{theorem}[quasi-best approximation]\label{thm:QO}
	Given an operator $J_h\in L(V_h; V)$ with \eqref{quasioptimalsmoother} and \eqref{eqn:H}, there exists  
	a constant $C_\qo>0$ (that exclusively  depends on $\alpha, M, C_{\rm J}, \Lambda_{\rm H}, $ and $\|J_h\|$)
	such that the exact solution $u\in V$ to \eqref{eqn:WP} and the discrete solution $u_h$ to \eqref{eq:discrete2}
	satisfy
\begin{align}\label{eqn:C_qo}\tag{{\bf QO}}
	\|u-u_{h}\|_h\leq C_\qo\min_{v_h\in V_h}\|u-v_h\|_h.
\end{align}
\end{theorem}
\begin{proof}
This is proven in {\cite[Thm. 2.3]{carstensen_lowest-order_2022}}  inspired by the seminal work \cite{veeser_zanotti1}.
\end{proof}
A stronger version ${\widehat{ \text{\bf (H)}}}$ of \eqref{eqn:H} in \cite[Sec.\ 6]{carstensen_lowest-order_2022} even %
leads to a priori error bounds in weaker (piecewise) Sobolev norms.

\subsection{Transfer operators} \label{sec:transfer}
The error analysis requires transfer operators with certain approximation properties between the three subspaces $V,
V_h, V_\nc$ of $\widehat V$.
Throughout this paper, assume there are three linear operators $I_h\in L(\Vnc;V_h), I_\nc\in L(V + V_h + \Vnc;\Vnc)$, and the
conforming companion operator $J\in L(V_\nc; V)$ and constants $\Lambda_h,\Lambda_\nc,\Lambda_\jc \ge 0 $
such that
\begin{align} \label{eqn:Ih_approx}
	\|v_\nc - I_h v_\nc\|_h &\le \Lambda_h \;\min_{v\in V}\trinl   v-v_\nc \trinr_\pw&&\text{ for all }v_\nc\in
	\Vnc,\\\label{eqn:Inc_approx}
	\|v_h- I_\nc v_h \|_h & \le \Lambda_\nc \;\min_{v\in V}\|v -v_h\|_h&&\text{ for all }v_h \in V_h,\\
	\trinl v_\nc- J v_\nc \trinr_\pw &\le \Lambda_\jc\;  \min_{v\in V}  \trinl v- v_\nc \trinr_\pw &&\text{ for all } v_\nc
	\in \Vnc.\label{eqn:J_approx}
\end{align}
{
\begin{minipage}{0.5\textwidth}
\vspace{-0.3in}
\begin{center}
\begin{longtable}{ |p{.2cm} p{3.7cm}|}
\hline 
 $I_{\rm nc}$ & {$\in L(V+V_{\rm nc}+V_h; V_{\rm nc})$}  \\  \hline
     $I_{h}$ & $\in L(V_{\rm nc}; V_{h})$   \\ \hline
 $J$ & $\in L(V_{\rm nc}; V)$  \\  \hline
$J_h$ & $\in L(V+V_{\rm nc}+V_h; V)$%
  \\  \hline
\caption{  Operators\label{long1}}\\
\end{longtable}
\end{center}
\end{minipage}%
\hfill
\begin{minipage}{0.4\textwidth}
\begin{tikzpicture}[node distance=3cm, auto]\label{fig:composition}
\node (A) {$V+V_{\rm nc} +V_h$};
\node(B) [right of=A] {$V$};
\node (C) [below of=A] {$V_{\rm nc}$};
\draw[->](A) to node {${J_h= J \circ I_{\rm nc}}$}(B);
\draw[->](A) to node [left] {$I_{\rm nc}$}(C);
\draw[->](C) to node [below=0.5ex] {$J$}(B);
\end{tikzpicture}

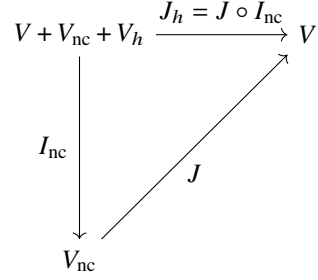
\captionof{figure}{Definition of $J_h$}
\end{minipage}%

}

\medskip \noindent 
Two immediate consequences on the abstract level at hand shall be utilized below.
\begin{lemma}[intermediate bound]
	\label{lem:IIJ_approx}
	Given any $v\in V$ and $v_\nc\in \Vnc$,  \eqref{eqn:Ih_approx}--\eqref{eqn:J_approx} imply
	\begin{align}\label{eqn:IIJ_approx}
		\trb{v - JI_\nc I_h v_\nc}&\leq (1+\Lambda_\jc)(1+\Lambda_\nc)(1+\Lambda_h)\trb{ v-v_\nc}_\pw.
	\end{align}
\end{lemma}
\begin{proof}
	Let $w_\nc\coloneqq I_\nc I_h v_\nc$ and $w_h\coloneqq I_h v_\nc$.
	The triangle inequality and \eqref{eqn:J_approx} show
	\begin{align*}
		\trb{v - &Jw_\nc}\leq \trb{ v-w_\nc}_\pw+\trb{(1-J)w_\nc}_\pw \leq
		(1+\Lambda_\jc)\trb{v-w_\nc}_\pw.
	\end{align*}
	Note $\trb{v-w_\nc}_\pw=\|v-w_\nc\|_h$ from \eqref{eqn:pw_h}. The triangle inequality and \eqref{eqn:Ih_approx}--\eqref{eqn:Inc_approx}  show
	\begin{align*}
		\|v-w_\nc\|_h &\leq \|v-w_h\|_h + \|(1-I_\nc)w_h\|_h\leq
								  (1+\Lambda_\nc)\|v - w_h\|_h,\\
		\|v-w_h\|_h &\leq \|v-v_\nc\|_h + \|(1-I_h)v_\nc\|_h\leq
					(1+\Lambda_h)\trb{v - v_\nc}_\pw.
	\end{align*}
	The combination of those estimates establishes \eqref{eqn:IIJ_approx}.
\end{proof}
\noindent The above transfer operators (see Figure \ref{fig:composition}) lead to a quasi-optimal smoother $J_h\coloneqq J \circ I_\nc\in
L(V_h;V)$.   Although $J_h$ maps $V+V_{\rm nc} +V_h \rightarrow V$,  its restriction to $V_h$ plays a central role  in the sequel.
\begin{lemma}[quasi-optimal smoother]\label{lem:J_h}
	Given any $v_h\in V_h$,  and $J_h\coloneqq
J \circ I_\nc\in L(V_h; V)$, 
	\eqref{eqn:Inc_approx}--\eqref{eqn:J_approx}  show \eqref{quasioptimalsmoother} with
	$C_\jc\coloneqq\Lambda_\nc+\Lambda_\jc+\Lambda_\jc\Lambda_\nc$.
\end{lemma}
\begin{proof}
	A triangle inequality with  $v_\nc\coloneqq I_\nc v_h$, and
	\eqref{eqn:J_approx} verify
	\begin{align*}
		\|v_h-Jv_\nc\|_h\leq \|v_h - v_\nc\|_h + \Lambda_\jc \left(\|v-v_h\|_h + \|v_h - v_\nc\|_h\right)
	\end{align*}
for an arbitrary  $v\in V$.  This and \eqref{eqn:Inc_approx} conclude the proof.
\end{proof}
\noindent Lemma \ref{lem:J_h} shows that $J_h$ %
is a
quasi-optimal smoother with the following property.
\begin{theorem}[quasi-best approximation {\cite{carstensen_lowest-order_2022}}] \label{thm:abstract_main} 
	Let 
$u \in V$ resp. $u_h \in V_h$ solve \eqref{eqn:WP} resp. \eqref{eq:discrete2}.
Suppose \eqref{eqn:H}, \eqref{eqn:pw_h}--\eqref{eq:a_ellipticity}, and
\eqref{eqn:Ih_approx}--\eqref{eqn:J_approx}. Then %
\begin{align*}
	\trinl  u- J_h u_h \trinr + \| u - u_h \|_h &\lesssim \min_{v_\nc\in\Vnc}\trinl u - v_\nc \trinr_\pw.
\end{align*}
\end{theorem}
\begin{proof}
	Lemma \ref{lem:J_h} and Theorem \ref{thm:QO} verify 
\eqref{eqn:C_qo} for $J_h$.
	{A triangle inequality, \eqref{eqn:Ih_approx}, and $\trinl \bullet \trinr_\pw=\|\bullet \|_h$ in $V+V_{\rm nc}$  verify}
	\begin{align*}
		\|u-u_h\|_h\leq C_\qo\|u - I_h v_\nc\|_h\leq C_\qo(1+\Lambda_h)\trb{u-v_\nc}_\pw
	\end{align*} for arbitrary $v_\nc\in V_\nc$.
The proof of Lemma \ref{lem:IIJ_approx} shows $\trb{v-J_h w_h}\leq (1+\Lambda_\jc)(1+\Lambda_\nc)\|v-w_h\|_h$ for
an arbitrary $v\in V, w_h\in V_h$. The combination with the previously displayed inequality concludes the proof.
\end{proof}

\section{Abstract a posteriori error analysis}%
\label{sec:Abstract a posteriori error analysis}
The abstract error identity in Subsection \ref{sub:Abstract error identity}  reveals that $\trb{Res}_*$ is a contribution to the error.
Subsection \ref{sub:Evolution of a posteriori error analysis} revisits the Crouzeix-Raviart and Morley FEM and recalls
known bounds thereof.
Subsection \ref{sub:Paradigm} explains a paradigm shift towards a universal error analysis that is explicit in the structure of the test
function through a quasi-optimal smoother and the property \eqref{eqn:H}.
\vspace{-0.2in}
\subsection{Abstract error identity for $F \in V^*$}%
\label{sub:Abstract error identity}

Given the exact solution $u\in V$ to \eqref{eqn:WP} and the discrete solution $u_h\in V_h$ to \eqref{eq:discrete2}, the
natural error $u-u_h\in V+V_h\subset \widehat V$ can be measured in the norm $\|\bullet\|_h$ from Subsection
\ref{sec:abstractSeconddiscretisation}.
This allows a well-known split with the residual $Res\coloneqq
F-a_\pw(u_h,\bullet)\in V^*$ \cite{CCHJ2007}.
\begin{theorem}[error identity]\label{thm:Pythagoras}  The exact solution $u\in V$ to \eqref{eqn:WP} and the discrete solution $u_h\in V_h$ to \eqref{eq:discrete2} satisfy
	\begin{align}\label{eqn:err_ident}
		\|u-u_h\|_h^2=\trb{Res}_*^2 + \min_{v\in V}\|v-u_h\|_h^2.
	\end{align}
\end{theorem}
\begin{proof}
	Let $w\in V$ be the Riesz representation of the linear and bounded functional $a_\pw(u_h, \bullet)\in V^*$ in the Hilbert space $(V, a_\pw)$, so that
	$a_\pw(u_h-w, \bullet) = 0$ in $V$. 
	This orthogonality shows that $w \in V$ is the best-approximation of $u_h\in V_h\subset \widehat V$ in the complete subspace $V$, i.e.,
	\begin{align}\label{eqn:proof1}
		\delta&\coloneqq \trb{w-u_h}_\pw = \min_{v\in V}\trb{v-u_h}_\pw
	\end{align}
and  allows for the Pythagoras identity
	\begin{align}
		\trb{u-u_h}_\pw^2&=\trb{u-w}^2 + \trb{w-u_h}_\pw^2.
	\end{align}
	The orthogonality also shows, for all $v\in V$, that
	\begin{align*}a(u-w, v) = a(u, v) - a_\pw(u_h, v) = Res(v)\end{align*}
	with $a(u, \bullet) = F$ in $V$ in the last step.
	In other words, $u-w$ is the Riesz representation of $Res\in V^*$ in the Hilbert space $(V,a)$ and the Riesz
	isomorphism reveals
	\begin{align}\label{eqn:proof3}
		\trb{u-w}&=\trb{Res}_*\coloneqq\sup_{v\in V\setminus\{0\}}\frac{Res(v)}{\trb{v}}.
	\end{align}
	The summary of \eqref{eqn:proof1}--\eqref{eqn:proof3} reads $\trb{u-u_h}_\pw^2 = \trb{Res}_{*}^2 + \delta^2$.
	Since ${j_h(\bullet, v) = j_h(v,\bullet) = 0}$, the proof concludes with $\|v-u_h\|_h^2=\trb{v-u_h}_\pw^2 + j_h(u_h,
	u_h)$ for
	any $v\in V$. 
\end{proof}
\begin{rem}[explicit a posteriori bounds]
The proof of Theorem \ref{thm:Pythagoras} is nothing but a Pythagoras identity and serves as an idealisation: 
While
$j_h(u_h, u_h)$ comes for free, the computation of $\trb{Res}_*$ or of $\delta=\min_{v\in V}\|v-u_h\|_h^2 - j_h(u_h,
u_h)$ is far too costly. 
Instead, the error identity rather serves as a guide to design individual upper bounds of $\delta$ and $\trb{Res}_*$.
The a priori error analysis of Section \ref{sec:abstractSeconddiscretisation} provides a
quasi-optimal smoother $J_h\in L(V_h;V)$.
Then 
\eqref{quasioptimalsmoother} shows
\begin{align}
	\label{eqn:proof4}
	\min_{v \in V}\|v-u_h\|_h&\leq\| u_h-J_hu_h\|_h\leq C_{\rm J}\min_{v \in V}\|v-u_h\|_h.
\end{align}
In the language of a posteriori error control, \eqref{eqn:proof4} asserts the reliability and efficiency of the a posteriori estimator
$\muold$ of the error $\min_{v\in V}\|v-u_h\|_h$. %
This ends the  discussion of $\muold$ and motivates the focus on bounds of $Res$ below.
\end{rem}
In order to understand the difference between the classical and the current treatment, the two simplest nonconforming schemes
will be discussed in the subsequent subsection.

\subsection{Crouzeix-Raviart and Morley FEM}%
\label{sub:Evolution of a posteriori error analysis}
This subsection motivates the abstract a posteriori error analysis by a recollection \cite{DEDRPCVV1996,
CCSBSJ2002,CCHJ2007,CCHJOA2007,VR2012,CCEMHRHWLC2012,CCDGJH14} for $m=1$ and 
 \cite{HuShi09,HuShiXu2012, CCDGHU13, DVNJSR2007} for $m=2$ of the simplest nonconforming
schemes for the $m$-harmonic equation  $(-\Delta)^m u = f$ for $m=1,2$ with right-hand function $f \in L^2(\Omega)$.
The weak solution seeks $u\in
V\coloneqq H^m_0(\Omega)\subset \widehat V\coloneqq H^m(\T)$ to
\begin{align}\label{eqn:WP2}
	a(u,v)&=(f, v)_{L^2(\Omega)} \text{ for all }v\in V
\end{align}
with the energy scalar product $a\coloneqq a_{\pw}|_{V\times V}$ and $a_{\pw}(\bullet, \bullet)\coloneqq (D^m \bullet , D^m \bullet)_{L^2(\Omega)}$
in $\widehat V$.
\subsubsection{Crouzeix-Raviart FEM}%
\label{ssub:Crouzeix-Raviart FEM}
\begin{figure}[]
	\centering
  \begin{tikzpicture}[line width=1pt,scale=4,line join=round,line cap=round]
    \draw (0,0) node[shape=coordinate] (A1) {A1};
    \draw (1,0) node[shape=coordinate] (A2) {A2};
    \draw ($(A1)!.5!(A2)$) node[shape=coordinate] (A12) {A12};
    \draw ($(A12) + (0,1)$) node[shape=coordinate] (A3) {A3};
    \draw ($(A2)!.5!(A3)$) node[shape=coordinate] (A23) {A23};    
    \draw ($(A3) + (1,0)$) node[shape=coordinate] (A4) {A4};
    \draw ($(A2)!.33!(A3)$) node[shape=coordinate] (nuBase) {nuBase};
    \draw ($(nuBase)!1!90:(A2)$) node[shape=coordinate,label=above:$\nu_E$] (nuTop) {nuTop};

    \draw (A2) node[shape=coordinate] (B1) {B1};
    \draw ($(A3) + (1,0)$) node[shape=coordinate] (B2) {B2};
    \draw (A3) node[shape=coordinate] (B3) {B3};

    \draw (A1) -- (A2) -- (A3) -- cycle;
    \draw (B1) -- (B2) -- (B3) -- cycle;
    \draw[->] (nuBase) -- (nuTop);

    \draw (nuBase) let 
                    \p1 = ($(nuTop)-(nuBase)$), 
                    \n2 = {veclen(\x1,\y1)}
                 in 
                 -- ($(nuBase)!.3!(A2)$) arc[start angle=-60,delta 
                   angle=90,radius={.28*\n2}];
    \draw[fill] ($(nuBase)!.15!45:(A2)$) circle[radius=.1pt];

    \node (E) at ($(A23) - .07*(1,1)$) {$E$};
    \node (Tplus) at ($(A12)!.33!(A3)$) (Tplus) {$T_+$};
    \node (Tminus) at ($(A23)!.40!(B2)$) (Tminus) {$T_-$};
	\node (Pplus) at ($(A1)- .05*(1,1)$) (Pplus) {$P_{\hspace{-.13em}+}$};
	\node (Pminus) at ($(A4) + .05*(1,1)$) (Pminus) {$P_{\hspace{-.13em}-}$};
     \node (Alabel) at ($(A2) - .05*(-1,1)$) (Pminus) {$A$};
     \node (Blabel) at ($(A3) + .05*(-1,1)$) (Pminus) {$B$};
  \end{tikzpicture}
  \caption{The interior edge patch $\omega(E)$ and normal $\nu_E=\pm\nu_{T_{\pm}}$	of $E=\partial T_+\cap \partial T_-$}
  \label{fig:Edge_patch}
\end{figure}
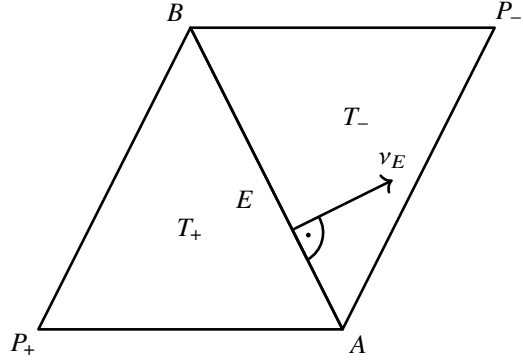
Let $u\in V\coloneqq H^1_0(\Omega)$ be the weak solution to the Poisson model problem, i.e., $u$ solves \eqref{eqn:WP2}
for $m=1$.
The Crouzeix-Raviart finite element space requires the definition of jumps across an edge $E\in\E$ in the triangulation
$\T$.
Let $\nu_T$ be the unit outer normal of $T\in\T$ and fix the orientation of the unit normal $\nu_E$ on every edge
$E\in\E$ with midpoint $\mathrm{mid}(T)$.
Every interior edge $E=\partial T_+\cap\partial T_-\in \E(\Omega)$ has exactly two neighbouring triangles 
$T_+, T_-\in \T$ as in Figure \ref{fig:Edge_patch}, labelled such that $\nu_E=\pm \nu_{T_\pm}|_E$, and the jump
of a piecewise Sobolev function $v\in H^1(\T)$ across $E$ reads $[v]_E\coloneqq v|_{T_+} - v|_{T_-}\in H^1(E)$.
On a boundary edge $E\in \E(\partial\Omega)$, the jump $[v]_E\coloneqq v $ is the unique trace of the function $v\in
H^1(\T)$.
Define the space
\begin{align*}\CR^1_0(\T)\coloneqq\left\{p\in P_{\hspace{-.13em}1}(\T)\ \middle|\
		\begin{array}{lc}
			[p]_E(\mathrm{mid} \: E)=0\text{ vanishes for every edge }E\in\E
	\end{array}\right\}
\end{align*}
of piecewise
affine polynomials over a given shape-regular triangulation $\T$ with continuity at the midpoints of the edges.
This space $V_\nc\coloneqq \CR_0^1(\T)$ comes with the natural interpolation operator $I_\CR:V+V_\nc\to V_\nc$ that maps $v\in
V+V_\nc$ to the unique function $I_\CR v\in V_\nc$ with $\int_E(v-I_\CR v)\ds=0$
for every edge $E\in\E$.
The classical formulation of the lowest-order nonconforming Crouzeix-Raviart FEM approximates the weak solution $u\in
H^1_0(\Omega)$ of \eqref{eqn:WP2} with the discrete solution $u_\CR\in  \CR^1_0(\T) \equiv V_\nc $ to
\begin{align}\label{eqn:DWP2}
	a_\pw(u_\CR,v_\CR)\equiv \int_{\Omega} \nabla_{\pw} u_{\CR} \cdot \nabla_{\pw} v_{\CR} \dx = (f,
	v_\CR)_{L^2(\Omega)}&\text{ for all }v_\CR\in  \CR^1_0(\T).
\end{align}
This is exactly \eqref{eq:discrete2} for the natural choice $\hatF\coloneqq F\equiv f\in L^2(\Omega)$ and without smoother
$Q\coloneqq\id$.
The semi-scalar product $a_\pw$ induces the piecewise energy norm
$\trb{\bullet}_\pw\coloneqq \|\nabla_\pw\bullet\|$ in $V+V_\nc \equiv H^1_0(\Omega) +\CR^1_0(\T)$ \cite{NNCC2020}.
In this particular example, the residual from Section~\ref{sub:Abstract error identity} reads \begin{align*}Res
	\coloneqq (f, \bullet )_{L^2(\Omega)}- a_\pw(u_{\CR}, \bullet)\in
V^*.\end{align*}
\subsubsection{Classical residual-based explicit error estimator }
\label{sub:Classical residual-based CR}
This approach follows \cite{CCHJ2007} and is closely related to the analysis of conforming schemes.
Let $I_\C:H^1_0(\Omega) + \CR^1_0(\T)\to S^1_0(\T)$ denote a quasi-interpolation operator onto the continuous piecewise
affine polynomials $S^1_0(\T)\coloneqq P_{\hspace{-.13em}1}(\T)\cap H^1_0(\Omega)$ with homogeneous boundary conditions.
Since \eqref{eqn:DWP2} holds, the definition of the residual shows $Res(w_\C) = 0 $ for any $w_\C\in S^1_0(\T) \subset \CR^1_0(\T) $, i.e.,
$S^1_0(\T)\subset \mathrm{ker}Res$ lies in the kernel of $Res\in V^*$, and an integration by parts with the test
function $w\coloneqq v-I_\C v$ shows, for $f \in L^2(\Omega)$, that
\begin{align*}%
	Res(v) &= Res(w) = (f, w)_{L^2(\Omega)} - \sum^{}_{T\in\T} \sum^{}_{E\in\E(T)} \int_E \nabla_\pw u_\CR\cdot \nu_E
w\ds\\
		   &=\left(f,w\right)_{L^2(\Omega)}-\sum^{}_{E\in\E(\Omega)} \int_E [\nabla_\pw u_\CR]_E\cdot \nu_E w\ds.
\end{align*}
The last step is a careful resummation over the edges: Each interior edge $E\in\E(\Omega)$ has two contributions (from $T_+$ and $T_-$) with
opposite signs from $\nu_{T_+}=-\nu_{T_-}$ on $E$.
No contributions arise from the boundary edges $E\in\E(\partial\Omega)$ because of $w|_{\partial\Omega}=0$.
Cauchy inequalities show
\begin{align}\label{eqn:Res2}
	Res(v) &\leq \|h_\T f\|\|h_\T^{-1}w\| + \sum^{}_{E\in\E(\Omega)} h_E^{1/2}\|[\nabla_\pw
		   u_\CR]_E\cdot\nu_E\|_{L^2(E)}h_E^{-1/2}\|w\|_{L^2(E)}.
\end{align}
The quasi-interpolation operator $I_\C$ from \cite{Clement1975, ScottZhang90} satisfies the stability estimates
 \[h_T^{-1}\|v-I_{\C} v\|_{L^2(T)}\leq C_{\rm apx}\; \|\nabla
v\|_{L^2(\omega(T))} \text{ in } T\in\T \]
 with a constant $C_{\rm apx}>0$ that exclusively depends on the shape regularity of $\T$.
 Here $\omega(T)$ denotes the layer-1 patch around $T\in\T$.
The trace inequality 
\cite[Eqn.\ (12.17)]{ern_finite_2021-1}
\begin{align*}
	h_E^{-1/2}\|v\|_{L^2(E)}\leq C_{\rm tr}\; \left(h_T^{-1}\|v\|_{L^2(T(E))} + \|\nabla v\|_{L^2(T(E))}\right)\quad \text{ for all }v\in V
\end{align*}
bounds the norms on the edge $E\subset \partial T(E)$ by norms
of some adjacent triangle $T(E)\in\T$ with a constant $C_{\rm tr}>0$ that exclusively depends on the shape-regularity of $\T$.
This and a final Cauchy inequality in $\ell^2$ for the sum in \eqref{eqn:Res2} show
\begin{align}\label{eqn:Res_bound_1}
	\trb{Res}_*\coloneqq \sup_{v\in V\setminus\{0\}}\frac{Res(v)}{\trb{v}} \lesssim \|h_\T f\| + \sqrt{\sum^{}_{E\in\E(\Omega)} h_E\|[\nabla
	u_\CR]_E\cdot \nu_E\|_{L^2(E)}^2}.
\end{align}
The jump term in \eqref{eqn:Res_bound_1} can be bounded by $\|h_\T f\|$ and the simpler 
\begin{align}\label{eqn:CR_Res_bound_simple}
	\trb{Res}_*\lesssim \|h_\T f\|
\end{align} estimate without normal jumps is possible.  
For any interior edge $E\in\E(\Omega)$, the edge-patch $\omega(E)\coloneqq
	\mathrm{int}(T_+\cup T_-)$ is the union of the two neighboring triangles $T_+, T_-\in \T$.
\begin{lemma}[bound without jumps]\label{lem:Res_bound1_no_jump}
	The normal jumps from \eqref{eqn:Res_bound_1} satisfy %
$$h_E^{1/2}\|[\nabla
u_\CR]_E\cdot \nu_E\|_{L^2(E)}\lesssim \|h_\T f\|_{L^2(\omega(E))}\quad\text{for any }E\in\E(\Omega).$$ %
\end{lemma}
\begin{proof}

 Recall the edge-oriented basis function  $\psi_E\in CR^1_0(\T)$ as the unique
	function in $\CR^1_0(\T)$ with $\psi_E(\mathrm{mid}\; E) = 1$ and $\psi_E(\mathrm{mid}\; F) = 0$ for every other edge
	$F\in\E\setminus \{E\}$.
	Since $\psi_E\in \CR^1_0(\T)$ is piecewise affine, its support $\overline{\omega(E)}$ is the edge-patch $\omega(E)$
with 
$\psi_E\equiv 1$ on $E=\partial
	T_{+}\cap\partial T_-$.
	This, an integration by parts for the interior edge $E\in\E(\Omega)$, and \eqref{eqn:DWP2} prove for $\beta \coloneqq [\nabla u_\CR]_E\cdot \nu_E\in\mathbb R$, 
	\begin{align*}
		\|[\nabla
		u_\CR]_E\cdot \nu_E\|_{L^2(E)}^2 &=\beta\!\! \int_E [\nabla
		u_\CR]_E {\cdot \nu_E} \; \psi_E\ds = \beta\!\! \int_{\omega(E)} \nabla u_\CR\cdot \nabla \psi_E\dx=\beta\, ( f, 
		\psi_E)_{L^2(\omega(E))}.
	\end{align*}
	The midpoint quadrature rule shows $\|\psi_E\|_{L^2(T)}^2 = |T|/3\approx h_T^2\approx h_E^2$ in 2D by shape-regularity. 
Since $\|\psi_E\|_{L^2(E)}^2 = |E| = h_E$, the previous displayed identity, a Cauchy inequality,  {and the definition of $ \beta$} verify
	\begin{align*}
		\|[\nabla
		u_\CR]_E\cdot \nu_E\|_{L^2(E)}^2\lesssim \|h_\T f\|_{L^2(\omega(E))}h_E^{-1/2}\|[\nabla
		u_\CR]_E\cdot \nu_E\|_{L^2(E)}.
	\end{align*}
	This concludes the proof of \eqref{eqn:CR_Res_bound_simple}. 
\end{proof}
\subsubsection{Bound from Crouzeix-Raviart interpolation}%
\label{ssub:Bound from interpolation}
The integration by parts formula on $T\in\T$ and the definition of the natural interpolation $I_\CR:V+V_\nc\to
V_\nc$ reveal
\begin{align*}
	\int_T \nabla(v-I_\CR v)\cdot \nabla p_1\dx= \sum^{}_{E\in\E(T)} [\nabla_\pw p_1]_E\cdot\nu_E\int_E(v-I_\CR v)\ds=0%
\end{align*}
for any $v\in V+V_\nc, p_1\in P_{\hspace{-.13em}1}(\T)$ 
and $(v-I_\CR v)\perp P_{\hspace{-.13em}1}(\T)$ is $a_\pw$-orthogonal to $P_{\hspace{-.13em}1}(\T)\supset \CR^1_0(\T)$.
This, \eqref{eqn:DWP2}, and the interpolation error estimate $\|h_\T^{-1}(v-I_\CR v)\|\leq \kappa_{\CR}\trb{v}$ from \cite[Sec.\ 4]{ccdg_2014} with $\kappa_{\CR}=(1/48
+ j_{1,1}^2)^{1/2}\leq 0.2983$ for the first positive root $j_{1,1}$ of the Bessel function of the first kind
show 
\begin{align}\label{eqn:Res_bound_2}
	\trb{Res}_*\coloneqq \sup_{v\in V  \setminus  \{0\}}\frac{Res(v)}{\trb{v}}=\sup_{v\in V \setminus \{0\}}\frac{(f, v-I_\CR
	v)_{L^2(\Omega)}}{\trb{v}}\leq \kappa_{\CR}\|h_\T f\|.
\end{align}
The difference to the bound in Lemma~\ref{lem:Res_bound1_no_jump} is not only the explicit control in terms of the smaller constant $\kappa_{\CR}$,
but above all, that the methodology directly controls $\trinl Res \trinr_*$ as in \cite[p.~317]{CCDGHU13} without jump terms. The latter also follows from \eqref{eqn:Res_bound_1} and Lemma \ref{lem:Res_bound1_no_jump}. 

\noindent The key observation is that this technique does \emph{not} need any conforming subspace $S^1_0(\T) \subset
\CR^1_0(\T)$ and this is a relevant advance for the application to the Morley FEM.
\subsubsection{Morley FEM}%
\label{ssub:Morley}
Let $u\in V\coloneqq H^2_0(\Omega)$ be the weak solution to the biharmonic equation $\Delta^2 u = f\in L^2(\Omega)$, i.e., $u$ solves \eqref{eqn:WP2}
for $m=2$.
Define the normal 
jump $[\partial v/\partial\nu_E]_E\coloneqq [\nabla v\cdot
\nu_E]_E$ of a function $v\in H^2(\T)$ along an edge $E\in\E$.
The Morley function space
\begin{align}\label{eqn:Morley}
	\M(\T)\coloneqq\left\{p\in P_{\hspace{-.13em}2}(\T)\ \middle|
		\begin{array}{lc}
			p(z)\text{ is continuous at every }z\in\V(\Omega)\text{ and } p|_{\mathcal{V}(\partial\Omega)}=0,\\
			{[\nabla_\pw p\cdot\nu_E]}_E(\mathrm{mid} \; E) = 0 \text{ vanishes for every edge }E\in\E
	\end{array}\right\}
\end{align}
comes with a natural interpolation operator $I_\M:H^2_0(\Omega)+\M(\T)\to \M(\T)$.
\begin{defn}[classical Morley interpolation {\cite{ccdg_2014,
	BrennerSungZhang13}}]\label{defccMorleyinterpolation_classic}
	Given any function $v \in H^2_0(\Omega)+ \M(\T)$, the Morley interpolation operator $I_\M:H^2_0(\Omega)+ \M(\T)\to \M(\T)$
	defines $I_\M v\in \M(\T)$ by
	\begin{align*}
		(v-I_\M v)(z)=0 
\quad\text{ for }z\in\mathcal{V}(\Omega)\text{ and } 
		\fint_E\frac{\partial (v - I_\M v)}{\partial \nu_E} \ds =0\quad\text{ for
 }E\in\E(\Omega).
	\end{align*}
\end{defn} 
\noindent This interpolation operator 
possesses the $a_\pw$-orthogonality property ${v - I_\M
v\perp_{a_{\pw}} P_{\hspace{-.13em}2}(\T)}$ for any $v\in H^2_0(\Omega)+ \M(\T)$.
The nonconforming Morley FEM approximates $u\in H^2_0(\Omega)$ with the unique discrete solution $u_\M\in 
\M(\T) =: V_\nc$ to
\begin{align}\label{eqn:DWP3}
	a_\pw(u_\M,v_\M)\coloneqq\int_\Omega D^2_\pw u_\M: D^2_\pw v_\M \dx= (f, v_\M)_{L^2(\Omega)} &\text{ for all }v_\M\in
	 \M(\T).
\end{align}
This represents \eqref{eq:discrete2} for $\widehat F\coloneqq F\equiv f\in L^2(\Omega)$ and $Q=\id$ while Section
\ref{sub:Morley method} considers $Q=\id$ and $Q=J_h$ simultaneously in a new a posteriori analysis and Section
\ref{sec:General sources} discusses general sources $F\in V^*$.
Note that $a_\pw$ is a scalar-product in $V+V_\nc\equiv H^2_0(\Omega)+\M(\T)$
\cite{NNCC2020}.
The 
residual from Section \ref{sub:Abstract error identity} reads $Res \coloneqq (f,  \bullet)_{L^2(\Omega)}  - a_\pw(u_{\M}, \bullet)\in
V^*$.
\subsubsection{Bounds from Morley interpolation}%
\label{ssub:Bounds on the residual for Morley}
An approach similar to Subsection \ref{sub:Classical residual-based CR} for the Crouzeix-Raviart FEM fails immediately because
$S^2_0(\T)\cap H^2_0(\Omega)$ is not rich enough: For many triangulations $S^2_0(\T)\cap H^2_0(\Omega)=\{0\}$ is
trivial, however not in general \cite[Sec.\ 3.3]{scott_c_2019}. 

However, the $a_\pw$-orthogonality ${v-I_\M v\perp_{a_{\pw}} P_{\hspace{-.13em}2}(\T)}$ for all $v\in V$ with the Morley interpolation
$I_\M$ allows the arguments from Subsection \ref{ssub:Bound from interpolation}
that lead in \cite{DVNJSR2007,HuShi09} to
\begin{align}\label{eqn:Res_bound_M}
	\trb{Res}_*\coloneqq \sup_{v\in V}\frac{Res(v)}{\trb{v}}=\sup_{v\in V}\frac{(f, v-I_\M
	v)_{L^2(\Omega)}}{\trb{v}}\leq \kappa_{\M}\|h_\T^2 f\|.
\end{align}
The interpolation error estimate $\|h_\T^{-2}(v-I_\M v)\|\leq \kappa_\M\trb{v}$ holds with
constant $\kappa_\M\leq 0.2575$ 
\cite[Sec.\ 4]{ccdg_2014}.

\subsection{Paradigm of unified a posteriori error analysis}%
\label{sub:Paradigm}
The discussion in this subsection departs from the error identity of Theorem \ref{thm:Pythagoras} that includes the dual
norm $\trb{Res}_*$ of the residual $Res\in V^*$.
Recall that $u\in V$ solves \eqref{eqn:WP} in $V$ and $u_h\in V_h$ solves
\eqref{eq:discrete2}.

Subsection \ref{sub:Abstract error identity} discussed the error identity \eqref{eqn:err_ident} with the dual norm of the
residual given as a supremum over all continuous test functions.
Since $u_h\not\in V$ in general, the computable (conforming) post-processing $J_hu_h\in V$ serves as its approximation
and motivates the error definition
$e\coloneq u-J_hu_h\in V$ on the continuous level and $I_h I_\nc u- u_h\in V_h$ on the discrete level.
The efficient error estimator
	$\|u_h-J_hu_h\|_h $
	from \eqref{eqn:proof4} is computable and a triangle inequality in the norm $\|\bullet\|_h$ and \eqref{eqn:pw_h} lead to 
	\begin{align}\label{eqn:u_u_h_mu}
	\|u-u_h\|_h\leq \trb{e} + \muold \text{ and } { \trb{e}\leq \|u-u_h\|_h + \muold.}
\end{align}

\noindent Recall $\hatF|_V=F$ from \eqref{eq:discrete2}. The first argument to establish an alternative abstract error
bound applies the continuous (resp.\ discrete)
equation \eqref{eqn:WP} (resp.~\eqref{eq:discrete2}) to the test function $J_he_h\in V$
(resp.\ $e_h\coloneqq I_h I_\nc e\in V_h$), namely
\begin{align}\label{eqn:key_identity}
	a_h(u_h, e_h) = \hatF(Qe_h) = a(u, J_he_h) - \hatF(J_he_h-Qe_h).
\end{align}
For $Q=J_h$, the last term vanishes and \eqref{eqn:key_identity} becomes the key identity $a_h(u_h, e_h) = a(u,
J_he_h)$.
The second argument is the link of $a_h(u_h, e_h)$ to $a(J_hu_h, J_he_h)$ by \eqref{eqn:H},
\begin{align}\label{eqn:H_corollary}
	a_h(u_h, e_h) - a(J_hu_h, J_he_h)\leq \Lambda_{\rm H}\|u_h-J_hu_h\|_h\|e_h\|_h.%
\end{align}
The (generalized) key identity \eqref{eqn:key_identity} shows that the left-hand side of \eqref{eqn:H_corollary} is equal to
$a(e,J_he_h)-\hatF(J_he_h-Qe_h)$.
This %
 and the abbreviation $w\coloneqq e-J_he_h$ show
\begin{align*}
	\trb{e}^2=a(e,w) + a(e, J_he_h)\leq F(w)-a(J_hu_h, w)+\hatF((J_h-Q)e_h)+\Lambda_{\rm H}\muold\|e_h\|_h
\end{align*}
with $a(u,w)=F(w)$ in the last step. 
This, the Cauchy inequality $a_\pw(u_h- J_hu_h, w)\leq\muold\trb{w}$ using \eqref{eqn:pw_h}, and the residual
$Res\coloneqq F - a_\pw(u_h,\bullet)\in V^*$
reveal
\begin{align}\label{eqn:e_reliable}
	\trb{e}^2\leq \left(\trb{w}+\Lambda_{\rm H}\|e_h\|_h\right)\muold + Res(w)+\hatF(J_he_h-Qe_h).
\end{align}
\begin{theorem}[alternative abstract error bound]\label{thm:a_posteriori_bound}
Let $J_h\in L(V_h; V)$ be a quasi-optimal
smoother and suppose \eqref{eqn:Ih_approx} -- \eqref{eqn:J_approx} and  \eqref{eqn:H}.
Then there exists  a constant $\const{cst:1}>0$ such that the error $e\coloneqq u-J_hu_h\in
	V$ for the solution  $u\in V$ to \eqref{eqn:WP} and  $u_h\in V_h$ to \eqref{eq:discrete2} satisfies
	\begin{align}\label{eqn:a_posteriori_bound}
		\|u-u_h\|_h^2 + \trb{e}^2\leq \const{cst:1}^2\left(\|u_h-J_hu_h\|_h^2 + Res(e-J_hI_hI_\nc e)+\hatF(J_he_h-Qe_h)\right).
	\end{align}
\end{theorem}
\begin{proof}
	Abbreviate $w\coloneqq e-J_he_h\in V$ with $e_h\coloneqq I_h I_\nc e\in V_h$.
	Lemma \ref{lem:IIJ_approx} leads to $\const{cst:2}^{-1}\trb{w}\le  \trb{e- I_\nc e}_{\rm pw} \leq (1+\|I_\nc\|)\trb{e}$
for $\const{cst:2}:= (1+\Lambda_{\rm J})(1+ \Lambda_{\rm nc})(1+\Lambda_h)$  and the 
operator norms control $\|e_h\|_h\leq \|I_\nc\|\|I_h\|\trb{e}$.
This, a Young inequality, and \eqref{eqn:e_reliable} show
\begin{align*} \frac{1}{2}\trb{e}^2\leq \frac{1}{2}\const{cst:3}^2\muold^2 + Res(w)+\hatF(J_he_h-Qe_h)\end{align*}
with 
\W{$\const{cst:3}\coloneqq\const{cst:2}(1+ \|I_{\rm nc} \|) +\Lambda_{\rm H} \|I_{\rm nc} \| \|I_h\| $}. This and \eqref{eqn:u_u_h_mu}
conclude the proof of \eqref{eqn:a_posteriori_bound} for $\const{cst:1}^2\coloneqq\max\{2+3\const{cst:3}^2, 6\}$.
\end{proof}
The equivalence $\muold \approx \min_{v\in V}\|v-u_h\|_h$ from \eqref{quasioptimalsmoother} provides 
	\begin{align*}%
		\|u-u_h\|_h^2 + \trb{e}^2\lesssim Res(e-J_hI_hI_\nc e)+\hatF(J_he_h-Qe_h)+\min_{v\in V}\|v-u_h\|_h^2
	\end{align*} as an equivalent
formulation of
\eqref{eqn:a_posteriori_bound}.
	The remaining parts of this paper discuss explicit bounds of the right-hand side of \eqref{eqn:a_posteriori_bound} for a simultaneous a
	posteriori analysis of five nonstandard FEMs for the biharmonic equation.

\section{Examples of {lowest}-order finite element schemes}%
\label{sec:Examples of second-order finite element schemes}
This section introduces the spaces and transfer operators for five lowest-order methods for the biharmonic equation.
\subsection{Three second-order finite element spaces}%
\label{sub:Three second-order finite element spaces}
Recall the space of piecewise polynomials $P_{\hspace{-.13em}k}(\T)$ of total degree at most $k\in \mathbb N$ from Subsection
\ref{sub:General notation}.
Let
$S^k(\T)\coloneqq P_{\hspace{-.13em}k}(\T)\cap C^0(\Omega)$ and 
$S^k_0(\T)\coloneqq\{p\in S^k(\T)\ |\ p_{|\partial\Omega}=0\}=P_{\hspace{-.13em}k}(\T)\cap H^1_0(\Omega).$
The associated $L^2$ projection $\Pi_k:L^2(\Omega)\to P_{\hspace{-.13em}k}(\T)$ is defined by
the $L^2$ orthogonality 
$(1-\Pi_{k})v\perp P_{\hspace{-.13em}k}(\T)$ for all $v\in L^2(\Omega)$.
Recall the nonconforming Morley space $\M(\T)$ from \eqref{eqn:Morley}.
Throughout the remaining parts of this paper on the biharmonic equation, specify $V_\nc\coloneqq \M(\T), V\coloneqq
H^2_0(\Omega)\subset \widehat V\coloneqq H^2(\T)$, and 
\begin{align} \label{eqn:V_h_def}
		V_h&\coloneqq \begin{cases}{}
			\M(\T)& \text{ for Morley},\\
				P_{\hspace{-.13em}2}(\T) &
			\text{ for dG or WOPSIP,} \\
				S^2_0(\T) &\text{ for $C^0$IP.} 
		\end{cases}
	\end{align}
\begin{center}
\begin{figure*}[t]
\begin{minipage}[h]{0.3\linewidth}
\begin{center}
\begin{tikzpicture}
\node[regular polygon, regular polygon sides=3, draw, minimum size=5cm]
(m) at (0,0) {};
\fill [black] (m.corner 1) circle (2pt);
\put(-3,-5){$T$}
\fill [black] (m.corner 2) circle (2pt);
\fill [black] (m.corner 3) circle (2pt);
\fill [black] (m.side 1) circle (2pt);
\fill [black] (m.side 2) circle (2pt);
\fill [black] (m.side 3) circle (2pt);
\phantom{\draw [-latex, thick] (m.side 2) -- ($(m.side 2)!0.5!90:(m.corner 2)$);}
\phantom{\draw[black] (m.corner 1) circle (.2cm);}
\end{tikzpicture}
\end{center}
\end{minipage}
\begin{minipage}[h]{0.3\linewidth}
\begin{center}
\begin{tikzpicture}
\node[regular polygon, regular polygon sides=3, draw, minimum size=5cm]
(m) at (0,0) {};
\fill [black] (m.corner 1) circle (2pt);
\put(-3,-5){$T$}
\fill [black] (m.corner 2) circle (2pt);
\fill [black] (m.corner 3) circle (2pt);
\draw [-latex, thick] (m.side 1) -- ($(m.side 1)!0.5!90:(m.corner 1)$);
\draw [-latex, thick] (m.side 2) -- ($(m.side 2)!0.5!90:(m.corner 2)$);
\draw [-latex, thick] (m.side 3) -- ($(m.side 3)!0.5!90:(m.corner 3)$);
\phantom{\draw[black] (m.corner 1) circle (.2cm);}
\end{tikzpicture}
\end{center}
\end{minipage}
\begin{minipage}[h]{0.3\linewidth}
\begin{center}
\begin{tikzpicture}
\node[regular polygon, regular polygon sides=3, draw, minimum size=5cm]
(m) at (0,0) {};
\fill [black] (m.corner 1) circle (2pt);
\draw[black] (m.corner 1) circle (.2cm);
\draw [thick] (0,0) -- (0,2.5);
\draw (0,0) -- (-2.1,-1.2);
\draw (0,0) -- (2.1,-1.2);
\draw (0,0) circle (2pt);
\put(1.2,0.6){${\rm mid}(T)$}

\fill [black] (m.corner 3) circle (2pt);
\draw[black] (m.corner 3) circle (.2cm);
\fill [black] (m.corner 2) circle (2pt);
\draw[black] (m.corner 2) circle (.2cm);
\draw [-latex, thick] (m.side 1) -- ($(m.side 1)!0.5!90:(m.corner 1)$);
\draw [-latex, thick] (m.side 2) -- ($(m.side 2)!0.5!90:(m.corner 2)$);
\draw [-latex, thick] (m.side 3) -- ($(m.side 3)!0.5!90:(m.corner 3)$);
\end{tikzpicture}
\end{center}
\end{minipage}
\caption{The Lagrange $P_{\hspace{-.13em}2}$, the Morley, and the HCT finite element (left to right)}\label{fig}
\end{figure*}
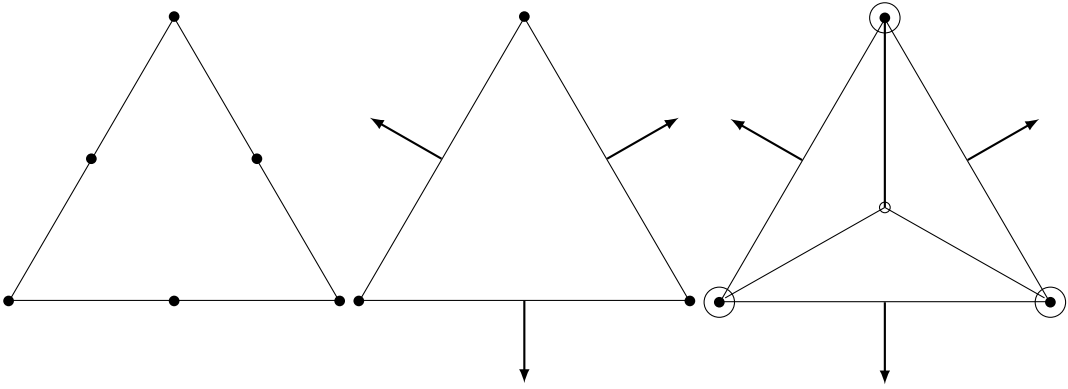
\end{center}
\subsection{Hilbert space of piecewise $H^2$ functions}%
\label{sub:Hilbert space of piecewise $H^2$ functions}
The semi-scalar product $a_\pw\coloneqq (D^2_\pw\bullet, D^2_\pw \bullet)_{L^2(\Omega)}$ in $\widehat V\coloneqq
H^2(\T)$ extends the energy scalar
product $a\coloneqq a_\pw|_{V\times V}$ and  the subspace $(\M(\T), a_\pw)$ is a Hilbert space.
Recall the jump $[v]_E$ resp.~normal jump $[\partial v/\partial\nu_E]_E$ across an edge $E\in\E$ of a piecewise function $v\in H^1(\T)$ resp.~$v\in
H^2(\T)$ from
Subsections \ref{ssub:Bound from interpolation} and \ref{ssub:Morley}.
Let $\V(E)$ denote the vertices of the edge $E\in\E$.
Define the semi-scalar product
$j_h:\widehat V\times \widehat V$, for any $v,w\in \widehat V$, by
\begin{align}\label{eqn:apw_jh}
	j_h(v, w)&\coloneqq
	\sum_{E \in \E}\left( \sum_{z \in {\mathcal V} (E)} \frac{[v]_E(z)}{h_E}\frac{[w]_E(z)}{h_E}
			 +
\fint_E  \jump{\frac{\partial v}{\partial\nu_E}}\!\!\!\mathrm ds\, \fint_E  \jump{\frac{\partial w}{\partial \nu_E}}\!\!\!\mathrm ds\right).
\end{align}
Since $j_h(v, \bullet)=0$ vanishes for any $v\in V+\M(\T)$,
$(H^2(\T), a_{\pw} + j_h)$ is a Hilbert space with the induced norm 
$\|\bullet\|_h$ from \eqref{eqn:pw_h}.
\begin{rem}[Completeness of $(\widehat V, a_\pw + j_h)$]\label{rem:Hilbert}
	It is clear \cite[Section 4.1]{carstensen_lowest-order_2022} that $(\widehat V, \|\bullet\|_h)$ is a normed linear
	space.
	Recall that $(H^2(\T), \|\bullet\|_{H^2(\T)})$ equipped with the piecewise $H^2$ norm
	$\|\bullet\|_{H^2(\T)}^2\coloneqq \sum^{}_{T\in\T} \|\bullet\|_{H^2(T)}^2$ is a Banach space.
	Let $Q:H^2(\T)\to P_{\hspace{-.13em}1}(\T)$ denote the $H^2$ orthogonal projection onto the finite dimensional space
	$P_{\hspace{-.13em}1}(\T)\subset
	H^2(\T)$ and set $X\coloneqq (1-Q) H^2(\T)$.
	The Bramble-Hilbert Lemma \cite[Lemma 11.9]{ern_finite_2021-1} asserts that $\trb{\bullet}_\pw$ is a norm on $X$ stronger
	than the piecewise $H^2$ norm $\|\bullet\|_{H^2(\T)}\lesssim \|\bullet\|_h$.
	Since \cite[Theorem 4.1]{carstensen_lowest-order_2022} shows that $\|\bullet\|_h\lesssim \|\bullet\|_{H^2(\T)}$ is
	also weaker than the piecewise $H^2$ norm, both norms are equivalent on $X$ and $X$ is complete.
	The direct sum $\widehat V=X \oplus P_{\hspace{-.13em}1}(\T)$ of two complete spaces is complete.
\end{rem}

\subsection{Classical and averaged Morley interpolation}%
\label{sub:Classical and averaged Morley interpolation}
The classical Morley interpolant from Subsection \ref{ssub:Morley} is defined for functions in $V + \M(\T)$ and has an
extension to piecewise $H^2$ functions.
Define the average $\langle\varphi\rangle_E\coloneqq\half\left(\varphi|_{T_+}+\varphi|_{T_-}\right)$ of $\varphi\in
H^1(\T)$ across an interior edge 
$E=\partial T_+\cap\partial T_-\in\E(\Omega)$ of the adjacent triangles  $T_+$ and $T_-\in\T$ as in Figure
\ref{fig:Edge_patch} and $\langle\varphi\rangle_E:=\varphi|_E$ 
along a boundary edge $E\in\E(\partial\Omega)$.
Let $\T(z)\coloneqq\{T\in\T\ |\ z\in T\}$ denote the $|\T(z)|\in\mathbb N$ many neighbouring triangles of $z\in T\in\T$.
\begin{defn}[Morley interpolation {\cite{carstensen_lowest-order_2022}}]\label{defccMorleyinterpolation}
	Given any piecewise function $v_{\pw} \in \widehat V$, the Morley interpolation operator $I_\M:\widehat V\to V_\nc$
	sets the degrees of freedom of the Morley finite element function $I_\M v_\pw\in V_\nc\coloneqq\M(\T)$ by
	\begin{align*}
	I_\M v_\pw(z)&\coloneqq 
 |\T(z)|^{-1}
\sum_{T \in \T(z)} (v_{\pw}|_{T})(z) &&\text{ for }z\in\mathcal{V}(\Omega), \\
	\fint_E\frac{\partial I_\M v_\pw}{\partial \nu_E} \ds &\coloneqq \fint_E \mean{\frac{\partial v_{\pw}}{\partial \nu_E}}
	\ds&&\text{ for
 }E\in\E(\Omega).
	\end{align*}
\end{defn} 
{ It is well known that there is a unique quadratic polynomial $I_\M v_\pw|_T\in P_{\hspace{-.13em}2}(T)$ that assumes the above values
$(I_\M v_\pw)(z)$ and $\fint_E\partial I_\M v_\pw/\partial\nu_E\ds$ at $z\in \mathcal{V}(T)$ and for all
$E\in\mathcal{E}(T)$.
Explicit formulas for the basis functions can be found in 
\cite{CCDGJH14}.}
This definition extends the classical Morley interpolation from Definition \ref{defccMorleyinterpolation_classic} to piecewise $H^2$ functions in
$\widehat V\equiv H^2(\T)$.
For any $v\in H^2_0(\Omega)+\M(\T)$, the $a_\pw$-orthogonality
\begin{align}\label{eqn:best_approx}
	a_{\pw}(v-I_\M v, w_2)  &=0 \text{ for all } w_2 \in P_{\hspace{-.13em}2}(\T)
\end{align}
 verifies the best-approximation property 
\begin{align} \label{eqn:best}
	\trb{v- I_\M v}_\pw = \min_{v_2\in P_{\hspace{-.13em}2}(\T)}\trb{v-v_2}_\pw.
\end{align}
This does not extend to discontinuous functions $v_h\in H^2(\T)$ in general. Recall $\|\bullet\|_h$ from \eqref{eqn:pw_h}.
\begin{theorem}[interpolation error {\cite[Thm.~4.3]{carstensen_lowest-order_2022}}]\label{thm:int_err}
Any piecewise smooth function $v_{\rm pw}\in H^2(\T)$ and its Morley interpolation $I_{\rm M} v_{\rm pw} \in \M(\T)$ from Definition \ref{defccMorleyinterpolation} satisfy
\begin{align*}
	(a)&&\|v_\pw - I_\M v_\pw\|_h
	   & \lesssim \| (1-\Pi_0)D^2_{\rm pw} v_{\rm pw}\|
+ j_h(v_{\rm pw}, v_{\pw})^{1/2},\\
	(b)&&
	\sum_{m=0}^2  h_\T^{m-2} | v_{\rm pw}- I_{\rm M} v_{\rm pw}|_{H^m(\T)}
	   &\approx\;\min_{w_\M\in \M(\T)}  \| v_{\rm pw}- w_\M \|_h\approx \|v_\pw - I_\M v_\pw\|_h.
\qed
\end{align*}
\end{theorem}
Since $I_\M\in L(\widehat V; V)$ is a bounded operator,
 the Cauchy inequality and the best-approximation property \eqref{eqn:best} verify \eqref{eqn:Inc_approx}
for $I_\nc\coloneq I_\M$ with $\Lambda_\nc\coloneqq\Lambda_\M\coloneqq 2 + \|I_\M\|_h$.
Indeed, for arbitrary $v_2\in P_{\hspace{-.13em}2}(\T)$ and $v\in V$, 
\begin{align}\label{eqn:Lambda_M}
	\|v_2 - I_\M v_2\|_h\leq \|v_2-v\|_h + \trb{v-I_\M v}_\pw + \|I_\M(v - v_2)\|_h\leq \Lambda_\M\|v-v_2\|_h.
\end{align}

\subsection{Transfer operator $I_h$}%
\label{sub:Transfer_operator_Morley}
The abstract setting from Section \ref{sec:abstractSeconddiscretisation} requires a transfer operator $I_h$ with
\eqref{eqn:Ih_approx} from
$\Vnc\coloneqq \M(\T)$ into $V_h$ defined in \eqref{eqn:V_h_def} for the different schemes.
The natural choice $I_h\coloneqq {\rm  id}$ for the Morley, dG, and WOPSIP method with $V_\nc\subseteq V_h$ fulfils \eqref{eqn:Ih_approx}
with $\Lambda_h=0$.
The situation is different for the $C^0$IP method with $V_\nc\not\subseteq V_h\coloneqq S^2_0(\T)$ and requires the
Lagrange interpolation
$I_h \coloneqq I_\C: \M(\T) \rightarrow S^2_0(\T)$ defined, for all $v_\M\in \M(\T)$, by
\begin{align} \label{eq:ic}
(I_\C v_\M)(z)=
\begin{cases}
	v_\M(z) &\text{for all } z \in \V,\\
	\mean{v_\M}(z)  &\text{for } z= \text{mid}(E), \; E \in \E(\Omega), \\
	0 &\text{for } z= \text{mid}(E), \; E \in \E(\partial \Omega).
\end{cases}
\end{align}
(It is well known that there exists a unique $I_\C v_\M|_T\in P_{\hspace{-.13em}2}(T)$ with prescribed values at the vertices and edge
midpoints from the unisolvence of the $P_{\hspace{-.13em}2}$ Lagrange finite element.)
Lemma 3.2 in \cite{CarstensenGallistlNataraj2015} establishes 
\eqref{eqn:Ih_approx} for the operator $I_h=I_\C$ with $\Lambda_h\approx 1$.

\subsection{Companion operator $J$}%
\label{sub:Companion operator $J$}
A conforming finite-dimensional subspace of $H^2_0(\Omega)$ is the Hsieh-Clough-Tocher  (HCT) \cite[Chap. 6]{Ciarlet} space 
$\mathrm{HCT}(\mathcal{T})\coloneqq\{v\in H^2_0(\Omega):\ v|_T\in P_{\hspace{-.13em}3}(\mathcal{K}(T))\text{ for all }T\in\mathcal{T}\}$
with the subtriangulation 
$\mathcal{K}(T):=\{\mathrm{conv}\{E, \mathrm{mid}(T)\}:\ E\in\mathcal{E}(T)\}$ 
of $T\in\T$ obtained by joining the vertices of $T$ with $\mathrm{mid}(T)$.
Figure \ref{fig} shows the degrees of freedom of the HCT finite element that extend those of the Morley element and
facilitate the design of a right-inverse to $I_\M:\widehat V\to \M(\T)$.

\begin{lemma}[right-inverse \cite{DG_Morley_Eigen,aCCP,veeser_zanotti1}]\label{companion} 
\label{lem:MorleyCompanion}
There exists a linear right-inverse $J: {\rm M}(\mathcal{T})\to HCT(\mathcal{T})+P_{\hspace{-.13em}8}(\mathcal{T}) \cap H^2_0(\Omega)$
for $I_\M:V\to\M(\T)$ 
and a constant   $\Lambda_\jc$, that exclusively depends on the shape regularity,
such that any $v_{\rm M} \in {\rm M}(\mathcal{T})$ satisfies
$$ \trinl v_\M- J v_\M \trinr_\pw \le \Lambda_\jc  \min_{v\in V}  \trinl v_\M- v \trinr_\pw. $$
\end{lemma}
See \cite[Section 5]{aCCP} for the definition of $J\in
L(V_\nc;V)$.
Note that Lemma \ref{lem:MorleyCompanion} verifies \eqref{eqn:J_approx} for the conforming companion $J$.
Recall from the previous subsections that $I_\nc\coloneqq I_\M\in L(\widehat V;V_\nc)$ and $I_h\in L(V_\nc; V_h)$ verify
\eqref{eqn:Inc_approx}--\eqref{eqn:Ih_approx}.
An immediate consequence of Lemma \ref{lem:J_h} is that $J_h\coloneqq J I_\M\in L(\widehat V; V)$
is a quasi-optimal smoother.
We refer to \cite{aCCP} for a 3D version.

\section[Building blocks for a posteriori estimators]{Building blocks for explicit residual-based a posteriori error estimators}%
\label{sec:Building block for explicit residual-based a posteriori error estimators}
This section establishes bounds on the error contributions in the right-hand side of \eqref{eqn:a_posteriori_bound}.
Recall the residual $Res\coloneqq F-a_\pw(u_h,\bullet)\in V^*$ from Section \ref{sec:Abstract a posteriori error
analysis} and set $V_\nc\coloneqq \M(\T)$ with interpolation operator $I_\nc\equiv I_\M$  and quasi-optimal smoother
$J_h\equiv JI_\M$ throughout the remaining parts
of this paper.
\subsection{Estimates for $1-\JIM I_h I_\M$ and $(1-\JIM) I_h I_\M$}%
\label{sub:Estimates for $id-JI_MI_hI_M$}
The linear operators $1-\JIM I_h I_\M:V\to V$ and $I_h I_\M-\JIM I_h I_\M:V\to \widehat V$ are stable in the energy norm.
\begin{lemma}[stability]\label{lem:5.1}
	Any $v\in V=H^2_0(\Omega)$ with $\widehat w= (1-\JIM I_h I_\M) v\in V$ or $\widehat w= (1-\JIM)I_hI_\M v\in \widehat V$ satisfies
	\begin{align*}
		\sum^{2}_{m=0} |h_\T^{m-2}\widehat w|_{H^m(\Omega)}^2+\sum^{}_{E\in\E(\Omega)}\left(
		\|h_E^{-3/2}\widehat w\|_{L^2(E)}^2+\|h_E^{-1/2}\nabla \widehat w\|_{L^2(E)}^2\right)\leq \const{cst:5}^2 \trb{v}^2.
	\end{align*}
\end{lemma}
\begin{proof}
	Since $J$ is a right-inverse of $I_\M$, the functions $v, v_\M\coloneqq I_\M v, v_h\coloneqq I_h v_\M, I_\M v_h$ and $\JIM v_h$ in $H^2(\T)$ are continuous
	at any vertex $z\in\mathcal{V}$ and coincide at $z\in\mathcal{V}$.
	Hence, $\widehat w|_{T}\in H^2(T)$ vanishes at the three vertices of the triangle $T\in\T$.
	It is textbook analysis \cite{Brenner, Ciarlet, Braess, ern_finite_2021-1} to derive the bounds
	\begin{align}\label{eqn:C_BH}
		\sum^{2}_{m=0} |h_\T^{m-2}\widehat w|_{H^m(T)}^2\leq C_{\rm BH}^2 |\widehat w|_{H^2(T)}^2
	\end{align}
	from an application of the Bramble-Hilbert lemma with a constant $C_{\rm BH}>0$
	and we refer to \cite[Sec.\ 3]{CarstensenGedickeRim}
	for explicit constants in terms of the maximal angles in the triangle $T\in \T$.
	The sum of all those estimates \eqref{eqn:C_BH} results in 
	\begin{align*}
		\sum^{2}_{m=0} |h_\T^{m-2}\widehat w|^2_{H^m(\Omega)}\leq C_{\rm BH}^2\trb{\widehat w}_\pw^2.
	\end{align*}
	The previous estimate, %
	$\trb{\widehat w}_\pw\leq\|\widehat w\|_h\leq \const{cst:w_h_to_v}\trb{v}$ with
	$
		\const{cst:w_h_to_v}\coloneqq\max\{1, C_\jc\}(1+\Lambda_\jc)(1+\Lambda_\nc)(1+\Lambda_h)
	$ from Lemma
	\ref{lem:IIJ_approx}--\ref{lem:J_h} and \eqref{eqn:best} conclude the proof of 
	\begin{align}
		\sum^{2}_{m=0} |h_\T^{m-2}\widehat w|_{H^m(\Omega)}^2\leq C_{\rm BH}^2\const{cst:w_h_to_v}^2\trb{v}^2.
	\end{align}
	Given any interior edge $E\in\E(\Omega)$ with adjacent triangle $T(E)\in\T$, the trace inequality 
\cite[Eqn.\ (12.17)]{ern_finite_2021-1}
provides a constant $C_{\rm tr}>0$ exclusively depending on the shape-regularity with
\begin{align*}%
	h_E^{-3/2}\|\widehat w\|_{L^2(E)}+h_E^{-1/2}\|\nabla \widehat w\|_{L^2(E)}\leq C_{\rm tr}\left(h_T^{-2}\|\widehat w\|_{L^2(T(E))} +
		h_T^{-1}|\widehat w|_{H^1(T(E))}
		+|\widehat w|_{H^2(T(E))}\right).
	\end{align*}
	This and the sum over the interior edges $E\in\E(\Omega)$ result in
	\begin{align*}
		\sum^{}_{E\in\E(\Omega)} \left(h_E^{-3/2}\|\widehat w\|_{L^2(E)}+h_E^{-1/2}\|\nabla \widehat w\|_{L^2(E)}\right)^2\leq
		3C_{\rm tr}^2\sum^{2}_{m=0}\sum^{}_{E\in\E(\Omega)} |h_T^{m-2}\widehat w|_{H^m(T(E))}^2\W{.}
	\end{align*}
	Since every triangle $T(E)\in\T$ is counted at most $3$ times (once for every edge $E\in \E(T(E))$) in the last sum,
	the claim follows with $\const{cst:5}\coloneqq (3 C_{\rm tr}+1)C_{\rm BH}\const{cst:w_h_to_v}$.
\end{proof}

\begin{cor}[bound for $F\in L^2(\Omega)$]\label{cor:F_bound}
	Any $F=f\in L^2(\Omega)$ and $v\in V=H^2_0(\Omega)$ with $\widehat w\coloneqq (1-\JIM I_hI_\M) v\in V$ or
	$\widehat w\coloneqq(1-\JIM)I_hI_\M v\in \widehat V$ satisfy
	\begin{align}
		\int_\Omega f\,\widehat w\;\mathrm dx\leq \const{cst:4}\|h_\T^2 f\|\trb{v}.
	\end{align}
\end{cor}
\begin{proof}
	This follows from Lemma \ref{lem:5.1} and a Cauchy inequality in $L^2(\Omega)$ in
	\begin{align*}
		\int_\Omega f \, \widehat w\;\mathrm dx\leq \|h_\T^2 f\|\|h_\T^{-2}\widehat w\|\leq \const{cst:4}\|h_\T^2 f\|\trb{v}.\qquad\qedhere
	\end{align*}
\end{proof}
\noindent Define the oscillations of $f\in L^2(\Omega)$ by $\mathrm{osc}_2(f, T)\coloneqq\|h_\T^2(f-\Pi_2 f)\|_{L^2(T)}$ and abbreviate
\begin{align*}\mathrm{osc}_2(f, \mathcal{S})\coloneqq \sqrt{\sum^{}_{T\in \mathcal{S}} \mathrm{osc}_2^2(f, T)}\end{align*} for a subset 
$\mathcal{S}\subseteq \T$ of triangles in $\T$.
The efficiency of the term $\|h_\T^2 f\|_{L^2(T)}$ is known, e.g., from \cite[Lemma 4.2 \& Remark 4.4]{BrenGudiSung10}; 
Section \ref{sec:General sources} treats a more general source $F\in V^*$.
\begin{lemma}[efficiency up to oscillations {\cite{BrenGudiSung10}}]\label{lem:f_efficiency}
	Let $u\in V$ be the weak solution to \eqref{eqn:WP} for a right-hand side $F=f\in L^2(\Omega)$. 
	Then
	$\|h_\T^2 f\|_{L^2(T)}\lesssim | u-I_\M u|_{H^2(T)} + \mathrm{osc}_2(f, T)$.\qed
\end{lemma}

\subsection{Error estimates for $a_\pw(v_h, w)$}%
\label{sub:Error estimates}
Recall the abbreviation $w\coloneqq v-\JIM I_h I_\M v$ for $v\in V$.
Since $J$ from Subsection \ref{sub:Companion operator $J$} is a right-inverse of the Morley interpolation $I_\M$ from
Subsection \ref{sub:Classical and averaged Morley interpolation}, the key observation for the situation $I_h=\id$ is
\begin{align}\label{eqn:IM_w_0}
	I_\M w=I_\M v-I_\M \JIM I_h I_\M v=I_\M(v-I_h I_\M v) = 0.
\end{align}
This is the case for the Morley, dG, and WOPSIP methods and, hence, the $a$-orthogonality of the Morley interpolation of
$w\in V$ and $I_\M w=0$ imply $a_\pw({ u_h}, w)=0$. 
For the $C^0$IP method with $V_h=S^2_0(\T)$ and $I_h=I_\C \ne \id$ from Subsection \ref{sub:Transfer_operator_Morley} the situation
differs and is the content of the remaining part of this subsection.
\begin{lemma}[bound for $a_\pw(v_h, w)$]\label{lem:C0IP_a}
	Any $v_h\in V_h$ and $v\in V$ with 
	$w\coloneqq v-\JIM I_h I_\M v$ satisfies
	\begin{align*}
		|a_\pw(v_h, w)|\leq \begin{cases}{}
			0&\text{if }I_h=\id,\\
			\const{cst:5}\sqrt{\sum^{}_{E\in\E(\Omega)} h_E\|[\partial_{\nu\nu}^2v_h]_E\|_{L^2(E)}^2}\trb{v}&\text{if
			}I_h=I_\C.
		\end{cases}
	\end{align*}
\end{lemma}
\begin{proof}
	With the remark succeeding \eqref{eqn:IM_w_0},  $(i)$ holds and it remains to prove $(ii)$. %
Since the piecewise Hessian $D^2_\pw v_h$ of $v_h\in S^2_0(\T)$ is piecewise constant, no volume contributions arise in
a piecewise integration by parts with the conforming test function $w\in V$.
A careful re-arrangement of the contributions along the boundary $\partial T$ of $T\in\T$ reveals
\begin{align}\label{eqn:a_pw_ibp}
	a_\pw(v_h, w)=\sum^{}_{E\in\E(\Omega)} \int_E\nabla w\cdot[D^2_\pw v_h]_E\nu_E\mathrm ds.
\end{align}
Recall from the proof of Lemma \ref{lem:5.1} that $w(z)=0$ vanishes at any vertex, whence
$\int_E\partial w/\partial s \: \mathrm ds=0$ on any edge $E\in \E$. Since the matrix $[D^2_\pw v_h]_E\in
P_{\hspace{-.13em}0}(E; \mathbb S)$ is constant, the split $\nabla w =
(\partial w/\partial s)\tau_E + (\partial w/\partial\nu_E) \nu_E$ along $E\in\E$ and the Cauchy inequality show
\begin{align*}
	\int_E\nabla w\cdot [D^2_\pw v_h]_E\nu_E\mathrm ds = \int_E\frac{\partial w}{\partial\nu_E}[\partial_{\nu\nu}^2
	v_h]_E\mathrm ds\leq h_E^{-1/2}\left\|\frac{\partial w}{\partial\nu_E}\right\|_{L^2(E)}h_E^{1/2}\left\|[\partial_{\nu\nu}^2
	v_h]_E\right\|_{L^2(E)}.
\end{align*}
Notice that the trace of $\nabla w\cdot \nu_E$ along $E$ is continuous for $w\in V$.
This, a Cauchy inequality in $\ell^2$, and $\|\partial w/\partial\nu_E\|_{L^2(E)}\leq\|\nabla
w\|_{L^2(E)}$ verify
\begin{align*}
	a_\pw(v_h, w)\leq
	\sqrt{\sum^{}_{E\in\E(\Omega)}
	h_E\|[\partial_{\nu\nu}^2v_h]_E\|_{L^2(E)}^2}\sqrt{\sum^{}_{E\in\E(\Omega)} h_E^{-1}\left\|\nabla
w\right\|_{L^2(E)}^2}.
\end{align*}
This and Lemma \ref{lem:5.1} conclude the
proof.
\end{proof}
The efficiency estimate of the jump contributions in Lemma \ref{lem:C0IP_a} is known, e.g., from the $C^0$IP method
\cite{BGS10}.
For any edge $E\in \E$, the sub-triangulation $\T(\omega(E))\coloneqq\{T\in\T\ |\ E\subset \partial T\}$ in the
edge-patch $\omega(E)\coloneqq \mathrm{int}(T_+\cup T_-)$ consists of one or two triangles.
\begin{lemma}[{\cite[Lemma 4.3]{BGS10}}]\label{lem:C0IP_eficiency}
	Let $u\in V$ solve \eqref{eqn:WP} for $F=f\in
L^2(\Omega)$.
Any $v_h\in P_{\hspace{-.13em}2}(\T)$ and any edge $E\in\E$ satisfy 
\begin{align*}
	h_E^{1/2}\|[\partial_{\nu\nu}^2v_h]_E\|_{L^2(E)}\lesssim |u-v_h|_{H^2(\T(\omega(E)))}  + \mathrm{osc}_2(f,
	\T(\omega(E))).
\end{align*}
\end{lemma}
\begin{proof}
	The proof of \cite[Lemma 4.3]{BGS10} for the jump $[\partial_{\nu\nu}^2v_h]_E$ of any $v_h\in P_{\hspace{-.13em}2}(\T)$
	shows that
	\begin{align*}
		h_E^{1/2}\|[\partial_{\nu\nu}^2v_h]_E\|_{L^2(E)}\lesssim |u-v_h|_{H^2(\T(\omega(E)))} + 
		\|h_\T^2f\|_{L^2(\omega(E))}.
	\end{align*}
	Lemma \ref{lem:f_efficiency} and $|u-I_\M u|_{H^2(T)}=\min_{v_h\in P_{\hspace{-.13em}2}(T)}|u-v_h|_{H^2(T)}$ as in \eqref{eqn:best} conclude the proof.
\end{proof}

\subsection{Estimate of $\|v_h - J_h v_h\|_h$}%
\label{sub:Estimate v_h}
This subsection discusses reliable and efficient bounds of $\|v_h - J_h v_h\|_h$ in terms of two different jump
terms that appear in the a posteriori analysis, e.g., in \cite{DVNJSR2007, HuShi09, BGS10, BrenGudiSung10}.
\begin{theorem}[reliability and efficiency of $\|v_h - J_h v_h\|_h$] \label{thm:mu_bound}
	Any $v_h\in V_h$ satisfies
	\begin{align*}
\min_{v\in V}\|v-v_h\|_h^2\approx	\|v_h - J_h v_h\|_h^2&\approx \sum_{E\in\E} h_E\|[D^2_\pw v_h]_E\tau_E\|_{L^2(E)}^2 + j_h(v_h,
		v_h)\\
							   &\approx\sum_{E\in\E}\left( h_E^{-3}\left\|\left[v_h\right]_E\right\|_{L^2(E)}^2 + h_E^{-1}
	\left\|\left[\frac{\partial v_h}{\partial \nu_E}\right]_E\right\|_{L^2(E)}^2\right).
	\end{align*}
	
\end{theorem}
\noindent The remaining parts of this subsection are devoted to the proof and depart with the following generalization
of \cite[Thm.~2.1]{HuShi09}.
\begin{lemma}[bound for $\|v_h-J_hv_h\|_h$]\label{lem:mu_bound}
	Any $v_h\in V_h$ satisfies
	\begin{align*}
		\const{cst:5}^{-1}\|v_h - J_h v_h\|_h^2\leq \sum_{E\in\E} h_E\|[D^2_\pw v_h]_E\tau_E\|_{L^2(E)}^2 + j_h(v_h,
		v_h).
	\end{align*}
\end{lemma}
\begin{proof}[Proof of Lemma \ref{lem:mu_bound}]
	Given any $v_h\in V_h$, set $v_\M\coloneqq I_\M v_h\in \M(\T)$.
	A triangle inequality and \eqref{eqn:pw_h} verify $\|v_h-J_hv_h\|_h\leq\|v_h-v_\M\|_h+\trb{v_\M-Jv_\M}_\pw$.
	It follows from \cite[Lem.~5.1]{aCCP} that %
	\begin{align*}
		\trb{v_\M-Jv_\M}_\pw^2&\lesssim \sum^{}_{E\in\E} h_E\|[D^2_\pw v_\M]_E\tau_E\|_{L^2(E)}^2.
	\end{align*}
	This, a triangle inequality, and the discrete trace inequality $h_E^{1/2}\|D^2_\pw (v_h-v_\M)\|_{L^2(E)}\lesssim
\|D^2(v_h-v_\M)\|_{L^2(T)}$
from \cite[Lem.~12.8]{ern_finite_2021-1} result in 
	\begin{align*}
		\trb{v_\M-Jv_\M}_\pw^2&\lesssim\sum^{}_{E\in\E}
		h_E\|[D^2_\pw v_h]_E\tau_E\|_{L^2(E)}^2+\trb{v_h-v_\M}_\pw^2.
	\end{align*}
	This and $\trb{v_h-v_\M}_\pw\leq\|v_h-v_\M\|_h\lesssim j_h(v_h, v_h)^{1/2}$ from \eqref{eqn:pw_h} and Theorem
	\ref{thm:int_err}.a with $D^2_\pw v_h\in P_{\hspace{-.13em}0}(\T)$ conclude the proof.
\end{proof}
The inverse inequality leads to an alternative upper bound in Lemma \ref{lem:mu_bound}.
\begin{lemma}[alternative bound]\label{lem:BrennerBound}
	Any $v_h\in V_h$ and any edge $E\in \E$ satisfy
	\begin{align*}
		h_E\|[D^2_\pw v_h]_E\tau_E\|_{L^2(E)}^2 +& \sum_{z\in \mathcal{V}(E)}\frac{|[v_h]_E(z)|^2}{h_E^2} + \left|\fint_E
		\jump{\frac{\partial v_h}{\partial\nu_E}}\!\!\!\mathrm ds\right|^2\\
												 &\leq \const{cst:6}\left(h_E^{-3}\|[v_h]_E\|_{L^2(E)}^2 + h_E^{-1}\|[\partial
												 v_h/\partial\nu_E]_E\|_{L^2(E)}^2\right).
	\end{align*}
\end{lemma}
\begin{proof}
	The split $D^2 v_h \cdot \tau_E = (\partial^2 v_h/\partial s\partial s) \tau_E + (\partial^2 v_h/\partial
	s\partial \nu_E) \nu_E$, the Cauchy inequality, and the linearity of the jump show
	\begin{align*}
		\|[D^2_\pw v_h]_E\cdot \tau_E\|_{L^2(E)}\leq \left\|\frac{\partial^2}{\partial s\partial
			s}\left[v_h\right]_E\right\|_{L^2(E)} +
		\left\|\frac{\partial}{\partial s}\left[\frac{\partial v_h}{\partial \nu_E}\right]_E\right\|_{L^2(E)}.
	\end{align*}
	The inverse inequality \cite[Lemma 12.1]{ern_finite_2021-1} states the existence of a constant $C_{\rm inv}>0$ with
	$|p|_{H^m(E)}\leq C_{\rm inv} h_E^ {-m}\|p\|_{L^2(E)}$ and $\|p\|_{L^p(E)}\leq C_{\rm inv}h_E^{1/p-1/q}\|p\|_{L^q(E)}$ for any
	$p\in P_{\hspace{-.13em}2}(E)$ and $m\in \mathbb N_0, 1\leq p,q\leq\infty$.
	Since $[v_h]_E$ and $[\partial v_h/\partial\nu_E]_E$ are quadratic polynomials on $E$, this shows
	\begin{align*}
		h_E^{1/2}\|[D^2_\pw v_h]_E\cdot \tau_E\|_{L^2(E)}&\leq C_{\rm inv}\left(h_E^{-3/2}\left\|\left[v_h\right]_E\right\|_{L^2(E)} +
	h_E^{-1/2}
	\|[\partial v_h/\partial\nu_E]_E\|_{L^2(E)}\right),\\
			\sum^{}_{z\in\mathcal{V}(E)} \frac{|[v_h]_E(z)|}{h_E}&\leq 2h_E^{-1} \|[v_h]_E\|_{L^\infty(E)}\leq 2C_{\rm inv}h_E^{-3/2}\left\|\left[v_h\right]_E\right\|_{L^2(E)},\\
		\left|\fint_E  \jump{\frac{\partial v_h}{\partial\nu_E}}\!\!\!\mathrm ds\right| &\leq h_E^{-1}\|[\partial
		v_h/\partial\nu_E]_E\|_{L^1(E)}\leq C_{\rm inv}h_E^{-1/2}\|[\partial
		v_h/\partial\nu_E]_E\|_{L^2(E)}.
	\end{align*}
	The sum of these terms squared and the Cauchy inequality $(A+B)^2 \leq 2 A^2 + 2B^2$ for $A, B\in\mathbb R$  conclude
	the proof with $\const{cst:6}\coloneqq 6C_{\rm inv}^2$.
\end{proof}
\begin{proof}[Proof of Theorem \ref{thm:mu_bound}]
	The reliability of the first estimator follows from Lemma \ref{lem:mu_bound}. This and Lemma \ref{lem:BrennerBound}
	provide 
	\begin{align*}
		\|v_h - \JIM v_h\|_h^2&\lesssim \sum_{E\in\E} h_E\|[D^2_\pw v_h]_E\tau_E\|_{L^2(E)}^2 + j_h(v_h,
		v_h)\\
							   &\lesssim\sum_{E\in\E}\left( h_E^{-3}\left\|\left[v_h\right]_E\right\|_{L^2(E)}^2 + h_E^{-1}
	\left\|\left[\frac{\partial v_h}{\partial \nu_E}\right]_E\right\|_{L^2(E)}^2\right).%
	\end{align*}
	Since the jumps $[\JIM v_h]_E$ and $[\partial \JIM v_h/\partial\nu_E]_E$ vanish for a conforming function 
	$\JIM v_h\in V$ on any edge $E\in\E$ and $\trb{v_h-\JIM v_h}_\pw\geq0$, the last term is bounded by
	\begin{align*}
		\trb{v_h - \JIM v_h}_\pw^2&+ \sum_{E\in\E}\left( h_E^{-3}\left\|\left[v_h - \JIM v_h\right]_E\right\|_{L^2(E)}^2 + h_E^{-1}
		\left\|\left[\frac{\partial (v_h - \JIM v_h)}{\partial \nu_E}\right]_E\right\|_{L^2(E)}^2\right)\\
		&\approx \|v_h - \JIM v_h\|_h^2
	\end{align*}
	with the equivalence of norms in $V+P_{\hspace{-.13em}2}(\T)$ from \cite[Thm.\ 4.1]{CarstensenGallistlNataraj2015} in the last step.
	This proves the equivalence of both estimators to $\|v_h - \JIM v_h\|_h\leq C_{\rm J}\min_{v\in
	V}\|v-v_h\|_h$ by the quasi-optimality \eqref{quasioptimalsmoother} of $\JIM$.
	The trivial estimate $\min_{v\in
	V}\|v-v_h\|_h\leq \|v_h - \JIM v_h\|_h$ concludes the proof.
\end{proof}
\section{Unified a posteriori error control}
\label{sec:Unified a posteriori error control}
This section reconsiders the biharmonic equation \eqref{eqn:WP2} with weak solution $u\in V\coloneqq H^2_0(\Omega)$ and
the discrete solution $u_h\in V_h$ of the Morley, dG, $C^0$IP, and WOPSIP schemes defined in Subsections \ref{sub:Morley
method}--\ref{sub:WOPSIP method}.
The presentation unifies the a posteriori error analysis of the well-known discretization schemes with orignal and modified  right-hand side.

\subsection{Discretisation of the biharmonic equation}%
\label{sub:Discretisation of the biharmonic equation}
Recall that $\widehat V\equiv H^2(\T)$ is a Hilbert space with scalar product $a_\pw + j_h$. 
Recall the discrete spaces $V_\nc\coloneqq \M(\T)$
and $V_h$ from Section \ref{sec:Examples of second-order finite element schemes}.
The weak solution $u\in V\coloneqq H^2_0(\Omega)$ to the biharmonic equation $\Delta^2 u = F\in V^*$ solves \eqref{eqn:WP}
with the energy scalar product $a\coloneqq a_\pw|_{V\times V}$ on $V$ and $a_\pw:\widehat V\times \widehat V\to\mathbb
R$ given in Subsection \ref{sub:Hilbert space of piecewise $H^2$ functions}.

Recall $J_h\coloneqq JI_\M\in L(\widehat V;V)$ from Section \ref{sec:Examples of second-order finite element schemes}.
Each method defines its particular discrete bilinear form $a_h:\left(V_h + \M(\T)\right)\times \left(V_h + \M(\T)\right)\to \mathbb R$ in the subsequent
subsections. 
The discrete solution $u_h\in V_h$ solves
\begin{align}\label{eqn:DWP_Vh}
	a_h(u_h, v_h) = (f, Q v_h)_{L^2(\Omega)}&\text{ for all }v_h\in V_h
\end{align}
with $f\in L^2(\Omega)$ and $Q\in\{\id, J_h\}$ in this section. %
The discrete problem \eqref{eqn:DWP_Vh} is a rewriting of \eqref{eq:discrete2} for $\hatF\coloneqq F\equiv f\in L^2(\Omega)$ without smoother $Q\coloneqq\id$ or with
 the quasi-optimal (by Lemma
\ref{lem:J_h}) smoother $Q\coloneqq J_h$.
Section \ref{sec:General sources}
discusses more general right-hand sides $F\in V^*$ with a natural extension $\hatF\in
H^2(\T)^*$.
The key assumption \eqref{eqn:H} from \cite{carstensen_lowest-order_2022} holds for the Morley, 
dG, and $C^0$IP discretisations. %
Hence, the a priori estimate from Theorem \ref{thm:abstract_main} holds for these methods and leads to the quasi-best approximation property 
\begin{align}
	\|u-u_h\|_h\lesssim \trb{u-I_\M u}_\pw =\min_{v_2\in P_{\hspace{-.13em}2}(\T)}\trb{u-v_2}_\pw.\label{eqn:quasi_best}
\end{align}
In particular, this shows equivalence of these methods from an a priori point of view.

\subsection{Morley FEM}%
\label{sub:Morley method}
The Morley FEM for the biharmonic equation \eqref{eqn:DWP_Vh} comes with $a_h\coloneqq a_\pw$.
The subsequent result recovers the equivalent a posteriori estimates from \cite[Thm.\ 2.2]{HuShi09} and \cite[Eqn.\ (3.2)]{DVNJSR2007}.
\begin{theorem}[a posteriori estimate]\label{thm:Morley_a_posteriori}
	The discrete Morley solution $u_h\in V_h$ to \eqref{eqn:DWP_Vh} and the exact solution $u\in V$ to \eqref{eqn:WP}
	with source $f\in L^2(\Omega)$ satisfy
	\begin{align*}
		\trb{u-u_h}_\pw^2+ \mathrm{osc}_2^2(f)&\approx \|h_\T^2 f\|_{}^2 + \sum_{E\in \E}h_E\|[D_\pw^2
		u_h]_E\tau_E\|_{L^2(E)}^2\\
		&\approx \|h_\T^2 f\|_{}^2 + \sum_{E\in\E}\left( h_E^{-3}\left\|\left[u_h\right]_E\right\|_{L^2(E)}^2 + h_E^{-1}
	\left\|\left[\frac{\partial u_h}{\partial \nu_E}\right]_E\right\|_{L^2(E)}^2\right).%
	\end{align*}
	The equivalence constants exclusively depend on the shape-regularity of $\T$.
\end{theorem}
\begin{proof}
	Set $\W{w\coloneqq e-\JIM e_h\in V}, e_h\coloneqq I_\M e\in \M(\T)$ and $\widehat w\coloneqq (1-J_h)e_h\in V+\M(\T)$ for $e\coloneqq u-\JIM u_h$ and recall $I_h=\id$ from Subsection
	\ref{sub:Transfer_operator_Morley}.
	Since $a_\pw(u_h, w) = 0$ from Lemma \ref{lem:C0IP_a}, the definition of the residual and
	$F(w)\lesssim\|h_\T^{2}f\|\trb{e}$ from Corollary \ref{cor:F_bound} show
	\begin{align*}
		Res(w) \coloneqq F(w) - a_\pw(u_h, w)=F(w)\lesssim \|h_\T^2f\|\trb{e}.
	\end{align*}
	Since $(f, \JIM e_h-Qe_h)_{L^2(\Omega)}=0$ vanishes for $Q=\JIM$, Corollary \ref{cor:F_bound} provides
	\begin{align}\label{eqn:fQ_bound}
		\int_\Omega f(\JIM e_h - Qe_h)\dx\lesssim \|h_\T^2 f\|\trb{e}
	\end{align}
	for $Q=\id$ and $Q=\JIM$.
	The two previously displayed estimates and \eqref{eqn:u_u_h_mu} verify \begin{align*}Res(w) + \int_\Omega f(\JIM
	e_h - Qe_h)\dx \lesssim \|h_\T^2 f\|(\trb{u-u_h}_\pw + \|u_h-\JIM u_h\|_h).\end{align*}
	The stability of the $L^2$ projection shows $\mathrm{osc}_2(f)\leq \|h_\T^2f\|$.
	Hence 
	Theorem \ref{thm:a_posteriori_bound} plus a weighted Young inequality result in
	\begin{align}\label{eqn:Morley_proof}
		\trb{u-u_h}_\pw^2 + \mathrm{osc}_2^2(f) \lesssim \|h_\T^2 f\|^2+ \|u_h-\JIM u_h\|_h^2.
	\end{align}
	Since $j_h(u_h,u_h) = 0$ for
	$u_h\in \M(\T)$, Theorem \ref{thm:mu_bound}
	bounds $\|u_h - \JIM u_h\|_h^2$ in \eqref{eqn:Morley_proof} by either of the jump terms.
	This proves the reliability for both estimators.
	The efficiency of $\|h_\T^2f\|$ follows from Lemma \ref{lem:f_efficiency} while Theorem \ref{thm:mu_bound}
	verifies the efficiency for all jump terms.
\end{proof}

\subsection{Discontinuous Galerkin 1} \label{sec:DGFEM1}
Recall the definition of the jump $\jump{\bullet}$ and average $\mean{\bullet}$ (applied componentwise to
matrix-valued functions) along an edge $E\in\E$ from Subsection \ref{sub:Evolution of a posteriori error analysis} and
\ref{sub:Classical and averaged Morley interpolation}.
\noindent   The bilinear form \begin{align}\label{eqn:ah_dg1}a_h(\bullet,\bullet)=a_\pw(\bullet,\bullet) +
b_h(\bullet,\bullet) + c_{\dg}(\bullet,\bullet)\end{align} for the discontinuous Galerkin method (dG)  \cite{Baker:1977:DG,
FengKarakashian:2007:CahnHilliard} depends on $-1\leq\Theta\leq 1$ and parameters $\sigma_1, \sigma_2>0$. 
For every $v_2, w_2\in P_{\hspace{-.13em}2}(\T)\supset V_\nc + V_h$,
\begin{subequations}
\begin{align}
   b_h(v_2, w_2)& \coloneqq -\Theta {\cal J}(v_2, w_2) -   {\cal J}(w_2, v_2), \label{eq:bhindG} \\
     {\cal J}(v_2, w_2) & \coloneqq  \sum_{E\in\E} 
	 \int_E \jump{\nabla_\NC v_2}  \cdot \mean{D^2_\NC w_2}\nu_E  \ds,\\\label{eqn:bhindGb}
	 c_\dg(v_{2},w_{2}) &\coloneqq\sum_{E\in\E}\left(
    \frac{\sigma_1}{h_E^3} \int_E \jump{v_2}\jump{w_2} \ds      
     +\frac{\sigma_2}{h_E} 
	 \int_E \jump{\frac{\partial v_2}{\partial \nu_E}}\jump{\frac{\partial w_2}{\partial \nu_E}} \ds\right).
\end{align}
\end{subequations}
This is the symmetric (resp.\ non symmetric) interior penalty Galerkin formulation for $\Theta=1$ (resp. $\Theta=-1$).
An appropriate choice \cite{FengKarakashian:2007:CahnHilliard,MozoSuli07} of the
parameters $\sigma_1,\sigma_2$ guarantees $V_h$-ellipticity \eqref{eq:a_ellipticity}.
Throughout this paper, \eqref{eq:a_ellipticity} is assumed for $\sigma_1=\sigma_2\approx 1$.
The following theorem recovers the known a posteriori error estimator from \cite{CCGMNN18} for the linear part.  

\begin{theorem}[a posteriori estimate]\label{thm:dG_a_posteriori}
	The discrete dG solution $u_h\in V_h$ to \eqref{eqn:DWP_Vh} with $a_h$ from \eqref{eqn:ah_dg1}  and the exact
	solution $u\in V$ to \eqref{eqn:WP} with $f\in
	L^2(\Omega)$ satisfy
	\begin{align*}
		\|u-u_h\|_h^2 + \mathrm{osc}_2^2(f)&\approx \|h_\T^2 f\|_{}^2 + \sum_{E\in \E}h_E\|[D_\pw^2 u_h]_E\tau_E\|_{L^2(E)}^2 +
		j_h(u_h, u_h)\\
		&\approx\|h_\T^2 f\|_{}^2 + \sum_{E\in\E}\left( h_E^{-3}\left\|\left[u_h\right]_E\right\|_{L^2(E)}^2 + h_E^{-1}
		\left\|\left[\frac{\partial u_h}{\partial \nu_E}\right]_E\right\|_{L^2(E)}^2\right).
	\end{align*}
	The
	equivalence constants exclusively depend on the shape-regularity of $\T$.
\end{theorem}
\begin{proof}
	Recall $I_h=\id$ from Subsection \ref{sub:Transfer_operator_Morley} and, thus, 
	the proof of the reliability and efficiency follows the proof of Theorem \ref{thm:Morley_a_posteriori} verbatim
	except for $j_h(u_h, u_h)\ne0$ in general.
	The additional term $j_h(u_h, u_h)$ from the reliability estimate of $\|u_h-J_hu_h\|_h$ in Lemma \ref{lem:mu_bound} enters
	the right-hand side of the first estimator.
	Since  $\|u-u_h\|_h^2$ bounds the efficient jump terms $j_h(u_h, u_h) = j_h(u-u_h, u-u_h)$  by definition in \eqref{eqn:pw_h},
	this concludes the proof.
\end{proof}
\begin{cor}\label{cor:dG_corollary}
	The discrete dG solution $u_h\in V_h$ to \eqref{eqn:DWP_Vh} with $a_h$ from \eqref{eqn:ah_dg1} and the exact solution $u\in V$ to \eqref{eqn:WP} and $f\in
	L^2(\Omega)$ satisfy
	\begin{align*}
		\trb{u-u_h}_\pw^2 + c_\dg(u_h, u_h)+  \mathrm{osc}_2^2(f)&\approx \|h_\T^2 f\|_{}^2 + c_\dg(u_h, u_h).
	\end{align*}
\end{cor}
\begin{proof}
	Since $\sigma_1=\sigma_2\approx1$, the jump contributions in the second estimator in Theorem
	\ref{thm:dG_a_posteriori} are equivalent to $c_\dg(u_h, u_h)$.
	Because $c_\dg(v, \bullet) = 0$ vanishes for any $v\in V$, the statement follows with the equivalence $\|u-u_h\|_h^2\approx \trb{u-u_h}_\pw^2 + c_\dg(u_h, u_h)$ from \cite[Thm.\ 4.1]{CarstensenGallistlNataraj2015}.
\end{proof}

\subsection{Discontinuous Galerkin 2} \label{sec:DGFEM2}
 The identity
$a(v,w)=(\Delta v, \Delta w)_{L^2(\Omega)}$ for $v,w\in V$ 
motivates
the alternative discontinuous Galerkin method from 
\cite{MozoSuli07,Georgoulis2011}
with
 discrete bilinear form
\begin{align}a_h=(\Delta_\pw\bullet, \Delta_\pw\bullet)_{L^2(\Omega)} + b_h + c_{\dg}\label{eqn:ah_dg2}.\end{align}
The semi-scalar product 
$c_\dg$ is \eqref{eqn:bhindGb} and $b_h$ reads, for any $v_2,w_2\in P_{\hspace{-.13em}2}(\T)\supset V_\nc+V_h$,
\begin{subequations}
\begin{align}
   b_h(v_2, w_2)& \coloneqq -\Theta {\cal J}(v_2, w_2) -   {\cal J}(w_2, v_2), \label{eq:bhindG2} \\
     {\cal J}(v_2, w_2) & \coloneqq  \sum_{E\in \E} \int_E 
		\jump{\frac{\partial v_2}{\partial\nu_E}}\mean{\Delta_\pw w_2} \ds
\end{align}
\end{subequations}
for $-1\leq\Theta\leq 1$. 
Appropriate 
parameters $\sigma_1,\sigma_2$ in $c_\dg$ guarantee $V_h$-ellipticity
\eqref{eq:a_ellipticity} of $a_h$ \cite{MozoSuli07}.
The bilinear form \eqref{eqn:ah_dg2} allows for \eqref{eqn:H}.
\begin{lemma}[quasi-best approximation]\label{lem:qb_dg2}
	The discontinuous Galerkin method with $a_h$ from \eqref{eqn:ah_dg2} satisfies \eqref{eqn:H} and the quasi-best
	approximation property \eqref{eqn:quasi_best}.
\end{lemma}
\newcommand{\tah}{{a_h}}
\newcommand{\tapw}{\widetilde{a_\pw}}
\newcommand{\tbh}{{b_h}}
\begin{proof}
	Given $v_h,w_h\in V_h$, 
	abbreviate $v\coloneqq J_hv_h, w\coloneqq J_hw_h\in V$ and $v_\M\coloneqq I_\M v_h, w_\M\coloneqq I_\M w_h\in \M(\T)$.
	Algebraic manipulations as in \cite[Eqn.~(6.15)]{carstensen_lowest-order_2022} reveal
	\begin{align}\label{eqn:tah_split}
		a_h(v_h, w_h) - a(v, w) &= (\Delta_\pw(v_h-v_\M), \Delta_\pw w_h)_{L^2(\Omega)} + \tbh(v_h-v_\M, w_h)\\
								&+ (\Delta_\pw v_\M, \Delta_\pw(w_h-w_\M))_{L^2(\Omega)}+\tbh(v_\M,
		w_h-w_\M)\notag\\
								&+c_\dg(v_h, w_h) + (\Delta_\pw v_\M, \Delta_\pw w_\M)_{L^2(\Omega)}%
								- a(v,w)\notag.
	\end{align}
	Cauchy inequalities, $\|\Delta_\pw \bullet\|\leq\sqrt{2}\trb{\bullet}$, the boundedness of $\tbh$, %
	and \eqref{eqn:Lambda_M} provide
	\begin{align}\label{eqn:ab_bound1}
		(\Delta_\pw(v_h-v_\M), \Delta_\pw w_h)_{L^2(\Omega)} + \tbh(v_h-v_\M, w_h)&\leq
		(\sqrt2+\|\tbh\|)\Lambda_\M\|v-v_h\|_h\|w_h\|_h.
	\end{align}
Recall the definition of the jump $\jump{\bullet}$ and average $\mean{\bullet}$ along an edge $E\in\E$ from Subsection \ref{sub:Evolution of a posteriori error analysis} and
\ref{sub:Classical and averaged Morley interpolation}
	and the product rule for jump terms $\jump{ab}=\mean{a}\jump{b}+\jump{a}\mean{b}$ for any
	$a,b\in H^1(\T)$. This and an integration by parts verify
	\begin{align}\label{eqn:ab_bound2}
		(\Delta_\pw v_\M&, \Delta_\pw(w_h-w_\M))_{L^2(\Omega)}+\tbh(v_\M,
		w_h-w_\M)\\
			  &=\sum^{}_{E\in \E} \int_E\left(\jump{\Delta_\pw v_\M}\mean{ \frac{\partial
							(w_h-w_\M)}{\partial\nu_E}}-\Theta\jump{ \frac{\partial
					v_\M}{\partial\nu_E}}\mean{\Delta_\pw (w_h-w_\M)}\right)\ds=0\notag
	\end{align}
	with $\int_E\mean{\partial(w_2-I_\M w_2)/\partial\nu_E}\ds=\int_E\jump{\partial v_\M/\partial\nu_E}\ds=0$ for any
	edge $E\in\E$ from the definition of $I_\M$ in the last step.
	Since the Morley interpolation $I_\M$ exactly interpolates the integral mean over an edge $E\in
	\E$ of the normal derivative of $w\equiv
	Jw_\M\in V$ (from $I_\M J=1$), an integration by
	parts for any $p_2\in P_{\hspace{-.13em}2}(\T)$ shows the orthogonality
	\begin{align*}
		(\Delta_\pw p_2,\Delta_\pw(w - w_\M))_{L^2(\Omega)}%
		=\sum^{}_{E\in\E} \mean{\Delta_\pw p_2}\int_E \partial((1-I_\M)Jw_\M)/\partial\nu_E\ds = 0.
	\end{align*}
	Since $a(v, w) = (\Delta v, \Delta w)_{L^2(\Omega)}$, this, a Cauchy inequality, and $\|\Delta_\pw\bullet\|\leq \sqrt 2\trb{\bullet}_\pw$ imply
	\begin{align*}
		(\Delta_\pw v_\M, \Delta_\pw w_\M)_{L^2(\Omega)}%
		- a(v, w)
		&= (\Delta_\pw(1-J)v_\M, \Delta_\pw Jw_\M)_{L^2(\Omega)}\\&\leq
		2(1+\Lambda_\M)\|J\|_h\|I_\M\|_h\|v-v_h\|_h\|w_h\|_h.
	\end{align*}
	This, the combination of \eqref{eqn:tah_split} with \eqref{eqn:ab_bound1}--\eqref{eqn:ab_bound2}, and $c_\dg(v_h,
	w_h)\leq\Lambda_{\rm c}\|v-v_h\|_h\|w_h\|_h$ for $\Lambda_{\rm c}\lesssim1$ from \cite[Sec.~7]{carstensen_lowest-order_2022} conclude the proof of
	\eqref{eqn:H}.
	The quasi-best approximation property \eqref{eqn:quasi_best} is a consequence of \eqref{eqn:H} and Theorem \ref{thm:abstract_main}.
\end{proof}
Since the dG formulations from Subsections \ref{sec:DGFEM1}--\ref{sec:DGFEM2} allow for \eqref{eqn:H} and utilize the
same space $V_h=P_{\hspace{-.13em}2}(\T)$, the a posteriori results from Subsection \ref{sec:DGFEM1} follow verbatim for the alternative
dG formulation in
this subsection.
\begin{theorem}[a posteriori estimate]\label{thm:dG_a_posteriori}
	The discrete dG solution $u_h\in V_h$ to \eqref{eqn:DWP_Vh} with $a_h$ from \eqref{eqn:ah_dg2} and the exact
	solution $u\in V$ to \eqref{eqn:WP} with $f\in
	L^2(\Omega)$ satisfy
	\begin{align*}
		\|u-u_h\|_h^2 + \mathrm{osc}_2^2(f)&\approx \|h_\T^2 f\|_{}^2 + \sum_{E\in \E}h_E\|[D_\pw^2 u_h]_E\tau_E\|_{L^2(E)}^2 +
		j_h(u_h, u_h)\\
		&\approx\|h_\T^2 f\|_{}^2 + \sum_{E\in\E}\left( h_E^{-3}\left\|\left[u_h\right]_E\right\|_{L^2(E)}^2 + h_E^{-1}
		\left\|\left[\frac{\partial u_h}{\partial \nu_E}\right]_E\right\|_{L^2(E)}^2\right).
	\end{align*}
	The
	equivalence constants exclusively depend on the shape-regularity of $\T$.\qed
\end{theorem}
The following corollary provides an a posteriori error estimator that comes without the jump term $\|[\Delta_\pw u_h]_E\|_{L^2(E)}$ over an interior edge
$E\in\E(\Omega)$ and
so refines the a posteriori result in \cite{Georgoulis2011}.  
\begin{cor}[\cite{Georgoulis2011}]\label{cor:dG_corollary}
	The discrete dG solution $u_h\in V_h$ to \eqref{eqn:DWP_Vh} with $a_h$ from \eqref{eqn:ah_dg2} and the exact
	solution $u\in V$ to \eqref{eqn:WP} with $f\in
	L^2(\Omega)$ satisfy
	\begin{align*}
		\trb{u-u_h}_\pw^2 + c_\dg(u_h, u_h)+  \mathrm{osc}_2^2(f)&\approx \|h_\T^2 f\|_{}^2 + c_\dg(u_h, u_h).\qed
	\end{align*}
\end{cor}
\subsection{$C^0$ interior penalty ($C^0$IP)} \label{sub:C0IP}
The bilinear form $a_h=a_\pw + b_h + c_{\ip}$ for $C^0$IP  \cite{BS05,CCGMNN18} utilizes $b_h$ from \eqref{eq:bhindG}
and
depends on the parameter $\sigma_\ip>0$ in
\begin{subequations}
\begin{align}
    c_\ip(v_2,w_2)& :=  
    \sum_{E \in \E} \frac{\sigma_\ip}{h_E} \int_E \jump{\frac{\partial v_2}{\partial \nu_E}} \jump{\frac{\partial w_2}{\partial \nu_E}} \ds \label{eq:chinC0IP} 
\end{align}
\end{subequations}
for $v_2, w_2\in V_h\coloneqq P_{\hspace{-.13em}2}(\T)$.
The scheme is a modification of the dG method in Section \ref{sec:DGFEM1} with trial and test functions restricted to
the continuous piecewise polynomials $V_h\coloneqq S^2_0(\T)$.
For $\sigma_\ip\approx 1$ sufficiently large but bounded, the bilinear form is coercive.
The abstract framework applies the transfer
operator $I_h=I_\C\in L(V_\nc;V_h)$ from Subsection \ref{sub:Transfer_operator_Morley}.

\begin{theorem}[a posteriori estimate]\label{thm:C0IP_a_posteriori}
	The discrete solution $u_h\in V_h$ of the $C^0$IP method to \eqref{eqn:DWP_Vh} and the exact
	solution $u\in V$ to \eqref{eqn:WP} with $f\in L^2(\Omega)$ satisfy
	\begin{align*}
		\|u-&u_h\|_h^2 +  \mathrm{osc}_2^2(f) \\
			   &\approx \|h_\T^2 f\|_{}^2 + \sum_{E\in \E}h_E\|[D_\pw^2 u_h]_E\tau_E\|_{L^2(E)}^2 +
								   \sum^{}_{E\in\E(\Omega)} h_E\|[\partial_{\nu\nu}^2u_h]_E\|_{L^2(E)}^2+
		j_h(u_h, u_h)\\
			   &\approx \|h_\T^2 f\|_{}^2 + \sum_{E\in\E}h_E^{-1}
		\left\|\left[\frac{\partial u_h}{\partial \nu_E}\right]_E\right\|_{L^2(E)}^2+
								   \sum^{}_{E\in\E(\Omega)} h_E\|[\partial_{\nu\nu}^2u_h]_E\|_{L^2(E)}^2.
	\end{align*}
	The equivalence constants exclusively depend on the shape-regularity of $\T$.
\end{theorem}
\begin{proof}
	Set $w\coloneqq (1-\JIM I_\C I_\M) e\in V$ and $\widehat w\coloneqq (Q-\JIM)I_\C I_\M e\in V+ S^2_0(\T)$ for $e\coloneqq u-\JIM
	u_h\in V$.
	The definition of the residual, Corollary \ref{cor:F_bound}, and Lemma \ref{lem:C0IP_a} with $I_h=I_\C$ show
	\begin{align*}
		Res(w)\coloneqq F(w) - a_\pw(u_h, w)&\lesssim \left(\|h_\T^2 f\| + \sqrt{\sum^{}_{E\in\E(\Omega)}
		h_E\|[\partial_{\nu\nu}^2u_h]_E\|_{L^2(E)}^2}\right)\trb{e}.%
	\end{align*}
Theorem \ref{thm:a_posteriori_bound} and the definition of the residual result in
$ \|u-u_h\|_h^2 \lesssim \|u_h-\JIM u_h\|^2 + Res(w)-\hatF(\widehat w)$.
Since the stability of the $L^2$ projection shows $\mathrm{osc}_2(f)\leq \|h_\T^2f\|$, this, the bound
$\hatF(\widehat w)\lesssim \|h_\T^2f\|\trb{e}$ from \eqref{eqn:fQ_bound}, and a weighted Young inequality reveal
	\begin{align*}
		\|u-u_h\|_h^2 +\mathrm{osc}_2^2(f)&\lesssim \|h_\T^2 f\|^2 + \|u_h-\JIM u_h\|^2_{\W{h}} +\sum^{}_{E\in\E(\Omega)}
		h_E\|[\partial_{\nu\nu}^2u_h]_E\|_{L^2(E)}^2.
	\end{align*}
	Theorem \ref{thm:mu_bound} bounds $\|u_h-\JIM u_h\|^2$ either in terms of $\sum^{}_{E\in\E} h_E\|[D^2_\pw
u_h]_E\W{\tau_E}\|_{L^2(E)}^2$ plus $j_h(u_h,
	u_h)$ or in terms of $\sum^{}_{E\in\E} h_E^{-1}\|[\partial
	u_h/\partial\nu_E]_E\|_{L^2(E)}^2$ (because $[u_h]_E \equiv 0$ for $u_h\in S^2_0(\T)$).
	This concludes the proof of the reliability.
	Lemma \ref{lem:C0IP_eficiency} provides the efficiency of the normal-normal jumps.
The efficiency of the remaining terms follows verbatim as in the proof of Theorem \ref{thm:dG_a_posteriori}.
\end{proof}
The following  corollary recovers the a posteriori result from \cite[Sections 3 and 4]{BGS10}.
\begin{cor}[\cite{BGS10}]
	The discrete $C^0$IP solution $u_h\in V_h$ to \eqref{eqn:DWP_Vh} and the exact solution
	$u\in V$ to \eqref{eqn:WP} with $f\in L^2(\Omega)$ satisfy
	\begin{align*}
		\trb{u-u_h}_\pw^2 + c_\ip(u_h, u_h)+  \mathrm{osc}_2^2(f)&\approx \|h_\T^2 f\|_{}^2 + c_\ip(u_h, u_h)+
								   \sum^{}_{E\in\E(\Omega)} h_E\|[\partial_{\nu\nu}^2u_h]_E\|_{L^2(E)}^2.
	\end{align*}
\end{cor}
\begin{proof}
Since $[v_h]_E=0$ for any $v_h\in S^2_0(\T)$, $c_\ip=c_\dg$ coincide in $S^2_0(\T)\times S^2_0(\T)$ and the proof follows verbatim that of Corollary \ref{cor:dG_corollary}; further details are omitted.
\end{proof}
\subsection{WOPSIP}%
\label{sub:WOPSIP method}
The weakly over-penalized symmetric interior penalty
(WOPSIP) scheme \cite{BrenGudiSung10} is a penalty method  
with the stabilisation term 
\begin{align}\label{eqn:cP}
	c_\w(v,w)&\coloneqq
	\sum_{E \in \E}h_E^{-2}\left( \sum_{z \in {\mathcal V} (E)} \frac{[v]_E(z)}{h_E}\frac{[w]_E(z)}{h_E}
			 +
\fint_E  \jump{\frac{\partial v}{\partial\nu_E}}\!\!\!\mathrm ds\, \fint_E  \jump{\frac{\partial w}{\partial
\nu_E}}\!\!\!\mathrm ds\right)
\end{align}
for piecewise smooth functions $v,w\in H^2(\T)$.
The difference of \W{$c_\w$} in \eqref{eqn:cP} to $j_h$ from \eqref{eqn:apw_jh} is the over-pernalisation by an additional negative power of the mesh size $h_E$.
This and $h_{\rm max}\coloneqq \mathrm{max}_{T\in\T}h_T$ establish
\begin{align}\label{eqn:cP_bound}
	j_h(v, v) \leq h_{\rm max}^2 c_\w(v,v) && \text{ for all }v\in \widehat V\coloneqq H^2(\T).
\end{align}
Hence \W{$\|\bullet\|_\w\coloneqq \left(\trb{\bullet}_\pw^2 + c_\w(\bullet, \bullet)\right)^{1/2}$} is
a norm in $\widehat V$ stronger than $\|\bullet\|_h$.
The WOPSIP method computes the 
discrete
solution $u_h\in V_h\coloneqq P_{\hspace{-.13em}2}(\T)$ to \eqref{eqn:DWP_Vh} with the bilinear form
$a_h\coloneqq a_\pw + c_\w$ and fits into the abstract setting with $V_\nc\coloneqq \M(\T)$.

The main difference to the methods under consideration above is the missing quasi-best approximation
property due to the penalisation.
Instead of this, the following a priori estimate for the energy norm 
\begin{align*}\trb{u-u_h}_\pw^2 + c_\w(u_h, u_h) \lesssim \trb{u-I_\M u}_\pw^2 + \trb{h_\T I_\M u}_\pw^2\end{align*}
holds with the extra term $\trb{h_\T I_\M u}_{\pw}^2$ \cite[Theorem 9.1]{carstensen_lowest-order_2022}.
This suggests that \eqref{eqn:H} does not hold, but the methodology of the a posteriori analysis of Subsection
\ref{sub:Paradigm} is still applicable. %
Indeed, the key assumption only enters in the error bound from Theorem \ref{thm:a_posteriori_bound} and a careful
analysis with $I_h=\id$ leads to \eqref{eqn:a_posteriori_bound}.
This allows the application of the developed tool chain and leads to a new a posteriori estimate \emph{without} the WOPSIP stabilisation
term \eqref{eqn:cP} but still with the weaker stabilization $j_h$.

\begin{theorem}[a posteriori estimate]\label{thm:WOPSIP_a_posteriori}
	The WOPSIP solution $u_h\in P_{\hspace{-.13em}2}(\T)$ to \eqref{eqn:DWP_Vh} and the exact solution $u\in V$
	to \eqref{eqn:WP} with $f\in L^2(\Omega)$ satisfy
	\begin{align*}
		\|u-u_h\|_h^2 + \mathrm{osc}_2^2(f)&\approx \|h_\T^2 f\|_{}^2 + \sum_{E\in \E}h_E\|[D_\pw^2 u_h]_E\tau_E\|_{L^2(E)}^2 +
		j_h(u_h, u_h)\\
		&\approx\|h_\T^2 f\|_{}^2 + \sum_{E\in\E}\left( h_E^{-3}\left\|\left[u_h\right]_E\right\|_{L^2(E)}^2 + h_E^{-1}
		\left\|\left[\frac{\partial u_h}{\partial \nu_E}\right]_E\right\|_{L^2(E)}^2\right).
	\end{align*}
	The equivalence constants exclusively depend on the shape-regularity of $\T$.
\end{theorem}
\begin{proof}
	Let $e\coloneqq u-J_h u_h\in V$  and recall $J_h\coloneqq JI_\M$ as well as $I_\M J=\id$ on $V_\nc$.
	The key assumption \eqref{eqn:H} enters the proof of Theorem \ref{thm:a_posteriori_bound} with
	\eqref{eqn:H_corollary}.
	This proof exploits that the transfer operator $I_h\coloneqq \id:V_\nc \to V_h$ is
	the identity and deduces \eqref{eqn:H_corollary} directly (and so circumvents \eqref{eqn:H}).
	Indeed, since $c_\w(\bullet, v_\nc) = 0$ for any $v\in V_\nc$ and $a_\pw(u_h, e_h - J_he_h) = 0$ by the orthogonality
	\eqref{eqn:best_approx} for $e_h\coloneqq I_h I_\M e=I_\M J_h e_h\in V_\nc$, 
	\begin{align*}
		a_h(u_h, e_h)&=a_\pw(u_h, e_h)=a_\pw(u_h, J_he_h)\quad\text{and}\\
		a_h(u_h, e_h) -a(J_hu_h, J_he_h)
					 &= a_\pw(u_h - J_hu_h, J_he_h) \leq
		\|J_h\|\trb{u_h - J_hu_h}_\pw\trb{e_h}_\pw
	\end{align*}
	follow 
	with a Cauchy inequality in the last step.
	Hence \eqref{eqn:H_corollary} even holds with the weaker norm $\trb{\bullet}_\pw\leq \|\bullet\|_h$.
	The remaining parts of the proof for Theorem \ref{thm:a_posteriori_bound} apply analogously and verify
\begin{align*}
	\|u-u_h\|_h^2\lesssim \|u_h-\JIM u_h\|^2_{\W{h}} + Res(w) - \int_\Omega f \widehat w\dx  = \|u_h-\JIM u_h\|^2_{\W{h}} +
	\int_\Omega f(w-\widehat w)\dx
\end{align*}
with $a_\pw(u_h, w) = 0$ from Lemma \ref{lem:C0IP_a} for $w\coloneqq e - J_h I_h
I_\M e$ and $\widehat w\coloneqq Qe_h-J_he_h\in V+V_h$.
The remaining arguments follow the proofs of Theorem \ref{thm:dG_a_posteriori} and Theorem \ref{thm:Morley_a_posteriori}
verbatim.
\end{proof}
The inclusion of the stabilisation term $c_\w$ on both sides of the error estimate in Theorem
\ref{thm:WOPSIP_a_posteriori} recovers the a posteriori estimate from \cite[Section 6]{BrenGudiSung10}.
\begin{cor}[{\cite{BrenGudiSung10}}]\label{cor:WOPSIP_corollary}
	The discrete WOPSIP solution $u_h\in V_h$ to \eqref{eqn:DWP_Vh} and the exact solution $u\in
	V$ to \eqref{eqn:WP} with $f\in L^2(\Omega)$ satisfy
	\begin{align*}
		\trb{u-u_h}_\pw^2 &+ c_\w(u_h, u_h)+  \mathrm{osc}_2^2(f)\\&\approx \|h_\T^2 f\|_{}^2 + \sum_{E\in\E}\left( h_E^{-3}\left\|\left[u_h\right]_E\right\|_{L^2(E)}^2 + h_E^{-1}
		\left\|\left[\frac{\partial u_h}{\partial \nu_E}\right]_E\right\|_{L^2(E)}^2\right) + c_\w(u_h,
		u_h).
	\end{align*}
	The equivalence constants exclusively depend on the shape-regularity of $\T$.\qed
\end{cor}
\begin{proof}
	This follows from Corollary \ref{cor:dG_corollary} with $c_\w(\bullet, v)=0$ for all $v\in V$ and $j_h(u_h,
u_h)\lesssim c_\w(u_h, u_h)$ from \eqref{eqn:cP_bound}; further details are omitted.
\end{proof}

\section{More general sources}%
\label{sec:General sources}
This section considers a class of rather general right-hand sides $F\in V^*$ and introduces an estimator for the
residual that is reliable and efficient up to a data approximation error.
\subsection{A general class of source terms}%
\label{sub:A class of right-hand sides}
Every functional in $F\in V^*\equiv H^{-2}(\Omega)$ has (non-unique) representations by volume loads 
$f_\alpha\in L^2(\Omega)$ for all 6 multi-indices
	$\alpha=(\alpha_1,\alpha_2)\in \N_0^2$ of order $|\alpha|\coloneqq \alpha_1+\alpha_2\leq 2$, written
	$(f_\alpha)_{|\alpha|\leq 2}\in L^2(\Omega)^6$, with
	\begin{align}
		\label{eqn:F_characterisation1}
		F(\varphi)\equiv \langle F, \varphi\rangle=\sum^{}_{|\alpha|\leq 2}
		(f_\alpha,\partial^\alpha\varphi)_{L^2(\Omega)}\qquad\text{ for all }\varphi\in H^2_0(\Omega).
	\end{align}
\begin{theorem}[characterization]\label{thm:F_char}
	Given any $F\in H^{-2}(\Omega)$ there exist $(f_\alpha)_{|\alpha|\leq 2}\in L^2(\Omega)^6$ such that
	\eqref{eqn:F_characterisation1} holds.
The norm of $F$ in $H^{-2}(\Omega)$ (the dual of $H^2_0(\Omega)$ endowed with the \W{full} Sobolev norm of
	$H^2(\Omega)$) is the minimum
	\begin{align*}
		\|F\|_{H^{-2}(\Omega)}=\min\left\{\sqrt{\sum^{}_{|\alpha|\leq 2} \|f_\alpha\|_{L^2(\Omega)}^2}\
		:\ (f_\alpha)_{|\alpha|\leq 2}\in L^2(\Omega)^6\text{ satisfies }\eqref{eqn:F_characterisation1}\right\}.
	\end{align*}
\end{theorem}
\begin{proof}
This is a natural generalization of the corresponding result for functionals in $H^{-1}(\Omega)$, e.g., 
\cite[Sec.~5.9, Thm.~1]{evans_partial_2010}; hence further details are omitted.
\end{proof}
\begin{rem}[characterization for semi-norm $\trb{\bullet}$]
	The norm representation of Theorem \ref{thm:F_char} is given in the (full) norm $\|\bullet\|_{H^{2}(\Omega)}$ of $H^{2}(\Omega)$.
	A corresponding assertion
	\begin{align}\label{eqn:F_char_semi}
		\trb{F}_{*}\coloneqq\sup_{v\in V}F(v)/\trb{v}=\min_{\sigma\in L^2(\Omega;\mathbb S)}\{\|\sigma\|_{L^2(\Omega)}\ :\ F=(\sigma,
		D^2\bullet)_{L^2(\Omega)}\}
	\end{align}
	follows from the Riesz representation theorem
	for the $H^2$ seminorm $\trb{\bullet}\equiv|\bullet|_{H^2(\Omega)}$ as well.
	The minimizer $\sigma=D^2u\in L^2(\Omega;\mathbb S)$
	in \eqref{eqn:F_char_semi} is the Hessian of the weak solution $u\in V$ to \eqref{eqn:WP}.
\end{rem}
A more general source term may include point forces
$\delta_z\in V^*$ at finitely many points $z\in A\subset\overline\Omega$ and line loads
 $(g_0, \bullet)_{L^2(\Gamma_0)}$, $(g_1, \partial_\nu\bullet)_{L^2(\Gamma_1)}$ along the hypersurfaces
$\Gamma_0,\Gamma_1\subset\Omega$ in addition to \eqref{eqn:F_characterisation1}.
The Dirac delta distribution $\delta_z(f)=f(z)$ evaluates $f\in V\subset C(\overline\Omega)$ at the atom $z$ and we
suppose for simplicity that the mesh is adapted in that $A\subset\mathcal{V}(\Omega)$ consists of interior vertices.
Recall the set $\T(z)\coloneqq\{T\in\T\ :\ z\in T\}$ of neighbouring triangles from Subsection \ref{sub:Classical and
	averaged Morley interpolation} and suppose that the mesh resolves $\Gamma_j=\bigcup
	\E(\Gamma_j)$ with $\E(\Gamma_j)\coloneqq\{E\in\E\ :\
\mathrm{int}(E)\subset\Gamma_j\}$ for $j=0,1$.

This section considers sources $F\coloneqq \hatF|_V\in V^*$ in terms of an extended source
 $\hatF\in \widehat V^*\equiv H^2(\T)^*$,
defined, for $\hatv\in H^2(\T)$, by 
\begin{align}\label{eqn:hatF_repr}
		\hatF(\hatv)\coloneqq
		&\;
		\sum^{}_{|\alpha|\leq2} (f_\alpha,\partial^\alpha_\pw \hatv)_{L^2(\Omega)} 
		+ \sum_{j=0,1}(g_j,
		\langle \partial_\nu^j\hatv\rangle_{\Gamma_j})_{L^2(\Gamma_j)}+ \sum_{z\in
		A}\sum_{T\in\T(z)}\beta_{T,z}\hatv|_T(z).
\end{align}
The given data in \eqref{eqn:hatF_repr} are Lebesgue functions $(f_\alpha)_{|\alpha|\leq 2}\in L^2(\Omega)^6$,  line loads $g_j\in
L^2(\Gamma_j)$ along the hypersurface $\Gamma_j\subset\bigcup\E(\Omega)$ for $j=0,1$, and point forces of intensity
$\beta_z=\sum^{}_{T\in\T(z)} \beta_{T,z}\in\mathbb R$ at $z\in A\subset\mathcal{V}(\Omega)$.
\begin{rem}[influence of $\hatF$]
	Since $\hatF\circ J_h=F\circ J_h$ holds (for the five schemes from Section \ref{sec:Unified a posteriori error
	control}) with the smoother $Q=J_h$, the discrete solution
	$u_h\in V_h$
	to \eqref{eq:discrete2} depends on $F$ but is independent of its representation \eqref{eqn:hatF_repr}.
	The classical scheme without smoother
	$Q=\id$, however, depends on the chosen data for the representation $\hatF$. %
\end{rem}

Throughout this section, we suppose that we have 
 piecewise smooth approximations $\tg_j\in L^2(\Gamma_j)$ of $g_j$ for $j=0,1$ and $\tf_\alpha\in H^{|\alpha|}(\T)$
of $f_\alpha$ in \eqref{eqn:hatF_repr} for
	$|\alpha|\leq 2$ to define an approximation $\hatFapx$ of
	$\hatF$ with piecewise smooth data. 
	The reason for this approximation is that we shall integrate by parts with piecewise smooth functions to
	reveal an efficient a posteriori upper error bound in the subsequent subsection.

\begin{defn}[data approximation error]\label{def:hatF_osc}
	The approximated source term reads
	\begin{align}\label{eqn:hatF_osc}
		\hatFh(\hatv)\coloneqq
		&\;\sum_{|\alpha|\leq 2}(\tf_\alpha,\partial^\alpha_\pw \hatv)_{L^2(\Omega)} 
		+ \sum^{}_{j=0,1} \sum_{E\in\E(\Gamma_j)}(\tg_j,
		\langle \partial_\nu^j\hatv\rangle_E)_{L^2(E)}%
\end{align}
for all $\hat v\in H^2(\T)$. The data approximation error
$\APX(F,\T)\coloneqq\big(\sum^{}_{T\in\T} \APX^2(F, T)\big)^{1/2}$ has, on the triangle $T\in\T$, the contribution
\begin{align*}
		\APX^2(F, T)\coloneqq \sum^{}_{|\alpha|\leq 2} \|h_T^{2-|\alpha|}(f_\alpha-\tf_\alpha)\|^2_{L^2(T)} +
		\sum^{}_{j=0,1} \sum^{}_{E\in\E(\Gamma_j)\cap \E(T)} \|h_E^{3/2-j}(g_j-\tg_j)\|^2_{L^2(E)}.
	\end{align*}
\end{defn}
The data approximation error generalizes data oscillations.
Let $\Pi_{E, k}:L^2(E)\to P_{\hspace{-.13em}k}(E)$ denote the $L^2$ projection onto $P_{\hspace{-.13em}k}(E)$ on the edge $E\in\E$.
\begin{example}[data oscillations]\label{rem:osc}
	The natural candidates for $(F_\alpha)_{|\alpha|\leq 2}$ and $G_0, G_1$ in \eqref{eqn:hatF_osc} are $L^2$ projections onto polynomials of degree at most $k\in\N_0$. Then the data approximation error
	$\APX^2(F,\T)$ becomes an oscillation term
	\begin{align*}
		\OSC^2(F, \T)\coloneqq\sum^{}_{|\alpha|\leq 2} \|h_\T^{2-|\alpha|}(1-\Pi_k)f_\alpha\|^2_{L^2(\Omega)} +
		\sum^{}_{j=0,1} \sum^{}_{E\in\E(\Gamma_j)} \|h_E^{3/2-j}(1-\Pi_{E,k})g_j\|^2_{L^2(E)}.
	\end{align*}

\end{example}
\begin{lemma}[data approximation error]\label{lem:F_osc}
	With the linear operators $J, I_\M, I_h$ from Table \ref{long1},%
	\begin{align*}
		\max\left\{\trb{(F-\hatFapx)(1-\JIM I_hI_\M)}_{*}, \trb{(\hatF-\hatFh)(1 - \JIM)I_hI_\M}_{*}\right\}\leq\const{cst:5} \APX(F, \T).
	\end{align*}
\end{lemma}
\medskip
\begin{proof}
	Recall that $w\coloneqq(1-\JIM I_hI_\M)v$ vanishes at the vertices for all $v\in V$ for all five schemes under
	consideration.
	This shows
	\begin{align*}%
		(F-\hatFapx)(w)&=\sum^{}_{|\alpha|\leq 2} (f_\alpha-F_\alpha, \partial^\alpha w)_{L^2(\Omega)} +
		\sum^{}_{j=0,1} (g_j-G_j,\partial_\nu^jw)_{L^2(\Gamma_j)}\leq \const{cst:5}\APX(F, \T)\trb{v}
	\end{align*}
	with a Cauchy inequality and the constant $\const{cst:5}$ from Lemma \ref{lem:5.1} in the last step. Analog arguments provide the asserted bound of
	$\trb{(\hatF-\hatFh)\circ(1 - \JIM)I_hI_\M}_{*}$.
\end{proof}

\subsection{Estimator for the residual}%
\label{sub:Estimator for the residual}
	The paradigm shift in this paper is that Theorem \ref{thm:a_posteriori_bound} provides an upper error bound
	with a specific structure
	of the test function as an element in $(1-J_hI_hI_\M)V$ for the residual part.
This subsection designs an estimator $\mu(\T)$ for the dual norm $\trb{Res\circ(1-\JIM I_hI_\M )}_{*}$ of the residual
that is reliable and efficient up to the data approximation error $\APX(F,\T)$ %
	\begin{align}\label{eqn:mu_res_bound}
	\trb{Res\circ(1-\JIM I_hI_\M )}_{*}\lesssim \mu(\T) + \APX(F, \T)\lesssim \trb{u-u_h}_{\W{\pw}}+\APX(F, \T).
	\end{align}
	The residual $Res\coloneqq F- a_\pw(u_h,\bullet)\in V^*$ includes the discrete solution $u_h\in V_h$ to
	\eqref{eq:discrete2} with or without smoother $Q\in\{\id, J_h\}$.
The analysis in this section for an upper bound of the dual norm $\trb{Res\circ(1-\JIM I_hI_\M )}_{*}$
	allows for a general discrete object $u_h\in V_h$; said differently, $u_h\in V_h$ is arbitrary in
	\eqref{eqn:mu_res_bound}.
	
	To define the estimator contributions in $\mu(\T)$,
	abbreviate $\tf_0\coloneqq \tf_{(0,0)}\in L^2(\Omega)$,
	\begin{align}
		\tf_1\coloneqq \begin{pmatrix}
			\tf_{(1,0)}\\\tf_{(0,1)}
		\end{pmatrix}\in H^1(\T;\mathbb{R}^2),\text{ and}\quad\tf_2\coloneqq
	\begin{pmatrix}
		\tf_{(2,0)}&\frac{1}{2}\tf_{(1,1)}\\\frac{1}{2}\tf_{(1,1)}&\tf_{(0,2)}
	\end{pmatrix}\in H^2(\T;\mathbb S).\label{eqn:F_abbr}
		\end{align}
		The extra factor $1/2$ in the definition of $F_2$ allows the simplification $(F_j, D^j_\pw\bullet)_{L^2(\Omega)}=\sum^{}_{|\alpha|=j} (F_\alpha,
		\partial^\alpha_\pw\bullet)_{L^2(\Omega)}$ for $j=0,1,2$.
Let the divergence \begin{align*}\div_\pw \tf_2\coloneqq \begin{pmatrix}
	\div_\pw (\tf_2)_1\\ \div_\pw (\tf_2)_2
\end{pmatrix}\in H^1(\T;\mathbb R^2)\end{align*}
of the matrix-valued function
$\tf_2\equiv\left((\tf_2)_1;(\tf_2)_2\right)\in H^2(\T;\mathbb S)$ apply row-wise.
	Recall $J, I_\M, I_h$ from %
	Section \ref{sec:Examples of second-order finite element schemes} and the special treatment of $I_h={\id}$ in
	Subsection \ref{sub:Error estimates}. %
	Define
\begin{align*}
	\mu_1^2(\T)&\coloneqq\|h_\T^2(\tf_0 - \div_\pw \tf_1+\div_\pw^2\tf_2)\|^2, \\
	\mu_2^2(\T)&\coloneqq \sum^{}_{E\in\E(\Omega)}
	h_E^{3}\|\tg_0+[\tf_1-\div_\pw \tf_2-\partial (\tf_2\tau_E)/\partial s]_E\cdot\nu_E\|_{L^2(E)}^2,\\
\mu_3^2(\T)&\coloneqq \sum^{}_{E\in\E(\Omega)}\begin{cases}{}
	h_E\|(1-\Pi_{E,0})(\tg_1+[\tf_2\nu_E]_E\cdot \nu_E)\|_{L^2(E)}^2&\text{if }I_h={\id},\\
h_E\|\tg_1+[(\tf_2-D^2_\pw u_h)\nu_E]_E\cdot \nu_E\|_{L^2(E)}^2&\text{if }I_h=I_\C,
\end{cases}\\
		\mu^2(\T)
		&\coloneqq \mu_1^2(\T) + \mu_2^2(\T) + \mu_3^2(\T).
\end{align*}
		Here 
		$G_j\in L^2(\Gamma_j)\subset L^2(\bigcup\E)$ is extended by zero to the entire skeleton for $j=0,1$.%
\begin{proposition}[reliability]\label{prop:mu_rel}
	The estimator
	$\mu(\T)\equiv \mu^2(\T)^{1/2}$ of the residual is reliable
		\begin{align*}
			\const{cst:5}^{-1}\W{\trb{\WW{Res\circ(1-\JIM I_hI_\M )}}_{*}}\leq \mu(\T) + \APX(F, \T).
		\end{align*}
\end{proposition}
\begin{proof}
	Given any $v\in V$, the function $w\coloneqq v-\JIM I_hI_\M v\in V$ vanishes at the vertices $z\in\V$.
	The split $\nabla w =
	(\partial w/\partial\nu_E) \nu_E + (\partial w/\partial s)\tau_E$ 
	along an edge $E=\mathrm{conv}\{A, B\}\in\E$ and an integration by parts with $w(A) = w(B) = 0$ verify
	\begin{align}
		([\tf_2\nu_E]_E, \nabla w)_{L^2(E)}
			&=([\tf_2\nu_E]_E,\nu_E \partial w/\partial\nu_E)_{L^2(E)}+ \left([\tf_2\nu_E]_E,
		\tau_E \partial w/\partial s\right)_{L^2(E)}\notag\\
			&=([\tf_2\nu_E]_E,\nu_E \partial w/\partial\nu_E)_{L^2(E)}- \left(\partial[\tf_2\tau_E]_E
			/\partial s, \nu_E w\right)_{L^2(E)}\label{eqn:S_ibp}
	\end{align}
	 with $\tau_E\cdot\tf_2\nu_E=\nu_E\cdot\tf_2\tau_E$ for all symmetric matrix-valued $\tf_2\in H^2(\T;\mathbb S)$ in the last
	step.
	An integration by parts and \eqref{eqn:S_ibp} lead to 
	\begin{align}
		\hatFapx(w)
		=&\,(\tf_0-\div_\pw \tf_1+\div_\pw^2 \tf_2,
		w)_{L^2(\Omega)}+\W{\sum_{E\in\E(\Omega)}}(\tg_1+[\tf_2\nu_E]_E\nu_E,\partial w/\partial\nu_E)_{L^2(E)}
		\nonumber\\
	  &+
	  \sum^{}_{E\in\E(\Omega)}(\tg_0+[\tf_1-\div_\pw \tf_2-\partial (\tf_2\tau_E)/\partial s]_E\cdot\nu_E, w)_{L^2(E)}.\label{eqn:Fapx_intbp}
	\end{align}
	Since $I_\M w=0$ for $I_h=\id$, the integral mean $\Pi_{E,0}(\partial w/\partial\nu_E)\equiv 0$ vanishes along any edge
  $E\in\E$.
  Hence, $(p_0, \partial w/\partial\nu_E)_{L^2(E)}=0$ is zero for any constant $p_0\in P_{\hspace{-.13em}0}(E)$.
  An integration by parts with \eqref{eqn:a_pw_ibp} for the piecewise constant Hessian $D^2_\pw u_h\in
  P_{\hspace{-.13em}0}(\T;\mathbb S)$, $q_0\coloneqq \Pi_{E,0}(\tg_2+[ \tf_2\nu_E]_E\cdot \nu_E)\in P_{\hspace{-.13em}0}(E)$, and the split of $\nabla w$ as in \eqref{eqn:S_ibp} result in
  \begin{align*}
	  a_\pw(u_h, w) = (D^2_\pw u_h, D^2 w)_{L^2(\Omega)}=\sum_{E\in\E(\Omega)}\begin{cases}{}
(q_0,
\partial w/\partial\nu_E)_{L^2(E)}  &\hspace{-1.8em}\text{if }I_h=\id,\\
([ D^2_\pw u_h\nu_E]_E\cdot \nu_E,
\partial w/\partial\nu_E)_{L^2(E)}  &\text{ else.}
	  \end{cases} 
  \end{align*}
  This and the Cauchy inequality reveal
  \begin{align*}
	  \hatFapx(w)- a_\pw(u_h,w)
			  &\leq \mu(\T)\sqrt{\|h_\T^{-2}w\|^2+\sum^{}_{E\in \E(\Omega)}
	  \bigg(\|h_E^{-3/2}w\|^2_{L^2(E)}+\|h_E^{-1/2}\frac{\partial w}{\partial\nu_E}\|^2_{L^2(E)}\bigg)}\\
	  &\leq \const{cst:5}\mu(\T)\trb{v}
	\end{align*}
	with the constant $\const{cst:5}$ from Lemma \ref{lem:5.1} in the last step.
	This and Lemma \ref{lem:F_osc} provide $Res(w) = (F-\hatFapx)(w) + \hatFapx(w) - a_\pw(u_h, w)\lesssim \mu(\T)+\APX(F,\T)$.
\end{proof}
\begin{proposition}[efficiency up to data approximation]\label{prop:mu_efficiency}
	Let $u\in V$ solve \eqref{eqn:WP} with the right-hand side $F\equiv \hatF|_V\in V^*$ given by \eqref{eqn:hatF_repr}.
	If $\tg\in P_{\hspace{-.13em}k}(\E)$ and $(\tf_\alpha)_{|\alpha|\leq 2}\in P_{\hspace{-.13em}k}(\T)^6$ are piecewise polynomials of degree at most
	$k\in\N_0$, 
	then the estimator $\mu(\T)$ of the residual 
	is efficient up to the data approximation error
	\begin{align*}
		\const{cst:mu_efficiency}^{-1}\mu(\T)\leq \trb{u-u_h}\W{_\pw} + \APX(F, \T).
	\end{align*}
	The constant $\const{cst:mu_efficiency}$ exclusively depends on the shape-regularity of $\T$ and on $k\in\mathbb
	N_0$.
\end{proposition}
Before the technical proof of Proposition \ref{prop:mu_efficiency} follows in Subsection \ref{sub:Proof of mu_eff}, the extension of the a posteriori analysis from Section \ref{sec:Unified a posteriori error
control} to $F\in V^*$ is in order.

\subsection{Application to lowest-order schemes}%
\label{sub:Examples}
This subsection extends the a posteriori error control from Section \ref{sec:Unified a posteriori error control} for the
right-hand side $F\in L^2(\Omega)$ to a general source $F\equiv\hatF|_V\in V^*$ from \eqref{eqn:hatF_repr}.
In fact, the efficient bounds of $\|u_h-\JIM u_h\|_h$ from Theorem \ref{thm:mu_bound} imply the following novel result
generalizing \cite[Thm.~6.2]{NNCC2020} for $Q=J_h$.
Let $(F_\alpha)_{|\alpha|\leq 2}\in P_{\hspace{-.13em}k}(\T)^6, (G_0,G_1)\in P_{\hspace{-.13em}k}(\E)^2$ be piecewise polynomials of degree at
most $k\in\mathbb N_0$ that enter Definition \ref{def:hatF_osc} for the data approximation error $\APX(F, \T)$.%
\begin{theorem}[a posteriori for $Q=J_h$]\label{thm:F_a_post}
	Let $u_h\in V_h$ solve \eqref{eq:discrete2} with $Q=J_h$ for any of the five discrete schemes from Section \ref{sec:Unified
	a posteriori error control} and let $u\in V$ solve
	\eqref{eqn:WP}. Then
	\begin{align*}
		\|u-u_h\|_{h}^2&+\APX^2(F, \T)\\
					   &\approx \mu^2(\T) + \sum_{E\in \E}h_E\|[D_\pw^2 u_h]_E\tau_E\|_{L^2(E)}^2 +
		j_h(u_h, u_h)+\APX^2(F, \T)\\
		&\approx\mu^2(\T) + \sum_{E\in\E}\left( h_E^{-3}\left\|\left[u_h\right]_E\right\|_{L^2(E)}^2 + h_E^{-1}
		\left\|\left[\frac{\partial u_h}{\partial \nu_E}\right]_E\right\|_{L^2(E)}^2\right)+\APX^2(F, \T).
	\end{align*}
	The hidden equivalence constants exclusively depend on the shape-regularity of $\T$ and on $k\in \mathbb
	N_0$.
\end{theorem}
\begin{proof}
	Theorem \ref{thm:a_posteriori_bound} provides $\|u-u_h\|_h\lesssim Res(w) + \|u_h-\JIM u_h\|_h$ for $w= v-\JIM I_hI_\M
	v$ and some $v\in V$ for Morley, dG, and $C^0$IP.
	Recall from the proof of Theorem
	\ref{thm:WOPSIP_a_posteriori} that this error bound also holds for the WOPSIP scheme even without the validity of
	\eqref{eqn:H} in full generality.
	The efficient bound $Res(w) \lesssim \mu(\T)\lesssim \trb{u - u_h}_\pw+\APX(F,\T)$ of the residual by the estimator
	$\mu(\T)$ from Proposition \ref{prop:mu_rel}--\ref{prop:mu_efficiency} and the efficient a posteriori control of
	$\|u_h-\JIM u_h\|_h$ from Theorem \ref{thm:mu_bound} conclude the proof.
\end{proof}
The original formulation \eqref{eq:discrete2} without a smoother, $Q=\id$, leads to
an additional term \begin{align*}
\hatF(e_h-J_he_h)=(\hatF-\hatFapx)(e_h-J_he_h)+\hatFapx(e_h-J_he_h)\end{align*} in the a posteriori error bound from Theorem
\ref{thm:a_posteriori_bound} and reflects the particular choice of the extended data $\hatF$ in the definition \eqref{eqn:hatF_repr}.
While the difference $\hatF-\hatFapx$ is bounded by the data approximation error $\APX(F, \T)$, the
non-conforming test function $e_h-J_he_h\not\subset V$ prevents an efficient control of the higher-order volume
sources in \eqref{eqn:hatF_osc} by residual terms through an integration by parts.

For Morley, dG, and WOPSIP, the critical terms are the intermediate sources $f_\alpha$ for $|\alpha|=1$ and the
proof below explains why those are omitted in the (reduced) model class of right-hand sides in \cite{NNCC2020}\X{\cite{carstensen_optimal_2022}}.
The following theorem generalizes \cite[Thm.~6.1]{NNCC2020} for $Q=\id$.
\begin{theorem}[a posteriori for $Q=\id$]\label{thm:Fh_a_post}
Suppose 
\begin{align}\label{eqn:F_data_ass1}
	F_\alpha&\coloneqq 0\quad\text{for all }|\alpha|=1,\quad
	F_\alpha\in P_{\hspace{-.13em}0}(\T)\quad\text{for all }|\alpha|=2
\end{align}
for Morley, dG, WOPSIP, and
\begin{align}\label{eqn:F_data_ass2}
F_\alpha&\coloneqq 0\quad\text{for all }|\alpha|=2,\quad
G_1\coloneqq0
\end{align}
for $C^0IP$.
	Let $u_h\in V_h$ solve \eqref{eq:discrete2} without smoother, $Q=\id$, for any of the five discrete schemes from Section \ref{sec:Unified
	a posteriori error control} and let $u\in V$ solve
	\eqref{eqn:WP}. Then
	\begin{align*}
		\|u-u_h\|_{h}^2&+\APX^2(F, \T)\\
					   &\approx \mu^2(\T) + \sum_{E\in \E}h_E\|[D_\pw^2 u_h]_E\tau_E\|_{L^2(E)}^2 +
					   j_h(u_h, u_h)+\APX^2(F, \T)\\
		&\approx\mu^2(\T) + \sum_{E\in\E}\left( h_E^{-3}\left\|\left[u_h\right]_E\right\|_{L^2(E)}^2 + h_E^{-1}
		\left\|\left[\frac{\partial u_h}{\partial \nu_E}\right]_E\right\|_{L^2(E)}^2\right)+\APX^2(F, \T).
	\end{align*}
	The hidden equivalence constants exclusively depend on the shape-regularity of $\T$ and on $k\in \mathbb
	N_0$.
\end{theorem}
\begin{proof}
	 For Morley, dG, and $C^0$IP, Theorem \ref{thm:a_posteriori_bound} with $e\coloneqq u-J_hu_h\in V$ and the split $\hatF(v_h)=\hatFapx(v_h) +
	(\hatF-\hatFapx)(v_h)$ for $v_h\coloneqq
	(1-J_h)I_hI_\M e\in V_h$ provide
	\begin{align}
		\|u-u_h\|_h^2
		&\lesssim \|u_h-J_hu_h\|_h^2 + Res(e-J_hI_hI_\M
		e)-(\hatF-\hatFapx)(v_h)-\hatFapx(v_h)\label{eqn:76_bound_a}\\
		&\lesssim \|u_h-J_hu_h\|_h^2 + \left(\mu(\T)+ \APX(F,\T)\right)\trb{u-J_hu_h}-\hatFapx(v_h)\label{eqn:76_bound}
	\end{align}
	with Lemma \ref{lem:F_osc} and Proposition \ref{prop:mu_rel} in the last step.
	The discussion in the proof of Theorem \ref{thm:WOPSIP_a_posteriori} implies \eqref{eqn:76_bound_a}--\eqref{eqn:76_bound} also for
	the 
	WOPSIP method. %
	The triangle inequality $\trb{u-J_hu_h}\leq \|u-u_h\|_h+\|u_h-J_hu_h\|_h$, \eqref{eqn:76_bound}, and a Young inequality 
verify
\begin{align}\label{eqn:proof_7.6}
		\|u-u_h\|_h^2
		&\lesssim \|u_h-J_hu_h\|_h^2 + \mu^2(\T)+ \APX^2(F,\T)-\hatFapx(v_h).
\end{align}
It remains to bound the extra term
$\hatFapx(v_h)$.
Recall the abbreviations $\tf_0,\tf_1, \tf_2$ from \eqref{eqn:F_abbr}.

The key step \W{towards an efficient control of $\hatFapx(v_h)$} is
	an integration by parts in \eqref{eqn:Fapx_intbp} that 
	collects the volume loads
	$F_0, F_1, F_2$ in the single residual term $\mu_1(\T)$ (resp.~the jumps in
	$\mu_2(\T), \mu_3(\T)$).
	A similar approach for the efficient bound of $\hatFapx(v_h)$ with the non-conforming test
	function $v_h\not\in V$ leads to additional terms from the product rule for jumps on the edge $E\in\E$, namely
	\begin{align*}
		\jump{F_1\cdot\nu_E v_h}&=\mean{F_1\cdot\nu_E}\jump{v_h}+\jump{F_1\cdot\nu_E}\mean{v_h},\\
		\jump{F_2\nu_E\cdot \nabla v_h}&=\mean{F_2\nu_E}\jump{\nabla v_h}+\jump{F_2\nu_E}\mean{\nabla v_h}.
	\end{align*}
	However, the
	average terms $\langle F_1\cdot \nu_E\rangle_E$
	and $\langle F_2\nu_E\rangle_E$
	over the edges $E\in\E$ are \emph{no} residuals and their efficiency is open; cf.~the partial efficiency
	result (excluding the average terms) in \cite[Thm.~7.2]{kim_morley_2021} or the omission of the efficiency analysis in
	\cite{CCGMNN18}.
	Instead, the assumptions \eqref{eqn:F_data_ass1}--\eqref{eqn:F_data_ass2} and the additional information on the
	structure of the test function $v_h\in (1-J_h)I_hI_\M V$ allows the efficient control of
	$\hatFapx(v_h)$.
\medskip

\noindent\emph{Case $I_h=\id$}:
Since $F_2\in P_{\hspace{-.13em}0}(\T;\mathbb S)=D^2_\pw P_{\hspace{-.13em}2}(\T)$ is piecewise constant, $I_\M v_h=0$ from $I_\M J_h=I_\M$ and \eqref{eqn:best_approx} verify the
$L^2$ orthogonality $v_h\perp F_2$.
This and \eqref{eqn:F_data_ass1} lead to
\begin{align}\label{eqn:F_h_id}
	\hatFapx(v_h)=(F_0, v_h)_{L^2(\Omega)}%
	+(G_0, 
	v_h)_{L^2(\Gamma_0)}+\sum^{}_{E\in\E(\Gamma_1)} ((1-\Pi_{E,0})G_1,
\partial_\nu v_h)_{L^2(E)}
\end{align}
with $\Pi_{E,0} \partial_\nu v_h\ds = 0$ for any $E\in\E$ from $I_\M v_h=0$.
\medskip

\noindent\emph{Case $I_h=I_C$}:
Since the test function $v_h\in V+S^2_0(\T)$ is $H^1$ conforming, \eqref{eqn:F_data_ass2} and an integration by parts
show
\begin{align}\label{eqn:F_h_IC}
	\hatFapx(v_h)=(F_0-\div_\pw F_1, v_h)_{L^2(\Omega)}+(G_0+[F_1]_E\cdot \nu_E, v_h)_{L^2(\Gamma_j)}.
\end{align}
Cauchy inequalities, Lemma \ref{lem:5.1}, and \eqref{eqn:F_h_id} for Morley, dG, WOPSIP and \eqref{eqn:F_h_IC} for
$C^0$IP result in $|\hatFapx(v_h)|\lesssim \mu(\T)\trb{e}$.
	This, \eqref{eqn:proof_7.6}, and a Young inequality provides
	\begin{align*}
			\|u-u_h\|_h^2
			&\lesssim \|u_h-J_hu_h\|_h^2 + \mu^2(\T)+ \APX^2(F,\T).
	\end{align*}
	Theorem \ref{thm:mu_bound} and the efficiency of $\mu(\T)$ from Proposition \ref{prop:mu_efficiency} conclude the proof.
\end{proof}
\begin{rem}[$\APX(F,\T)$ in Theorem {\ref{thm:Fh_a_post}}]\label{rem:apx_id}
	Since Theorem \ref{thm:Fh_a_post} requires $(F_\alpha)_{|\alpha|=1}\equiv0$ to vanish for the Morley, dG, and WOPSIP
	methods, the data approximation error $\APX(F,\T)$ includes the term
	$
		\sqrt{\sum^{}_{|\alpha|=1} \|h_\T f_\alpha\|_{L^2(\Omega)}^2}.
	$
	This term is linear in the mesh-size and converges with the expected rate for lowest-order schemes.
This term may even be of higher order if the triangulation is quasi-uniform and $\Omega$ is non-convex with a
reduced convergence rate $\|u-u_h\|_h=\mathcal{O}(h_{\rm max}^\sigma)$ of the schemes.
	However, it is \emph{not} a classical (higher-oder) data oscillation term if $f_\alpha\ne 0$ does not vanish
	for all $|\alpha|=1$.
	The assumption \eqref{eqn:F_data_ass2} for $C^0$IP leads to the term
	$
		\sqrt{\sum^{}_{|\alpha|=2} \|f_\alpha\|_{L^2(\Omega)}^2} %
	$
	independent of the mesh-size
	in the data approximation error $\APX(F, \T)$. Hence a meaningful interpretation of the a posteriori estimate in Theorem
	\ref{thm:Fh_a_post} for $C^0$IP \W{requires $\|f_\alpha\|_{L^2(\Omega)}$ small} for all $|\alpha|=2$.
\end{rem}
\begin{rem}[smoother vs.~no smoother]
	Since Theorem \ref{thm:F_a_post} for the smoother $Q=J_h$ applies to any choice of data approximations, 
	Remark \ref{rem:osc} shows that the data approximation error $\APX(F, \T)$ can be replaced by data oscillations
	$\OSC(F, \T)$ of arbitrary
	order.
	This provides a novel reliable and efficient
	a posteriori error bound for any right-hand side $F\in V^*$ of the form \eqref{eqn:hatF_repr} up to data
	oscillations.

	For no smoother $Q=\id$, additional requirements on the data approximations \eqref{eqn:F_data_ass1} for Morley, dG, and WOPSIP 
	(resp.~\eqref{eqn:F_data_ass2} for $C^0$IP) in Theorem \ref{thm:Fh_a_post} seem necessary for an
	efficient error control. 
	However, Remark \ref{rem:apx_id} explains that this either restricts the admissible data in \eqref{eqn:hatF_repr} or leads to 
	terms in the data approximation error $\APX(F, \T)$ that are no oscillations.

\end{rem}
\begin{rem}[$F\in L^2(\Omega)$]
	Theorems \ref{thm:F_a_post}--\ref{thm:Fh_a_post} for source terms $F\equiv f\in L^2(\Omega)$ with 
	$F_0\coloneqq \Pi_2 f\in P_{\hspace{-.13em}2}(\T)$ (and $F_\alpha=f_\alpha\equiv0$ for all $|\alpha|=1,2$ as well as $G_0=G_1\equiv0$)
	imply the a posteriori results
		 of Theorems \ref{thm:Morley_a_posteriori}, \ref{thm:dG_a_posteriori},
	\ref{thm:C0IP_a_posteriori}, and \ref{thm:WOPSIP_a_posteriori}. %
	Indeed, the Pythagoras theorem $\|f\|_{L^2(T)}^2=\|f-\Pi_2f\|_{L^2(T)}^2+\|\Pi_2 f\|_{L^2(T)}^2$ for the triangle
	$T\in\T$ verifies
	\begin{align*}
		\mu(\T)^2+\APX^2(F, \T) = \|h_\T^2f\|^2 + 
\sum^{}_{E\in\E(\Omega)}\begin{cases}{}
	0&\hspace{-2.1em}\text{if }I_h={\id},\\
h_E\|[D^2_\pw u_h\nu_E]_E\cdot \nu_E\|_{L^2(E)}^2&\text{else}.\end{cases}
	\end{align*}
	Since $j_h(u_h,u_h)=0$ vanishes for all Morley solutions $u_h\in \M(\T)$ and 
	every $C^0$IP solution $u_h\in S^2_0(\T)$ has zero jump $[u_h]_E\equiv 0$ along an edge $E\in\E$,
	Theorems \ref{thm:F_a_post}--\ref{thm:Fh_a_post} recover the corresponding results from Section \ref{sec:Unified a
	posteriori error control}.
\end{rem}

\subsection{Proof of Proposition \ref{prop:mu_efficiency}}
	\label{sub:Proof of mu_eff}
	This proof applies the bubble-function methodology \cite{Verfurth}. Recall
	\begin{align}\label{eqn:Du_h_T}
		\|D^2(u-I_\M u)\|_{L^2(T)}=\min_{v_h\in P_{\hspace{-.13em}2}(\T)}\|D^2(u-v_h)\|_{L^2(T)}\leq\|D^2(u-u_h)\|_{L^2(T)}
	\end{align}
	for any $T\in\T$ from the best-approximation property \eqref{eqn:best}.\medskip\\
	\emph{Step 1 (efficiency of the volume contribution)}. Let 
	$\varpi\coloneqq \tf_0-\div_\pw \tf_1+ \div_\pw^2\tf_2\in P_{\hspace{-.13em}k}(T)$
	abbreviate the volume contribution of $\mu(\T)$ for some $T\in\T$.
	The element bubble-function $b_T=27\varphi_1\varphi_2\varphi_3\in P_{\hspace{-.13em}3}(T)\cap H^1_0(T)$
	with $\|b_T\|_{L^\infty(T)}=1$ is given in terms of the three barycentric coordinates
	$\varphi_j\in P_{\hspace{-.13em}1}(T)$ for $j=1,2,3$.
	Since $a_\pw(I_\M u, b_T^2\varpi)=0$ from \eqref{eqn:best_approx} and $I_\M(b_T^2\varpi)\equiv0$, the 
	equivalence of the weighted norm $\|b_T\varpi\|_{L^2(T)}\approx\|\varpi\|_{L^2(T)}$ and an 
	integration by parts without boundary terms from $b_T^2\varpi\in H^2_0(T)$ show	
	\begin{align*}
		\|\varpi\|_{L^2(T)}^2
			&\approx(\tf_0-\div_\pw \tf_1+ \div_\pw^2 \tf_2, b_T^2\varpi)_{L^2(T)}\\
			&=\hatFapx(b_T^2\varpi)=a_\pw(u-I_\M u, b_T^2\varpi) + (\hatFapx(b_T^2\varpi) - F(b_T^2\varpi))\\
			&\lesssim \left(\|D^2(u-I_\M u)\|_{L^2(T)} + \APX(F, T)\right)\,\|D^2(b_T^2\varpi)\|_{L^2(T)}
	\end{align*}
	with $(F-\hatFapx)(v)=(F-\hatFapx)(v-J_hI_hI_\M v)\lesssim\APX(F, T)\|D^2 v\|_{L^2(T)}$ for $v\in H^2_0(T)$ from $I_\M v\equiv 0$ plus Lemma
	\ref{lem:F_osc} and a Cauchy inequality in the last step.
	This and the inverse inequality $h_T^2\|D^2(b_T^2\varpi)\|_{L^2(T)}\lesssim
	\|b_T^2\varpi\|_{L^2(T)}\leq\|\varpi\|_{L^2(T)}$ from \cite[Lem.~12.1]{ern_finite_2021-1} conclude the proof of the local efficiency
of the volume contributions, namely 
\renewcommand{\th}{\varrho}
\begin{align}\label{eqn:eff_T}
		h_T^2\|\tf_0-\div_\pw \tf_1+ \div_\pw^2 \tf_2\|_{L^2(T)}\lesssim \|D^2(u-I_\M u)\|_{L^2(T)} + \APX(F, T).
	\end{align}
	\emph{Step 2 (set-up for an interior edge)}.
	For any interior edge $E=\mathrm{conv}\{A,B\}=T_+\cap T_-\in\E(\Omega)$, let 
	$\varphi_{P_{\hspace{-.13em}+}},\varphi_{A,+},\varphi_{B,+}\in P_{\hspace{-.13em}1}(\mathbb R^2)$
	(resp.~$\varphi_{P_{\hspace{-.13em}-}},\varphi_{A,-},\varphi_{B,-}\in
	P_{\hspace{-.13em}1}(\mathbb R^2)$) denote the
	barycentric coordinates of $T_+=\mathrm{conv}\{P_{\hspace{-.13em}+},A,B\}$
	(resp.~$T_-\coloneqq\mathrm{conv}\{P_{\hspace{-.13em}-}, A, B\}$) seen as
	globally defined affine functions.
	The edge bubble-function reads
	$b_E\coloneqq
	16\varphi_{A,+}\varphi_{A,-}\varphi_{B,+}\varphi_{B,-}\in P_{\hspace{-.13em}4}(\mathbb R^2)\cap H^1_0(\omega(E))$ and $b_{T_\pm}\coloneqq
	27\varphi_{A,\pm}\varphi_{B,\pm}\varphi_{P_{\hspace{-.13em}\pm}}\in P_{\hspace{-.13em}3}(\mathbb R^2)\cap H^1_0(T_\pm)$ denotes the element bubble-function in
	$T_\pm\in\T$ with $\nu_{T_\pm}|_E=\pm\nu_E$.
	Let $\partial_\nu\coloneqq\partial/\partial{\nu_E}$ abbreviate the normal derivative and
	recall that the gradient
	$\nabla\varphi_{P_{\hspace{-.13em}+}}=-\th_E^{-1}\nu_E$ of the barycentric coordinate $\varphi_{P_{\hspace{-.13em}+}}$ scales like
	$h_E^{-1}$ with the height $\th_E\coloneqq 2|T_+|/|E|\approx h_E$ from shape-regularity.
	\W{The function $b_E^2$ has been utilised in the literature before, e.g., in 
\cite[p.~788]{BGS10} with its scaling properties; the usage of $b_{T_\pm}^2$ is standard.}
	The product rule and $\varphi_{T_+}|_E\equiv0$ verify \begin{align}\partial_\nu(\varphi_{T_+}b_E^2) =
	-\th_E^{-1}b_E^2\quad\text{on }E.\label{eqn:trace_vf_E}\end{align}
	Given $p\in \N_0$, any polynomial $q\in P_{\hspace{-.13em}p}(E)$ on the edge $E$ defines a unique polynomial on the straight line $L$ that extends
$E\subset L$.
The extension of $q$ from $L$ to $\mathbb R^2$ by constant values along the normal $\nu_E$ defines a polynomial
$\widehat q\in P_{\hspace{-.13em}p}(\mathbb R^2)$ on
$\mathbb R^2$ of the same degree.
Let $\Pi_LP_{\hspace{-.13em}\pm}\in L$ denote the projection of the vertex $P_{\hspace{-.13em}\pm}\in T_\pm$ opposite to $E$ onto $L$ along the normal direction $\nu_E$.
The maximal value $\|\widehat q\|_{L^\infty(\omega(E))}$ is attained on the line segment $\widehat L\coloneqq
\mathrm{conv}\{E,\Pi_LP_{\hspace{-.13em}\pm}\}\subset L$ and the shape-regularity controls the ratio $|\widehat L|/|E|\geq 1$. 
Hence
\begin{align}\label{eqn:L_inv_ext}
	\|\widehat q\|_{L^\infty(w(E))}=\|q\|_{L^\infty(\widehat L)}\leq C\|q\|_{L^\infty(E)}
\end{align}
follows 
with some constant $C\approx 1$ that exclusively depends on the shape-regularity of the triangulation $\T$ and on
$p$.\medskip\newline
\emph{Step 3 (efficiency of the first jump contribution)}. \W{This step establishes the local efficiency} of the term
	$\vartheta_E\coloneqq\tg_1+[\tf_1-\div_\pw \tf_2-\partial \tf_2\tau_E/\partial s ]_E\cdot\nu_E \in
	P_{\hspace{-.13em}k}(E)$ in the form
		\begin{align}\label{eqn:eff_F1}
		h_E^{3/2}\|\vartheta_E\|_{L^2(E)}\lesssim \|D^2_\pw(u-I_\M u)\|_{L^2(\omega(E))} + \APX(F, T_+)+ \APX(F, T_-).
	\end{align}
	Let $\xi_E\in P_{\hspace{-.13em}{2k}}(E)$ denote the (unique) Riesz representation of the functional
	$\th_E(\partial_\nu b_E^2, \bullet)_{L^2(E)}$ in the vector space $P_{\hspace{-.13em}{2k}}(E)$ with respect to the weighted scalar product
	$\left(b_E^2\bullet, \bullet\right)_{L^2(E)}$, i.e., 
	\begin{align}\label{eqn:xi_def}
		\big(b_E^2\xi_E, p_{2k}\big)_{L^2(E)}=\th_E\big(\partial_\nu b_E^2,
		p_{2k}\big)_{L^2(E)}&&\text{ for all }p_{2k}\in P_{\hspace{-.13em}{2k}}(E).
	\end{align}
	This, the equivalence of the weighted norm $\|b_E\xi_E\|_{L^2(E)}\approx\|b_E\|_{L^2(E)}$, and $h_E\approx \th_E$
	show
	\begin{align*}
		\|\xi_E\|_{L^2(E)}^2\approx \|b_E\xi_E\|_{L^2(E)}^2= h_E(\partial_\nu b_E^2, \xi_E)_{L^2(E)}\lesssim
		h_E\|\partial_\nu
		b_E^2\|_{L^2(E)}\|\xi_E\|_{L^2(E)}
	\end{align*}
	with a Cauchy inequality in the last step.
	Hölder's inequality and an inverse estimate \cite[Lem.\ 12.1]{ern_finite_2021-1} 
	lead to \begin{align*}\|\partial_\nu b_E^2\|_{L^2(E)}\leq
h_E^{1/2}\|\partial_\nu b_E^2\|_{L^\infty(E)}=h_E^{1/2}\|\nabla b_E^2\|_{L^\infty(E)}\lesssim
h_{E}^{-1/2}\|b_E^2\|_{L^\infty(T)}\leq h_E^{-1/2}.\end{align*}
This proves $\|\xi_E\|_{L^2(E)}\lesssim h_E^{1/2}$ and another inverse inequality provides $\|\xi_E\|_{L^\infty(E)}\lesssim h_E^{-1/2}\|\xi_E\|_{L^2(E)}\lesssim 1$.
Let $\hatg_E\in P_{\hspace{-.13em}k}(\mathbb R^2)$ and $\hatx_E\in P_{\hspace{-.13em}{2k}}(\mathbb R^2)$ denote the
extension of $\vartheta_E\in P_{\hspace{-.13em}k}(E)$
and $\xi_E\in P_{\hspace{-.13em}{2k}}(E)$ to $\mathbb R^2$ as in Step 2.
	This, \eqref{eqn:trace_vf_E}, and \eqref{eqn:xi_def} verify the $L^2$ orthogonality 
	\begin{align*}
		\partial_{\nu}((b_E^2+\varphi_{T_+}b_E^2\hatx_E)\hatg_E)=(\partial_\nu (b_E^2)
	-\th_E^{-1}b_E^2\hatx_E)\hatg_E\perp P_{\hspace{-.13em}k}(E)&&\text{ in }L^2(E).\end{align*} 
	Let $\xi_{T_\pm}\in P_{\hspace{-.13em}{2k}}(T_\pm)$ be the unique solution to
	\begin{align*}
		(b_{T_\pm}^2\xi_{T_\pm}, p_{2k})_{L^2(T_\pm)}=(b_E^2+\varphi_{T_+}b_E^2\hatx_E, p_{2k})_{L^2(T_\pm)}&&\text{for all }p_{2k}\in
		P_{\hspace{-.13em}2k}(T_\pm).
	\end{align*}
	An inverse inequality and \eqref{eqn:L_inv_ext} show $\|\xi_{T_\pm}\|_{L^\infty(T_\pm)}\lesssim 1$.
	The definition of $\xi_{T_\pm}$ verifies that the function $\psi_E\coloneqq
	(b_E^2+\varphi_{T_+}b_E^2\hatx_E-b_{T_+}^2\chi_{T_+}\hatx_{T_+}-b_{T_-}^2\chi_{T_-}\hatx_{T_-})\hatg_E\in
	H^2_0(\omega(E))$ is $L^2(T_\pm)$ orthogonal to $P_{\hspace{-.13em}k}(T_\pm)$.
	Since $b_{T_\pm}^2\in H^2_0(T_\pm)$ vanishes on $E$, the normal derivative $\partial_\nu\psi_E|_E\equiv
	\partial_\nu((b_E^2+\varphi_{T_+}b_E^2\hatx_E)\hatg_E)|_E\perp
	P_{\hspace{-.13em}k}(E)$ is $L^2(E)$ orthogonal to $P_{\hspace{-.13em}k}(E)$.
This, \eqref{eqn:S_ibp}, and an integration by parts show
	\begin{align*}
		0&=(\tf_0 - \div_\pw\tf_1+\div_\pw^2\tf_2, \psi_E)_{L^2(\omega(E))} + (\tg_2-[\tf_2\nu_E]_E\cdot\nu_E,
		\partial_\nu\psi_E)_{L^2(E)}\\
   		 &=(\tf_0, \psi_E)_{L^2(\omega(E))} + (\tf_1-\div_\pw\tf_2, \nabla\psi_E)_{L^2(\omega(E))} \\
   		 &\quad- ([F_1 - \div_\pw \tf_2]_E\cdot\nu_E, \psi_E)_{L^2(E)} + (\tg_2-[\tf_2\nu_E]_E\cdot \nu_E,
   		\partial_\nu\psi_E)_{L^2(E)} \\
		 &=\hatFapx(\psi_E) - (\vartheta_E, \psi_E)_{L^2(E)}.
	\end{align*}
	The Morley interpolation $I_\M\psi_E\equiv0$ of $\psi_E\in H^2_0(\omega(E))$ vanishes from $\Pi_{E,0}\partial_\nu\psi_E=0$ 
and $a_\pw(I_\M u, \psi_E) = 0$ follows from
	\eqref{eqn:best_approx}.
Since $(\psi_E-b_E^2\vartheta_E)|_E\equiv0$ is zero on $E$, the equivalence
	$\|\vartheta_E\|_{L^2(E)}\approx\|b_E\vartheta_E\|_{L^2(E)}$ results in
	$
		\|\vartheta_E\|_{L^2(E)}^2
		\approx (\vartheta_E, b_E^2\vartheta_E)_{L^2(E)}=(\vartheta_E,\psi_E)_{L^2(E)}=\hatFapx(\psi_E)$.
	With $a(u, \psi_E) = F(\psi_E)$ from \eqref{eqn:WP}, this shows
	\begin{align*}
		\|\vartheta_E\|_{L^2(E)}^2&=a_\pw(u-I_\M u, \psi_E) + (\hatFapx(\psi_E) - F(\psi_E))\\
		&\leq  \left(\|D^2_\pw(u-I_\M u)\|_{L^2(\omega(E))} + \APX(F, T_-)+ \APX(F,
		T_+)\right)\,\|D^2\psi_E\|_{L^2(\omega(E))}
	\end{align*}
	The inverse inequality, $\|\psi_E\|_{L^2(\omega(E))}\lesssim \|\hatg_E\|_{L^2(\omega(E))}$, and
	\eqref{eqn:L_inv_ext} provide
\begin{align}\label{eqn:F1_inverse}
	h_E^2\|D^2\psi_E\|_{L^2\omega(E))}\lesssim \|\psi_E\|_{L^2(\omega(E))}\lesssim
\|\hatg_E\|_{L^2(\omega(E))}\lesssim h_E^{1/2}\|\vartheta_E\|_{L^2(E)}.\end{align}
This verifies the efficiency \eqref{eqn:eff_F1} of the jump contributions $\vartheta_E$.\qed\medskip\newline
\emph{Step 4 (efficiency of the second jump contribution)}.
	The local efficiency of the remaining term follows with similar arguments.
	Since the Hessian $D^2_\pw u_h$ of $u_h\in P_{\hspace{-.13em}2}(\T)$ is piecewise constant,
	the stability of the $L^2$ projection results in
	\begin{align*}
		\|(1-\Pi_{E,0})(\tg_1+[\tf_2\nu_E]_E\cdot \nu_E)\|_{L^2(E)}\leq\|\tg_1+[(\tf_2-D^2_\pw u_h)\nu_E]_E\cdot \nu_E\|_{L^2(E)}.
	\end{align*}
It is therefore sufficient to prove the local efficiency of the term 
$\zeta_E\coloneqq\tg_1+[(\tf_2-D^2_\pw u_h)\nu_E]_E\cdot\nu_E\in P_{\hspace{-.13em}k}(E)$, namely
\begin{align}\label{eqn:eff_F2}
		h_E^{1/2}\|\zeta_E\|_{L^2(E)}\lesssim \|D^2_\pw(u-u_h)\|_{L^2(\omega(E))} + \APX(F, T_+)+ \APX(F, T_-).
	\end{align}
	Indeed, let $\varrho_{T_\pm}\in P_{\hspace{-.13em}{2k}}(T_\pm)$ be the unique solution to
	\begin{align*}
		(b_{T_\pm}^2\varrho_{T_\pm}, p_{2k})_{L^2(T_\pm)}=(\varphi_{T_+}b_E^2, p_{2k})_{L^2(T_\pm)}&&\text{for all }p_{2k}\in
		P_{\hspace{-.13em}{2k}}(T_\pm).
	\end{align*}
	Observe that $\psi_2\coloneqq -(\varphi_{T_+}b_E^2-b_{T_+}^2\chi_{T_+}\varrho_{T_+}-
	b_{T_-}^2\chi_{T_-}\varrho_{T_-})\widehat\zeta_E\in H^2_0(\omega(E))$
	is $L^2$ perpendicular %
	to $P_{\hspace{-.13em}k}(T_\pm)$ with zero trace $\psi_2|_E\equiv 0$ on $E$.
	This and an integration by parts show 
	\begin{align*}
		&0=(\tf_0 - \div_\pw \tf_1+\div_\pw^2\tf_2, \psi_2)_{L^2(\omega(E))}+ (\vartheta_E, \psi_2)_{L^2(E)}\\
		 &=%
		 \hatFapx(\psi_2)-(\tg_1+[\tf_2\nu_E]_E\cdot \nu_E, \partial_\nu\psi_2)_{L^2(E)}=\hatFapx(\psi_2)-a_\pw(u_h,\psi_2)
		 -(\zeta_E,\partial_\nu\psi_2)_{L^2(E)}
	\end{align*}
	with $([D^2_\pw u_h\nu_E]_E\cdot \nu_E, \partial_\nu\psi_2)_{L^2(E)}=a_\pw(u_h, \psi_2)$ in the last step.
	The equivalences
	$\|\zeta_E\|_{L^2(E)}\approx \|b_E\zeta_E\|_{L^2(E)}$ and $h_E\approx \widetilde h_E$, \eqref{eqn:trace_vf_E}, and
	$\partial_\nu(\psi_2+\varphi_{T_+}b_E^2\widehat\zeta_E)|_E\equiv 0$ provide
	\begin{align*}
		h_E^{-1}\|\zeta_E\|_{L^2(E)}^2
		&\approx \th_{E}^{-1}(\zeta_E,b_E^2\zeta_E)_{L^2(E)}=- (\zeta_E,
		\partial_\nu(\varphi_{T_+}b_E^2\widehat\zeta_E))_{L^2(E)}=(\zeta_E, \partial_\nu\psi_2)_{L^2(E)}\\
		&=\hatFapx(\psi_2)-a_\pw(u_h, \psi_2)=a_\pw(u-u_h,\psi_2) + (\hatFapx(\psi_2)-F(\psi_2))
	\end{align*}
	with $a(u,\psi_2)=F(\psi_2)$ from \eqref{eqn:WP} in the last step.
	The remaining steps follow Step 2.1 and utilize $\|D^2\psi_2\|_{L^2(\omega(E))}\lesssim h_E^{-3/2}\|\zeta_E\|_{L^2(E)}$ from an inverse
	inequality as in \eqref{eqn:F1_inverse}; further details are omitted.
	The combination of the local efficiency results \eqref{eqn:eff_T}, \eqref{eqn:eff_F1}, and \eqref{eqn:eff_F2} with
	\eqref{eqn:Du_h_T} concludes the proof.%
	\qed

\section*{Acknowledgements}
The research of the first \W{two authors have} been supported by the Deutsche Forschungsgemeinschaft in the Priority Program 1748
under the project \emph{foundation and application of generalized mixed FEM towards nonlinear problems
in solid mechanics} (CA 151/22-2).  This paper has been supported by   SPARC project 
(id 235) {\it the mathematics and computation of plates} and SERB POWER Fellowship SPF/2020/000019.
\W{The second author is supported by the \emph{Berlin Mathematical School, Germany}.}

{\footnotesize
\bibliographystyle{amsplain}
\bibliography{ref_combined,Numerik}}
	\W{
		\newpage
		\appendix
		\section[A posteriori error control of piecewise polynomial sources]{A posteriori error control of a piecewise
		polynomial source in $H^{-2}(\Omega)$}
	\label{apx:A posteriori error control of a piecewise polynomial source}}
This appendix provides an alternative view on the reliable and efficient estimator from Section \ref{sec:General
sources} as lower and upper bounds for the dual norm of a piecewise polynomial source in $H^{-2}(\Omega)$.
Suppose the piecewise polynomials $\Lambda_0\in P_k(\T), \Lambda_1\in P_k(\T;\mathbb R^2)$, and $\Lambda_2\in
P_k(\T;\mathbb{S})$ define the linear functional $\Lambda\in H^{-2}(\Omega)$ by
\begin{align}\label{eqn:F_repr_apx}
	\Lambda(v)\coloneqq \int_\Omega (\Lambda_0 v + \Lambda_1\cdot \nabla v + \Lambda_2:D^2 v)\;\mathrm{d} x\qquad \text{for all
	}v\in H^2_0(\Omega).
\end{align}
Recall the transfer operators $I_\M,I_h,J_h$ for the five quadratic discretization schemes of Section \ref{sec:Examples of second-order
finite element schemes} listed in Table \ref{tab:spaces}.
A reliable and efficient estimator $\mu^2(\T) \coloneqq \mu_1^2(\T) + \mu_2^2(\T) + \mu_3^2(\T)$ of the functional
$\Lambda$ is given by
\begin{align*}
	\mu_1^2(\T)&\coloneqq\|h_\T^2(\Lambda_0 - \div_\pw \Lambda_1+\div_\pw^2\Lambda_2)\|^2, \\
	\mu_2^2(\T)&\coloneqq \sum^{}_{E\in\E(\Omega)}
	h_E^{3}\|[\Lambda_1-\div_\pw \Lambda_2-\partial (\Lambda_2\tau_E)/\partial s]_E\cdot\nu_E\|_{L^2(E)}^2,\\
\mu_3^2(\T)&\coloneqq \sum^{}_{E\in\E(\Omega)}
\begin{cases}{}
	h_E\|(1-\Pi_{E,0})[\Lambda_2\nu_E]_E\cdot \nu_E\|_{L^2(E)}^2& \text{if }I_h=\id,\\
	h_E\|[\Lambda_2\nu_E]_E\cdot \nu_E\|_{L^2(E)}^2&\text{if }I_h=I_\C.
\end{cases}
\end{align*}
\begin{theorem}[reliability and efficiency]\label{thm:rel_eff_apx}
	There exist positive constants $C_{\rm{rel}}, C_{\rm{eff}}>0$ that exclusively depend on the shape regularity of
	$\T$ and on the polynomial degree $k\in\N_0$ such that
	\begin{align*}
		C_{\rm rel}^{-1}\trb{\Lambda\circ(1-J_hI_hI_\M)}_*\leq\mu(\T)\leq C_{\rm eff}\trb{\Lambda}_*.
	\end{align*}
\end{theorem}
\begin{proof}
	The discussion in Subsection \ref{sub:Estimator for the residual} applies to $Res\coloneqq F\coloneqq\Lambda$
	and $u_h\coloneqq0$ with $\APX(F,
	\T)=0$.
	In this particular case, Proposition \ref{prop:mu_rel} provides the first inequality
\begin{align*}
	\trb{\Lambda\circ(1-J_hI_hI_\M)}_*\leq C_{\rm rel}\mu(\T)
\end{align*}
with $C_{\rm rel}=\const{cst:5}$. 
	Let $u\in H^2_0(\Omega)$ denote the Riesz representation of $a(u,\bullet)=\Lambda\in H^{-2}(\Omega)$ with the isometry
	$\trb{\Lambda}_*=\trb{u}$ in the Hilbert space $(H^2_0(\Omega), a)$ and $\trb{\bullet}\equiv
	a(\bullet,\bullet)^{1/2}$.
	Then 
	the efficiency estimate
\begin{align*}
	\mu(\T)\leq C_{\rm eff}\trb{\Lambda}_*
\end{align*}
follows from Proposition \ref{prop:mu_efficiency} with $C_{\rm eff}=\const{cst:mu_efficiency}$.
\end{proof}
Theorem \ref{thm:rel_eff_apx} allows for a direct
application to the linearization of semilinear problems in \cite{carstensen_semilinear202+}.
It can be further generalized in various directions, e.g., in the spirit of Section~\ref{sec:General
sources} that considers the a posteriori error analysis of the linear biharmonic
problem for a more general class of functionals in $H^{-2}(\Omega)$ including
line and point loads.
The reliability requires only piecewise smoothness of $\Lambda_0,\Lambda_1,\Lambda_2$ so that the traces and derivatives
in $\mu_1,\mu_2,\mu_3$ exist, while the efficiency may require extra oscillation terms (as in
\eqref{eqn:mu_res_bound}).

\end{document}